\documentclass [english,12pt,leqno]{article}
\usepackage[latin1]{inputenc}

\usepackage{fullpage}

\usepackage{tipa}
\usepackage{dsfont}
\usepackage{tikz-cd}
\usepackage{float}
\usepackage{graphicx}

\usepackage{amsxtra}
\usepackage{amsmath}
\usepackage{amssymb}
\usepackage{amsfonts}
\usepackage[all,2cell]{xy}
\usepackage{mathrsfs}
\usepackage{amsthm}
\usepackage{enumitem}
\usepackage{hyperref}
\usepackage{subfigure}
\usepackage{mathabx}
\usepackage{euscript}
\usepackage[bbgreekl]{mathbbol}
\usepackage{pgfplots}

\pgfplotsset{width=7cm, compat=newest}

\usetikzlibrary{patterns}
\usetikzlibrary{calc}

\usepackage[OT2,T1]{fontenc}
\DeclareSymbolFont{cyrletters}{OT2}{wncyr}{m}{n}
\DeclareMathSymbol{\Sha}{\mathalpha}{cyrletters}{"58}

\makeatletter
\hypersetup{
unicode=true,           
pdftoolbar=true,        
pdfmenubar=true,        
pdffitwindow=false,     
pdfstartview={FitH},    
pdftitle={\@title},     
pdfauthor={\@author},   
pdfsubject={},          
pdfcreator={},          
pdfproducer={},         
pdfkeywords={},         
pdfnewwindow=true,      
colorlinks,             
linkcolor=black,        
citecolor=black,        
filecolor=black,        
urlcolor=black          
}
\makeatother

\newtheorem{thm}[equation]{Theorem}

\newtheorem{lem}[equation]{Lemma}
\newtheorem{prop}[equation]{Proposition}

\newtheoremstyle{example}{\topsep}{\topsep}%
 {}
 {}
 {\bfseries}
 {.}
 {2pt}
 {\thmname{#1}\thmnumber{ #2}\thmnote{ #3}}

\theoremstyle{example}

\newtheorem{Defi}[equation]{Definition}
\newtheorem{rem}[equation]{Remark}
\newtheorem{rems}[equation]{Remarks}

\newtheorem{exas}[equation]{Examples}
\newtheorem{ex}[equation]{Example}

\newtheorem{asses}[equation]{Assumptions}

\numberwithin{equation}{section}

\newcommand{\vertiii}[1]{{\left\vert\kern-0.25ex\left\vert\kern-0.25ex\left\vert #1 
    \right\vert\kern-0.25ex\right\vert\kern-0.25ex\right\vert}} 

\def\CC{\mathbb{C}}
\def\DD{\mathbb{D}}
\def\EE{\mathbb{E}}
\def\FF{\mathbb{F}}

\def\PP{\mathbb{P}}
\def\RR{\mathbb{R}}
\def\ZZ{\mathbb{Z}}
\def\QQ{\mathbb{Q}}
\def\TT{{\mathbb{T}}}

\def\aen{\mathfrak{a}}

\def\hen{\mathfrak{h}}

\def\Ben{\mathfrak{B}}

\def\Ac{\mathcal{A}}

\def\Cc{\mathcal{C}}

\def\Dc{\mathcal{D}}
\def\Ec{\mathcal{E}}
\def\Fc{\mathcal{F}}
\def\Gc{\mathcal{G}}
\def\Ic{{\mathcal{I}}}

\def\Hc{\mathcal{H}}

\def\Pc{\mathcal{P}}

\def\Sc{\mathcal{S}}
\def\Tc{\mathcal{T}}

\def\Vc{\mathcal{V}}

\def\ab{\mathbf{a}}

\def\Ano{{\on{Ano}}}

\def\b1{{\mathbf{1}}}

\def\ba{{{\bf a}}}
\def\bb{{\bf {b}}}
\def\bd{{\bf {d}}}
\def\bI{{{\bf I}}}

\def\bm{{{\bf m}}}
\def\bn{{{\bf n}}}

\def\bz{{{\bf z}}}

\def\Cat{{\on{Cat}}}

\def\CM{{\on{CM}}}
\def\CMen{{\mathfrak {CM}}}
\def\codim{{\on{codim}}}

\def\Cone{{\on{Cone}}}

\def\cu{{\on{cu}}}

\def\db{{\mathbf{d}}}
\def\Dec{{\on{Dec}}}
\def\del{{\partial}}

\def\eps{{\varepsilon}}
\def\Exit{{\on{Exit}}}
\def\Ext{{\on{Ext}}}

\def\Fun{{\on{Fun}}}

\def\gr{{\on{gr}}}
\def\Gro{{\on{Gro}}}

\def\Hom{\operatorname{Hom}\nolimits}
\def\hra{\hookrightarrow}

\def\Id{\operatorname{Id}\nolimits}

\def\Im{{\on{Im}}} 

\def\JPS{{\on{JPS}}}
\def\JS{{\on{JS}}}

\def\k{\mathbf k}
\def\Ker{\operatorname{Ker}\nolimits}

\def\lla{\longleftarrow}
\def\lra{\longrightarrow}

\def\Map{{\on{Map}}}
\def\Mar{{\on{Mar}}}

\def\Mor{{\on{Mor}}}

\def\Ob{\operatorname{Ob}\nolimits}
\def\on{\operatorname}
\def\ol{\overline}
\def\oo{{\infty}}

\def\op{{\on{op}}}
\def\OR{{\on{OR}}}

\def\Perv{{\on{Perv}}}
\def\phi{{\varphi}}

\def\Ran{{\on{Ran}}}

\def\Set{{\on{Set}}}

\def\Sh{{\on{Sh}}}

\def\Sup{{\on{Sup}}}
\def\supp{{\on{supp}}}
\def\Supp{{\on{Supp}}}
\def\STop{{\on{STop}}}

\def\Sym{{\on{Sym}}}

\def\Top{{\Tc op}}

\def\Tot{{\on{Tot}}}

\def\ul{\underline}
\def\un{{\on{un}}}

\def\vbn{{\overset{\rightarrow}{\bn}}}
\def\Vect{\on{Vect}}

\def\wt{\widetilde}

\def\+{{\oplus}}
\def\- {{\,\setminus\,}}
\def\={{\,\simeq\,}}
\def\x{{\otimes}}
\def\1baa{{\mathbf{1}}}
\def\2{{\mathbf{2}}}
\def\(({(\hskip -1mm (}
\def\)){)\hskip -1mm )}
\def\be{\begin{equation}}
\def\ee{\end{equation}}
\def\ed{\end{document}}

\def\0bar{{\ol {0}}}
\def \1bar{{\ol {1}}}

\title{  PROBs  and perverse sheaves II. Ran spaces and $0$-cycles with coefficients}
\author{  Mikhail Kapranov, Vadim Schechtman}

\begin{document}

\maketitle

 \thanks{\em To Ezra Getzler for his sixtieth birthday}

 \begin{abstract}
 We consider the space $Z(\CC,L)$ of $0$-cycles on the complex line $\CC$ with coefficients in a commutative
 monoid $L$ subject to certain conditions. Such spaces include the symmetric products (for $L=\ZZ_+$)
 and the Ran space (for $L=\TT=\{ {\tt{True}}, {\tt{False}} \}$ being the Boolean algebra of truth values). We
 describe the appropriately defined category of perverse sheaves on $Z(\CC,L)$ in terms of the braided
 category (PROB) generated by the components of the universal $L$-graded
  bialgebra. We give another description in terms of so-called Janus sheaves which are objects
  of mixed functoriality (data covariant in one direction and contravariant in the other) on a category
  formed by certain matrices with entries in $L$. The matrices in question are analogs of contingency tables
  familiar in statistics.

   \end{abstract}

\tableofcontents

\addtocounter{section}{-1}

\vskip .6cm

\section{Introduction}

\paragraph{Description of the main result.}

This is the second part of the series, see \cite{KS-prob}, where we develop a new
type of ``quiver description'' of the categories of perverse sheaves on certain configuration
spaces. The ``quiver'' in our description is a block (direct summand) of the braided
monoidal category (colored PROB) governing graded bialgebras of appropriate type.
We refer to \cite{KS-prob} and references therein for discussion of the terms PROP and PROB. 

\vskip .2cm

More precisely, we start with a commutative monoid $(L,+,0)$ with operation written additively
and  unit element $0$, satisfying natural monotonicity conditions (Assumptions \ref{ass:L})
and consider the space $Z(\CC,L)$ formed by $0$-cycles (or divisors)
$\bz=\sum_{z\in\CC} n_z\cdot z$ on the complex line
$\CC$ with coefficients  $n_z\in L$. This space splits into connected components $Z(\CC,L)_n$
labelled by $n\in L$ (the degree of the $0$-cycle), and we consider each $Z(\CC,L)_n$
separately. Among the $Z(\CC,L)_n$ we find:
\begin{itemize}

\item The symmetric products $\Sym^n(\CC) = Z(\CC,\ZZ_+)_n$, corresponding to $L=\ZZ_+$. This is the case considered in  \cite{KS-prob}. 

\item The Ran space $\Ran(\CC)=Z(\CC,\TT)_\1bar$. It corresponds to $L=\TT=\{\0bar, \1bar\}$
being the Boolean algebra of truth values $\0bar= {\tt{False}}$, $\1bar =  {\tt{True}}$,
with the operation $+ =\vee =  {\tt{or}}$ being the Boolean addition. 

\item Various flag Ran spaces and spaces of colored divisors, see \S\ref{sec:0-cyc}\ref {par:0-cyc}
\end{itemize}

\noindent The space $Z(\CC,L)_n$ may be infinite-dimensional but it has a natural structure of
an inductive limit of stratified spaces with strata being complex manifolds. So fixing a base field 
$\k$,
 we can consider perverse sheaves (of $\k$-vector spaces, with middle perversity) on
 $Z(\CC,L)_n$ supported on a finite-dimensional skeleton, similarly to the approach common, e.g.,
 in the theory of the Geometric Satake Correspondence \cite{mirkovich-vilonen}. We denote this category
 $\Perv^{<\oo}(Z(\CC,L)_n)$. 
 
 \vskip .2cm 
 
 On the other hand, for any $\k$-linear braided monoidal category $(\Vc, \otimes, \b1, R)$
 we can consider $L$-graded bialgebras $A=\bigoplus_{n\in L} A_n$ in $\Vc$, with $A_0=\b1$
 and all (co)multiplications involving $A_0$ being the identities. Note that we can avoid speaking
 of direct sums or assuming their existence in $\Vc$, formulating everything in terms of
 data involving the individual $A_i$. With this understanding, there is a {\em universal
 $L$-graded bialgebra} $\aen$ belonging to a similar universal braided category $\Ben^L$
 which is generated (as a braided category) by formal generating objects $\ba_n$, 
 $n\in L$ and by generating morphisms which are the components  of the multiplication and
 comultiplication, subject only to the axioms of an $L$-graded bialgebra, see \S 
 \ref{sec:bi-PROB}\ref{par:bi-PROB} We refer to $\Ben^L$ as the
 {\em  (colored) PROB of $L$-graded bialgebras}. Our main result, Theorem 
 \ref{thm:Bn-perv}, says that $\Perv^{<\oo}(Z(\CC,L)_n)$ is identified with the category of
 functors $\Ben_n^L\to\Vect_\k$ which vanish on all but finitely many objects. In fact,
 we consider a nominally more general case where $\Vect_\k$ is replaced by an
 arbitrary abelian category.
 
 \paragraph{The main tool: contingency matrices.}
 Our main tool relating the topology of $Z(\CC,L)$ and $L$-graded bialgebras is
 the concept of a contingency matrix with values in $L$. By this we mean, generalizing
 \cite{KS-cont} , a matrix $M=\|m_{ij}\|$ over $L$ of arbitrary size $p\times q$
 without zero rows or columns. Remarkably, they appear in two contexts;
 \begin{itemize}
 \item They parametrize cells $U_M$ of a natural {\em contingency cell decomposition}
 $\Sc^{(2)}$ of $Z(\CC,L)$. More precisely, we we record the real and imaginary parts of
 every complex number separately, then every finite subset of $\CC$ can be seen
 as lying naturally on a finite grid. The multiplicities $n_z\in L$ of any $0$-cycle
 $\bz=\sum_{z\in \CC} n_z\cdot z$ arrange therefore into a contingency matrix $M(\bz)$,
 see Fig. \ref{fig:S-a-Im}
 and \S \ref{sect:cont}\ref{par:S2} The cell $U_M$ consists of $\bz$ with $M(\bz)=M$. 
 
 \item They appear naturally (and classically for $L=\ZZ_+$) in formulas expressing
 the compatibility between multiplication and comultiplication in graded bialgebras as the level
 of graded components, see \S \ref{sec:bi-PROB}\ref{par:L-gr-bi}
 \end{itemize}
 
 \noindent 
In general, the cell decomposition $\Sc^{(2)}$ is not quasi-regular, i.e., informally, a bigger
 cell can approach a smaller cell from several different directions. This type of behavior
 is formalized by Mac Pherson's concept of the exit path category
  \cite{treumann, curry, curry-patel}.
 We describe such category for $\Sc^{(2)}$ in terms of the two augmented semi-simplicial
 structures on the set of all contingency matrices, given by adding rows and by adding columns.
 See Proposition \ref{prop:exit-S-2=gro} (where we consider $0$-cycles in arbitrary $\RR^p$, not just in
 $\CC=\RR^2$). This allows us to describe the category of $\Sc^{(2)}$-constructible
 sheaves as representations of the exit path category. 
 
 Note that perverse sheaves on $Z(\CC,L)_n$ are, in particular, $\Sc^{(2)}$-constructible 
 complexes since the natural ``complex stratification by multiplicities'' $\Sc^{(0)}$
 is refined by $\Sc^{(2)}$. 
 
 \paragraph{Second main result: Janus sheaves.} 
 Our proof of Theorem
\ref{thm:Bn-perv} is based on another result, Theorem 
\ref{thm:perv=cont} which is interesting in its own right
and analogous to the main result of \cite{KS-hW}. It
describes the category
 $\Perv^{<\oo}(Z(\CC,L)_n, \Vc)$
in terms of combinatorial diagrams of bivariant nature (covariant in one direction,
contravariant in the other)
which we call {\em Janus sheaves}.
Such diagrams are analogous to mixed Bruhat sheaves introduced in \cite{KS-hW}
in the context of Coxeter groups.  A nontrivial intersection
between the two concepts appears for $L=\ZZ_+$, when Janus sheaves
are the same as mixed Bruhat sheaves for the symmetric group $S_n$ acting on $\RR^n$. 
In this case we have a {\em  third appearance} of contingency matrices: as parametrizing
the parabolic Bruhat decompositions  for the space of partial classical ($GL_n$-) flags,
see  \cite{chriss-ginzburg} \S 4.3. 

\vskip .2cm

If we view a perverse sheaf $\Fc\in \Perv^{<\oo}(Z(\CC,L)_n)$ as a complex of
$\Sc^{(2)}$-constructible sheaves  so as a complex $F^\bullet$ of functors
on the exit path category, then the Janus sheaf corresponding to $\Fc$  can be obtained
from $F^\bullet$ by an appropriate formal version of  {\em relative Verdier duality}.
Here by ``relative'' we mean duality taken along the fibers of the imaginary part map
$\Im: Z(\CC,L) \to \ZZ(\RR,L)$ to the space of $0$-cycles on $\RR$. A remarkable
fact is that the perverse t-structure on constructible complexes is taken by such
duality into the standard t-structure on the derived category, so the Janus sheaf
associated to a perverse sheaf $\Fc$ is a non-derived object, i.e., consists of
vector spaces and not complexes. 

\vskip .2cm

More generally, it seems to us that the two main points of view on perverse sheaves used
in algebraic geometry: as complexes of sheaves and as ``linear algebra data''
with values in an abelian category, are almost invariably connected by an appropriate
relative version of Verdier duality. We plan to develop this point of view in a future paper.

\paragraph{Ran spaces and logic?} Specializing to Ran spaces,
 we have a cell decomposition $\Sc^{(2)}$ of $\Ran(\CC)$ labelled by {\em Boolean 
 contingency matrices} $M=\| m_{ij }\|$, $m_{ij}\in\TT=\{\0bar, \1bar\}$. The cell
 attachment structure (i.e., the category of exit paths)  of this cell decomposition
 is described in terms of Boolean addition of rows and columns of such matrices. This
 adds an extra touch to the ``logical flavor'' of Ran spaces in general. Thus, for any
 finite connected CW-complex $X$ the space $\Ran(X)$ formed by finite nonempty
 subsets $I\subset X$, possesses an idempotent semigroup structure $(I,J)\mapsto I\cup J$.
 This implies that $\Ran(X)$ is weakly contractible \cite[\S 3.4.1]{BD}, a fact whose algebro-geometric
 analog (for $X$ an algebraic curve) plays an important role in the study of the moduli
 stack of $G$-bundles on $X$, see \cite{gaitsgory-lurie}. 
 
 \vskip .2cm
 
 It is therefore tempting to ``marry'' Ran spaces with more general algebraic structures
 appearing in mathematical logic. For example, considering $0$-cycles with coefficients
 in the ``algebra of truth values of the $N$-valued logic'', see Example \ref{exas:L}(f), leads
 to the flag Ran space formed by flags of finite subsets
  $X\supset I^1\supset\cdots\supset I^{N-1}$. Further, 
 one can likely weaken Assumption \ref{ass:L}(b) so as to allow
more general  distributive lattices of interest in logic such as, e.g., 
 the lattice of vector subspaces in  $\FF_q^N$  with addition
 being the sum of subspaces (``quantum logic'').
  The corresponding  versions of Ran spaces seem worth studying. 
  
  \paragraph{Further directions: non-finitistic perverse sheaves and factorization.}
  It is desirable (and we believe possible) to develop a more general concept of perverse
  sheaves on ind-stratified spaces such as $Z(\CC,L)$ (in particular, on ind-schemes
  such as the affine Grassmannian). We mean ``perverse sheaves'' not necessarily supported
  on a finite skeleton and not represented as limits of perverse subsheaves thus supported.
  One reason for attempting this is to incorporate factorization structures, i.e., to speak about
  {\em factorizing systems} $(\Fc_n)_{n\in L}$ of perverse sheaves on all the
  $Z(\CC,L)_n$ at once. Such systems should correspond to actual $L$-graded
  bialgebras $A$, i.e., to braided monoidal functors from the PROB $\Ben^L=\bigoplus_{n\in L}
  \Ben^L_n$. Note that for a general $A$ (with all $A_n\neq 0$) each $\Fc_n$
  must be supported on the entire $Z(\CC,L)_n$ and so cannot be supported on a finite
  skeleton unless $Z(\CC,L)_n$ is finite-dimensional (as in the case $L=\ZZ_+$,
  $Z(\CC,\ZZ_+)_n =\Sym^n(\CC)$ treated in \cite{KS-shuffle, KS-prob}). 
  
  Note that the naive approach of considering all constructible complexes on $Z(\CC,L)_n$
  and imposing the standard geometric perversity conditions  
  does not seem to work as it is not clear whether this gives a t-structure. 
  Still we believe that there is a natural and satisfactory way to address this issue and hope to take
  it up in the future. 
  
  \paragraph{Organization of the paper.} In \S \ref{sec:0-cyc} we recall the basic properties of  the space $Z(X,L)$ of
  $0$-cycles  in a space $X$ with coefficients in an abelian monoid $L$. Taken by itself, this generality is not new: much
  more general ``configuration spaces with summable labels'' were studied by P. Salvatore  \cite{salvatore}
  and we refer to his work for foundational material. We focus, however, on monoids satisfying certain
  monotonicity properties (Assumptions \ref{ass:L}) which allow us to treat $Z(X,L)$ for $X$ a manifold as an
  (ind-)stratified space and to speak about perverse sheaves on it later on. These assumotions also allow
  us to identify, for connected $X$,  the
  connected components $Z(X,L)_n$ with  elements $n\in L$ (the degree of  $0$-cycles).

  \S \ref{sect:cont} is devoted to the concept of contingency matrices with coefficients in a monoid $L$. We consider
  ``$p$-dimensional contingency matrices'' (i.e., arrays) $\|m_{i_1,\cdots, i_p}\|$ for any $p$. 
  The set of such matrices with {\em content} (sum of all the elements) equal to a given $n\in L$, is 
    organized   into a
  $p$-fold augmented semi-simplicial set $\ul\CM^p_n(L)$ by adding together adjacent slices in each of the $p$ directions.

  In \S \ref{sec:cont-cells} we construct a cell decomposition $\Sc^{(2)}$ of $Z(\RR^p, L)_n$ labelled by $p$-dimensional
  contingency matrices  of content $n$ and describe the exit path category of this cell decomposition as the Grothendieck
  construction of $\ul\CM^p_n(L)$ (Proposition \ref{prop:exit-S-2=gro}).

  In \S \ref{sec:X=C} we specialize to the case $X=\CC=\RR^2$ and consider some related stratifications of
  $Z(\CC,L)$: the Fox-Neuwirth-Fuchs cell decomposition $\Sc^{(1)}$ and the decomposition $\Ic$ into
  imaginary strata $Z_\bn^\Im$.

  In \S \ref{sec:bi-PROB} we define the concept of $L$-graded bialgebras in braided monoidal categories and
  consider the PROB $\Ben^L$ generated by the components of the universal such bialgebra. In the case
  $L=\TT$, relevant for Ran spaces, we relate $\TT$-graded bialgebras with more familiar ungraded bialgebras
  with unit and counit.

  In \S \ref{sec:perv-0-cyc} we define the category $\Perv^{<\oo}(Z(\CC,L)_n)$ and formulate the main result
  (Theorem  \ref{thm:Bn-perv}) relating it with the block $\Ben^L_n$ of $\Ben^L$. We also prove some
  technical results describing various features of  and operations on constructible complexes in terms of their combinatorial description via the category of exit paths.

  In \S \ref{sec:janus} we introduce the concept of Janus sheaves and give an example of a Janus sheaf associated to
  an $L$-graded bialgebra (Proposition \ref{prop:bialg-Janus}). We also formulate our second main result,
  Theorem \ref{thm:perv=cont} relating Janus sheaves with perverse sheaves.

  \S  \ref{sec:janus=cousin} is devoted to the proof of  Theorem \ref{thm:perv=cont}. Our main main tool, similar to
  \cite{KS-hW}, is the concept of a Cousin-type complexes which on one hand, are algebraically in a one-to-one
  correspondence with Janus (pre)sheaves and on the other, appear as Cousin resolutions of perverse sheaves.

  Finally, in \S \ref{sec:uni-janus} we deduce Theorem \ref{thm:Bn-perv} from Theorem  \ref{thm:perv=cont}.
  This is done by considering the {\em universal Janus sheaf} taking values in a braided category $\CMen(L)$
  specially constructed for this purpose. We then identify $\Ben^L$ and $\CMen(L)$ as braided categories
  (Theorem \ref{thm:Janus=CM}) by constructing a natural $L$-graded bialgebra in $\CMen(L)$.

  \paragraph{Acknowledgements.} 
  
  We are grateful to  D. Gaitsgory for correspondence about the issues related to defining perverse sheaves on the Ran
  space.

  This work was supported by the  World Premier International Research Center Initiative (WPI Initiative), 
 MEXT, Japan.  The work of M.K. was also supported by the JSPS  KAKENHI grant 20H01794.


\section{$0$-cycles with coefficients in a monoid}\label{sec:0-cyc}

\paragraph{Spaces of $0$-cycles.}\label{par:0-cyc}

Let $(L,+)$ be a commutative monoid (semigroup with unit element, denoted $0$). Let $X$ be a topological space
which we assume homeomorphic to an open subset of a finite CW-complex.  By a $0$-{\em cycle} on $X$ with coefficients
in $L$ we will mean a formal linear combination $\bz = \sum_{x\in X} n_x\cdot x$  of points $x\in X$ 
with multiplicities $n_x\in L$ such that $n_x=0$ for almost
all $x$.   Denote by $Z(X,L)$ the set of all $0$-cycles on $X$
with coefficients in $L$. This set has a natural topology induced by the topology of $X$ and the additive structure on $L$
(``when points merge together, their multiplicities are added''), see \cite{salvatore}.  
  With respect to this topology and the obvious addition of $0$-cycles, 
$Z(X,L)$ is a  commutative topological monoid with unit element $0$, the $0$-cycle with all $n_x=0$. 

\begin{exas}\label{exas:L}
(a) Suppose that $L$ is an abelian group. In this case the classical Dold-Thom theorem says that
$
  \pi_i Z(X,L) \= H_i(X,L)
$
for all $i\geq 0$. In particular, taking $X=S^n$ to be the $n$-dimensional sphere, we get $Z(S^n,L) \= L\times K(L,n)$,
where $K(L,n)$ is the Eilenberg-MacLane space. 

\vskip .2cm

(b) If $X$ is connected, then for any monoid $L$ as above we have $\pi_0 Z(X,L)\= L$. The identification is given by
associng to a $0$-cycle $\bz=\sum n_x\cdot x$ its {\em degree} $\deg(\bz)=\sum n_x \in L$. That is, the component
$Z(X,L)_n$ corresponding to $n\in L$ consists of cycles of degree $n$. 

\vskip .2cm

(c) Let $L=\ZZ_+$. Then $Z(X,\ZZ_+)=\bigsqcup_{n\geq 0} \Sym^n(X)$ is the disjoint union of the symmetric powers of $X$.
Thus $Z(X,\ZZ_+)_n=\Sym^n(X)$. 

\vskip .2cm

(d) More generally, let $L=\Lambda_+$ be the dominant (co)weight lattice of a reductive group. Then $L\= \ZZ_+^I$ where
$I$ is the set of fundamental (co)weights. If $X$ is an algebraic curve over $\CC$, then $Z(X,\Lambda_+)$
is the space of $I$-colored divisors on $X$  considered in  \cite{finkelberg-mirkovic} \cite {gaitsgory-fact}. 

\vskip .2cm

(e) Let $L=\TT = \bigl\{ \0bar = {\tt{False}}, \1bar = {\tt{True}}\bigr\} $ be the Boolean algebra of truth values, with 
$+=\vee$ being the Boolean addition: $\0bar+x=x$, $\1bar+x=\1bar$. Suppose $X$ connected. Then $Z(X,\TT)$ 
consists of two components: $Z(X,\TT)_\0bar = \{0\}$ and $Z(X,\TT)_\1bar = \Ran(X)$ being the {\em Ran space}
of $X$, i.e., the space of all nonempty finite subsets $I\subset X$ equipped with natural (Vietoris) topology, see \cite[\S 3.4.1]{BD}.

\vskip .2cm

(f) More generally, let $N\geq 2$ be an integer and 
\[
L=\TT_N = \bigl\{0, 1/(N-1), 2/(N-1), \cdots, 1\bigr\}
\]
 be the
``algebra of truth values of  $N$-valued logic'', with the operation $x+y=\max(x,y)$
(maximum with respect to the order induced from $\RR$). Then $Z(X, \TT_N)$ is the  {\em flag Ran space} of $X$, i.e., the space formed by flags of finite subsets
$X\supset I^1\supset\cdots\supset I^{N-1}$ with the topology induced from $(\{0\}\sqcup \Ran(X))^{N-1}$, in which
it is embedded as a closed subset given by the incidence conditions. Explicitly, let us put $I^0=X$, $I^N=\emptyset$. Then
a flag $I^1\supset\cdots\supset I^{N-1}$ corresponds to the $0$-cycle $\sum n_x\cdot x$ where
$n_x=k/(N-1)$, if $x\in I^k\- I^{k+1}$. 
\end{exas}

\vfill\eject

\paragraph{ Skeleta of $Z(X,L)$.} 
From now on we impose the following.

\begin{asses}\label{ass:L}
\begin{itemize}
\item[(a)] $L$ is finitely generated as a monoid. 

\item[(b)] $L$ admits a total order $\leq$ compatible with the addition such that $0$ is the minimal
element and each lower interval $L_{\leq n}= \{m\in L \,| m\leq n\}$ is finite. 
\end{itemize}
\end{asses}

\begin{exas}
Assumptions \ref{ass:L} exclude part (a) of Examples \ref{exas:L} but allow  parts (c)-(f). Indeed:

\vskip .2cm

 (a) An abelian group other than
$\{0\}$ does not admit a compatible order with $0$ being minimal.

\vskip .2cm

(c)(d)
 The free monoid $\ZZ_+^d$ admits a compatible order as in Assumption \ref{ass:L} (b).
Indeed, let us choose positive real numbers $\alpha_1,\cdots, \alpha_d$ which are
linearly independent over $\QQ$. Then the map
\[
\ZZ_+^d \lra \RR, \quad (a_1,\cdots, a_d) \mapsto \sum a_i \alpha_i 
\]
is an embedding. Pulling back the order on $\RR$ under this embedding, we get an order
on $\ZZ_+^d$ satisfying the assumption.

\vskip .2cm

(e)(f) The monoid  $\TT_N$ (of which $\TT=\TT_2$ is a particular case) is finite and comes equipped with
the order induced from $\RR$ which satisfies the assumption. 
\end{exas}

\begin{prop}\label{prop:L-finite-part}
 For any $d\geq 1$ and $n\in L$ the set
  \[
 \{(n_1, \cdots, n_d)\in L^d \,| n_1+\cdots + n_d=n\}
 \]
 is finite. 
 \end{prop}

\noindent{\sl Proof:} The equality $n_1+\cdots + n_d=n$ implies that each $n_i\leq n$, so
our set is embedded into the Cartesian power $L_{\leq n}^d$ which is finite. 
 \qed
 
 \vskip .2cm
 
 We now study the space $Z(X,L)$ in more detail. First, note that
 our assumptions imply that
 $0$ is an isolated point of $Z(X,L)$.

 \vskip .2cm
 
 Second, let us choose once and for all a system of generators $(\xi_k)_{k\in K}$ of $L$
 so we have a surjective morphism of monoids $\ZZ_+^K\to L$. Let $X_K=K\times X$
 be the disjoint union of $K$ copies of $X$; its points are pairs $(k,x)$, $k\in K$, $x\in X$. 
 For each $r\geq 0$  we have the
  continuous map
 \be\label{eq:map-p-l}
 p_l: X_K^l \to Z(X,L), \quad ((k_1, x_1), \cdots, (k_l,  x_l)) \mapsto \sum\,  \xi_{k_i} \cdot x_i. 
 \ee
 
 \begin{prop}\label{prop:p_n-proper}
 Each $p_l$ is a proper map with finite fibers. 
 \end{prop}
 
 \noindent{\sl Proof:} The fact that $p_r$ has finite fibers follows from Proposition 
 \ref{prop:L-finite-part}.
 The fact that $p_r$ is proper follows from this and from the definition of topology on
 $Z(X,L)$. \qed
 
 \vskip .2cm
 
 Proposition \ref{prop:p_n-proper} implies that the image of $p_l$ is a closed subset of $Z(X,L)$. 
 We denote by $Z_l(X,L)$ the union of the images of the $p_i, i\leq l$. Thus we have an
 exhaustion
 \be\label{eq:exhaustion1}
 Z_0(X,L)\subset Z_1(X,L) \subset \cdots \subset Z(X,L)\,=\, \varinjlim\nolimits^\Top_l \, Z_l(X,L)
 \ee
 by closed subsets so that the topology on $Z(X,L)$ is that of their inductive limit.
 We will refer to  $Z_l(X,L)$ as the {\em $l$th skeleton} of $Z(X,L)$. 
 
 
 \paragraph{Strata of $Z(X,L)$ labelled by partitions.}\label{par:strata-part}
 
  From now on we further assume that $X$ is a $C^\oo$-manifold 
   of dimension $p$ (without boundary). 
  
  \vskip .2cm
 
 Let $n\in L$ be a nonzero element. A {\em partition} of $n$ is, by definition, a sequence $\alpha=
 (\alpha_1,\cdots, \alpha_d)$ of nonzero elements such that $\sum n_i=n$, considered up to
 permutations.  
 We denote $\Pc^L_n$ the set of all partitions of $n$. We put 
 $\Pc^L_0=\{\emptyset\}$
 (the set consisting of the empty partition)  and $\Pc^L=\bigsqcup_{n\geq 0} \Pc^L_n$. 
 We introduce a partial order $\leq$ on $\Pc^L$ putting
 \[
 \beta=(\beta_1,\cdots, \beta_c) \leq \alpha=(\alpha_1,\cdots, \alpha_d)
 \]
 if $\beta$ can be obtained from $\alpha$ by partial summation, i.e., there is a surjection
 $f: \{1,\cdots, d\}\to \{1,\cdots, c\}$ such that
  $\beta_i = \sum_{j\in f^{-1}(i)} \alpha_j$. 
 
 \vskip .2cm
 
 Let $\bz=\sum_{x\in X} n_x\cdot x\in Z(X,L)$. Writing the set $\{x\, |\,  n_x\neq 0\}$ as 
 $\{x_1,\cdots, x_d\}$
 in some order, we get a partition $\alpha_\bz = (n_{x_1}, \cdots, n_{x_d})$ of $n=\deg(\bz)$. 
 For a given partition $\alpha\in\Pc^L$ 
 we denote $Z(X,L)_\alpha\subset Z(X,L)$ the set
 of $\bz$ such that $\alpha_\bz=\alpha$. Thus we obtain a disjoint decomposition
 \be
 Z(X,L) = \bigsqcup _{\alpha\in\Pc^L} Z(X,L)_\alpha
 \ee
 into locally closed subsets. We denote this decomposition $\Sc^{(0)}$
 and refer to the $Z(X,L)_\alpha$ as the {\em strata} of $\Sc^{(0)}$.

 \begin{prop}\label{prop:start-S0}
 (a) Each stratum  $Z(X,L)_\alpha$, $\alpha=(\alpha_1,\cdots, \alpha_d)$ is a $C^\oo$-manifold of dimension $pd$. 
 
 \vskip .2cm
 
 (b) The closure $\ol{Z(X,L)_\alpha}$ is the union of the $Z(X,L)_\beta$
 for $\beta\leq\alpha$. 
 
 \vskip .2cm
 
 (c) With respect to the decomposition $\Sc^{(0)}$, each  skeleton $Z_l(X,L)$, has a natural structure of a Whitney stratified space
 and each  $Z_{l}(X,L)\hra Z_{l'}(X,L)$, $r\leq r'$ is 
 a closed embedding of stratified spaces. 
 
  \end{prop}
  
   Thus $(Z(X,L), \Sc^{(0)})$ has the structure of an {\em ind-stratified space},
 a filtered inductive limit of stratified spaces and their closed stratified 
 embeddings. 
 
 \vskip .2cm
  
  \noindent{\sl Proof:} (a) Let us write  
  us write $\alpha\in \Pc^L_n$ as consisting of $d_1$ repetitions of $n_1$, $d_2$ repetitions
  of $n_2$, ..., $d_q$ repetitions of $n_q$, where $n_1,\cdots, n_q$
  are distinct. So $d_1n_1+\cdots + d_q n_q=n$. Then it is immediate that
  $Z(X,L)_\bn$ is  identified with the  open subset in $\Sym^{d_1}(X)\times \cdots\times  \Sym^{d_q}(X)$ formed by tuples
  \[
  \bigl( \{x_1,\cdots x_{d_1}\}, \{x_{d_1+1},\cdots, x_{d_1+\cdots + d_2}\}, \cdots,
  \{x_{d_1+\cdots + d_{q-1}+1},\cdots, x_{d_1+\cdots + d_n}\}\bigr)
  \]
 such that all the $x_i$ are distinct. 
 
 \vskip .2cm
 
 (b) This is clear from the definition of the topology on $Z(X,L)$, 
 see \S \ref{sec:0-cyc}\ref{par:0-cyc}
 
 \vskip .2cm
 
 (c)  We consider the Cartesian power $X_K^l$  with its natural stratification
 by multiplicities (patterns of coincidence of points from $X_K$).
 The map $p_l$ of Proposition \ref{prop:p_n-proper} takes this stratification
 to the decomposition of $\Im(p_l)$ induced by $\Sc^{(0)}$.
 Since $p_l$ is a finite map, the identifications on $X_K^l$ induced
 by the $p_l$ define the stratified structure on $\Im(p_l)$.
 Since $Z_l(X,L)$ is the union of the $\Im(p_i)$, $i\leq l$,
 we conclude the same for $Z_l(X,L)$. \qed
 
 \begin{ex}\label{ex:partitions-T}
 Let $L=\TT = \{\0bar, \1bar\}$, so $Z(X,\TT)_\1bar = \Ran(X)$. 
  The set $\Pc^\TT_\1bar$
 consists  of partitions $\1bar^d = (\1bar, \cdots, \1bar)$  ($d$ times) of arbitrary length $d\geq 1$. 
 The corresponding stratum $\Ran(X)_{\1bar^d}$ is identified with
 $\Sym^d_\neq (X)$, the open part in $\Sym^d(X)$ consisting of
 collections of $d$ distinct points.  Further, 
 the natural set of generators
 of $\TT$ consists of one element $\1bar$. This gives the skeleta
 $\Ran_r(X) = \{ I\subset X: |I|\leq r\}$. Each such skeleton is a stratified space,
 and the ind-stratified structure on $\Ran(X)$ is given by
 \[
 X= \Ran_1(X) \subset \Ran_2(X)\subset \cdots \subset \Ran(X),
 \quad \Ran_r(X)\- \Ran_{r-1}(X) \=  \Sym^r_\neq (X). 
 \]

 \end{ex}


 \section{Contingency matrices with values in a monoid}\label{sect:cont}
 
 \paragraph{$p$-dimensional contingency matrices.}\label{par:cont}
 Let $L$ be a commutative monoid satisfying Assumptions
 \ref{ass:L} and let $p,  d_1,\cdots, d_p\geq 1$ be integers.
 By a {\em $p$-dimensional   matrix} of size
 $d_1\times\cdots\times d_p$ with values in $L$ we will mean
 simply an array
 \[
 M=\| m_{i_1,\cdots, i_p}\|, \quad m_{i_1,\cdots, i_p}\in L, \,\,\, 1\leq i_\nu
 \leq d_\nu
 \]
 of elements of $L$. Thus for $p=1$ we have vectors $\bm=(m_1,\cdots, m_d)$,
 for $p=2$ matrices $M=\|m_{ij}\|$ in the usual sense and so on.
 The {\em content} of a $p$-dimensional matrix $M$ is defined as
 $|M|=\sum m_{i_1,\cdots, i_p}\in L$. 
 
 \vskip .2cm
 
 Let $1\leq \nu\leq p$ and $1\leq i\leq d_\nu$. The {\em $i$th slice}\
 of $M$ in the direction $\nu$ is the $(p-1)$-dimenional matrix
 $\sigma_i^\nu(M)$ formed by elements $m_{i_1, \cdots, i_{\nu-1}, i,
 i_{\nu+1}, \cdots, i_p}$ with the $\nu$th index fixed and equal to $p$. 
 Thus for $p=2$ we obtain the rows and columns of a usual matrix.
 For $p=1$ we define the $i$th slice $\sigma_i^1(M)$ to be the $i$th
 element of the vector $M=(m_1,\cdots, m_d)$. 
 
 \vskip .2cm

 We say that $M$ is a {\em contingency matrix}, if no slice $\sigma_i^\nu(M)$
 consists entirely of zeroes. Thus a $1$-dimensional contingency matrix
 is a vector $\bm=(m_1,\cdots, m_d)$ with all $m_i\neq 0$
 (we will also refer to such vectors as {\em ordered partitions}).
 A $2$-dimensional contingency matrix is a matrix $M=\|m_{ij}\|$
 with all the rows and columns being nonzero. 
 We denote by $\CM_n(d_1,\cdots, d_p; L)$ the set of contingency matrices
 of size $d_1\times\cdots\times d_p$ with values in $L$ which have content
 $n$.  Proposition \ref{prop:L-finite-part} implies that each $\CM_n(d_1,\cdots, d_p; L)$ is a finite set. 
  We  put
 \[
 \CM_n^p(L) = \bigsqcup_{d_1,\cdots, d_p\geq 1} \CM_n(d_1,\cdots, d_p; L).
 \]
 Here $n\neq 0$.
 We also introduce formally the {\em empty $p$-dimensional contingency matrix}
$\emptyset$ of content $0$ and  put
 $  \CM^p_0(L)=\{\emptyset\}$. We then put
 \[
 \quad \CM^p(L) = \bigsqcup_{n\in L} \CM_n^p(L). 
 \]
   
  \vskip .2cm
  
  Let $M\in \CM_n(d_1,\cdots, d_p; L)$. For $1\leq\nu\leq p$ and $0\leq i\leq
  d_{\nu-1}$ we define the contingency matrix
   \[
  \del_i^\nu M\in\CM_n^{p}(d_1,\cdots, d_{\nu-1}, d_\nu-1, d_{\nu+1},\cdots,
  d_p; L)
\]
  by putting
  \be\label{eq:def-del}
  (\del_i^\nu M)_{i_1,\cdots, i_p} =
  \begin{cases}
  m_{i_1,\cdots, i_p}, & \text{if } i_\nu\leq i,
    \\
  m_{i_1,\cdots, i_{nu-1}, i+1, i_{\nu+2}, \cdots, i_p} +  
   m_{i_1,\cdots, i_{nu-1}, i+2, i_{\nu+2}, \cdots, i_p}, & \text{if } i_\nu =i+1,
   \\
   m_{i_1,\cdots, i_{\nu-1}, i_\nu -1, i_{\nu+1}, \cdots, i_p}, &
   \text{if } i_\nu \geq i+2, 
  \end{cases}
  \ee
  i.e., by adding the $(i+1)$st and $(i+2)$nd slices in direction $\nu$. 
  The fact that $\del_i^\nu$ is again a contingency matrix, follows from
  Assumptions   \ref{ass:L}. 
  
  We refer to the operations $\del_i^\nu$ as {\em contractions}. 
  It is clear that the  $\del_i^\nu$ commute for different $\nu$
  and satisfy the (semi-)simplicial relations for the same $\nu$:
  \be\label{eq:del-simpl}
  \del_i^\nu\del_j^\mu = \del_j^\mu \del_i^\nu, \,\, \mu\neq\nu,
  \quad
  \del_i^\nu\del_j^\nu = \del_{j-1}^\nu \del_i^\nu, \,\, i<j. 
  \ee
  This extends the definition given in
   \cite{KS-cont} \S 2 for the case $p=1,2$ and $L=\ZZ_+$. 
   
   \vskip.2cm
   
   To formulate this more concisely, let $\Delta_+$ be the {\em augmented semi-simplex category}
   whose objects are all finite ordinals, including $\emptyset$ and morphisms are monotone
   embeddings. As usual, we denote the ordinal $\{0,1,\cdots, m-1\}$, $m\geq 0$, by $[m]$
   and $\emptyset$ by $[-1]$. Eq. \eqref{eq:del-simpl} can be restated as:
   
   \begin{prop}
   The maps $\del_i^\nu$ make $\CM_n^p(L)$ into a $p$-fold augmented semi-simplicial set.
   More precisely, they define a functor
   \[
   \ul\CM_n^p(L) = \CM_n^p(\bullet+2, \cdots, \bullet+2; L): (\Delta_+)^\op\lra\Set. \qed
   \]
   \end{prop}
   
   For future convenience introduce the following.
   
   \begin{Defi}\label{def:margins}
   Let $M\in \CM_n^p(d_1,\cdots, d_p; L)$. For $\nu=1,\cdots, p$ the {\em $\nu$th margin}
   of $M$ is the ordered partition $\Mar^\nu(M) \in \CM_n^1(d_\nu; L)$ with components
   \[
   \Mar^\nu(M)_j  = \sum m_{i_1, \cdots, i_{\nu-1}, j, i_{\nu+1}, \cdots, i_p},
   \]
   the summation being over all collections of 
   $i_\lambda \in \{1,\cdots, d_\lambda\}$, $\lambda\neq \nu$. 
   \end{Defi}
   
   \paragraph{Contingency matrices as objects of a category: the Grothendieck construction.}
   \label{par:gro}
   Let $\Cc$ be a category and $F: \Cc \to\Set$ be a functor. By the {\em Grothendieck
   construction} of $F$ we will mean the category
   $\Gro(F)$ whose objects are pairs $(c,a)$ with $c\in\Ob(\Cc)$ and $a\in F(c)$. 
   A morphism $(c,a)\to(c', a')$ in $\Gro(F)$ consists of a morphim $\phi: c\to c'$ in $\Cc$ such that $F(\phi)(a)=a'$. Composition of morphisms in $\Gro(F)$ is induced by the composition of
   morphisms in $\Cc$. This definition is a particular case of the standard more general
   concept defined for (pseuso-)functors with values in categories, not just sets (discrete
   categores), see, e.g., [...].

   Alternatively, for $c\in\Ob(\Cc)$ let $h^c=\Hom_\Cc(c,-): \Cc \to\Set$ be the convariant
   functor represented by $c$. Then elements of $F(c)$ are in bijections with natural
   transformations $h^c\to F$. Thus $\Gro(F)$ consists of all natural transformations
   of this kind. 
   
      \vskip .2cm
      
      We will apply this to $\Cc = (\Delta_+^p)^\op$ and $F=X:  (\Delta_+^p)^\op\to\Set$
      being a $p$-fold augmented semi-simplicial set. Objects of $\Delta_+^p$ being tuples
      $[\db]=([d_1], \cdots, [d_p])$, $\db=(d_1,\cdots, d_p)$, $d_i\geq -1$, we denote by
      \[
      \Delta^p_+[\db] = \Hom_{\Delta^p_+}(-, [\db]): (\Delta_+^p)^\op\lra\Set
      \]
   the functor represented by $[\db]$ and call it the {\em standard multisimplex} of dimension
   $\db$. Thus $\Gro(X)$ can be viewed as consisting of all morphisms $\Delta^p_+[\db]\to X$,
   i.e., of  ``all  multisimplices of $X$''. 
   
   We further specialize to $X=\ul\CM_n^p(L)$ being the $p$-fold augmented semi-simplicial
   set of $p$-dimensional contingency matrices of content $n$. Thus 
   $\Ob(\Gro(\ul\CM_n^p(L)))=\CM_n^p(L)$ is the set of contingency matrices and morphisms
   are given by contractions and their compositions. That is, $\Mor(\Gro(\ul\CM_n^p(L)))$
   is generated by the morphisms $M\to \del_i^\nu M$ for $M\in \CM_n(d_1,\cdots, d_p)$,
   $\nu=1,\cdots, p$ and $i=0,\cdots, d_p-2$; we will call these  {\em elementary morphisms}. 
   
   For $\nu=1,.\cdots, p$ we denote by $\Gro^\nu(\ul\CM_n^p(L))$ the subcategory in
   $\Gro(\ul\CM_n^p(L))$ with the same objects and morphisms induced only by the $\nu$th copy
   of $\Delta_+^\op$, i.e., by elementary morphisms $M\to \del_i^\nu M$ with this fixed $\nu$.
   It is clear that any morphism in $\Gro(\ul\CM_n^p(L))$ can be uniquely written as the
   composition of $p$ morphisms, one from each of the $\Gro^1(\ul\CM_n^p(L))$, ...,
   $\Gro^p(\ul\CM_n^p(L))$ (in each of the $p^!$ possible orders). 
   
   \begin{ex}
    Note that $\emptyset$ is an initial object of $\Delta_+$. Therefore $(n)$, the $1$-entry
    matrix of size $1\times\cdots\times 1$, is a final object in $\Gro(\ul\CM_n^p(L))$:
    each $M\in\CM_n^p(L)$ has a canonical morphism $M\to (n)$ given by the total contraction
    (summation of all the entries).
   \end{ex}
   
   \begin{ex}\label{ex:S-2-sym}
   Let $L=\ZZ_+$. In this case the category $\Gro(\ul\CM_n^p(\ZZ_+))$ reduces to a poset, i.e.,
   for any two $M,M'\in \CM_n^p(\ZZ_+)$ there is at most one morphism $M\to M'$ in
   $\Gro(\ul\CM_n^p(\ZZ_+))$.
   
   This fact is immediately seen for $p=1$, when elements of $\CM_n^1(\ZZ_+)$, i.e.,
   unordered partitions of $n\in \ZZ_+$, are in bijection with faces of the $(n-1)$-dimensional
   simplicial cone $\RR_+^{n-1}$, and the contraction $\del_i=\del_i^1$ correspond to taking
   faces of the faces of this cone. The case $p>1$ is reduced to this. 
   Indeed, given a contraction $M\to M'$ of $p$-dimensional contingency matrices we have,
   for each $\nu=1,\cdots, p$, the contraction $\Mar^\nu(M)\to \Mar^\nu(M')$ of unordered
   partitions obtained by taking the margins in the $\nu$th direction (Definition 
   \ref{def:margins}). By the above, such a contraction corresponds to a unique morphism
   of $\Delta_+$. The $p$-tuple of such morphisms is a morphism of $\Delta_+^p$ which is
   the only possible morphism that can correspond to $M\to M'$. 
   \end{ex}
   
   \begin{ex}
   Let $L=\TT$. In this case already for $p=1$ the category $\Gro(\ul\CM_\1bar^1(\TT))$ is
   not reduced to a poset. In fact, this category is identified with $\Delta_+^\op$, with the
   ordered partitions $(\1bar)$, $(\1bar, \1bar)$, $(\1bar, \1bar, \1bar)$ etc. corresponding
   to the ordinals $[-1]=\emptyset$, $[0]$, $[1]$ etc. 
   \end{ex}
   
   \paragraph{Anodyne morphisms of contingency matrices.} \label{par:ano}
   An elementary morphism
   $M\to\del_i^\nu M=M'$ in $\Gro(\ul\CM_n^p(L))$ will be called {\em anodyne} if, in
    each sum in the second line of 
   \eqref{eq:def-del}, at least one summand is equal to $0$. In this
   case the sets of nonzero entries of $M'$ and $M$ are in bijection. 
  More generally, a morphism $M\to M'$ will be called
  anodyne, if it can be decomposed into a chain of elementary
   anodyne morphisms. In this case again, the sets of nonzero entries
  of $M'$ and $M$ are in bijection. 
  
  We denote by $\Ano$ resp. $\Ano^\nu$ the set of anodyne morphisms in 
   $\Gro(\ul\CM_n^p(L))$ resp.  $\Gro^\nu(\ul\CM_n^p(L))$. We further denote
   $\equiv_\Ano$ resp. $\equiv^\nu_\Ano$ the equivalence relation of $\CM_n^p(L)$
   generated by the binary relation $R_\Ano$ resp. $R^\nu_\Ano$, where
   $M \, R_\Ano N$ resp. $M\, R^\nu_\Ano N$ means that there is an anodyne morphis,
   $M\to N$ in $\Gro(\ul\CM_n^p(L))$ resp.  $\Gro^\nu(\ul\CM_n^p(L))$. 
   Thus, contingency matrices from each equivalence class of $\equiv_\Ano$ have the same set of
   nonzero entries.


 \section{The contingency cell decomposition of the space of $0$-cycles in $\RR^p$}\label{sec:cont-cells}
 
 \paragraph{The contingency cell decomposition $\Sc^{(2)}$.}\label{par:S2}
 
  We now consider the case $X=\RR^p$ being
 the Euclidean space. Let $\bz=\sum_{x\in\RR^p} n_x\cdot x$ be a point of 
 $Z(\RR^p,L)$. We associate to $\bz$ a $p$-dimensional contingency matrix
 $M(\bz)=\|m_{i_1,\cdots, i_p}(\bz)\|$ as follows.  Each $x\in \RR^p$ is
 a vector $(x_1, \cdots, x_p)$. For $1\leq\nu\leq p$ let $I_\nu  
 \subset \RR$ be the set of the $\nu$th coordinates $x_\nu$ of all the
 $x\in\RR^p$ such that $n_x\neq 0$. Let $d_\nu = |I_\nu|$ and let us write
 $I_\nu = \{ t^{(\nu)}_1 < \cdots < t^{(\nu)}_{d_\nu}\}$.  By definition,
 $n_x\neq 0$ implies that $x\in I_1\times\cdots\times I_p$. 
  We define the entries of $M(\bz)$ to be the multiplicities of the 
  points of $ I_1\times\cdots\times I_p$ in $\bz$:
  \[
  m_{i_1,\cdots, i_p} \,=\, n_{(t^{(1)}_{i_1}, \cdots, t^{(p)}_{i_p})} \,\in \, L. 
  \]
  By construction, $M(\bz)$ is a contingency matrix. 
  Figure \ref{fig:S-a-Im} which we reproduce with notational changes
   from \cite {KS-cont}, illustrates this construction for $p=2$.

 \begin{figure}[H]
 \centering
 \begin{tikzpicture}[scale=.4, baseline=(current bounding box.center)]

 \draw[->] (-2,0) -- (13,0); 
 
 \draw[->] (0,-2) -- (0,11); 
 
 \draw[dashed] (0,2) -- (11,2); 
 
  \draw[dashed] (0,4) -- (11,4); 
  
   \draw[dashed] (0,7) -- (11,7); 
   
      \node at (-1,2){$t^{(2)}_1$};   \node at (-1,4){$t^{(2)}_2$};  
      \node at (-1,5.5){$\vdots$}; 
       \node at (-1,7){$t^{(2)}_{d_1}$}; 
   
   \node at (2,7){$\bullet$}; 
   \node at (2,2){$\bullet$}; 
  \node at (8,2){$\bullet$}; 
    \node at (5,4){$\bullet$};  
    \node at (8,7){$\bullet$};  \node at (8,4){$\bullet$};

    \draw[dashed] (2,0) -- (2,9);  
        \draw[dashed] (5,0) -- (5,9);  
            \draw[dashed] (8,0) -- (8,9);  
    
    \node at (2,-1){$t^{(1)}_1$};   \node at (5,-1){$t^{(1)} _2$};   \node at (8,-1){$t^{(1)}_{d_1}$}; 
    \node at (6.5, -1){$\cdots$}; 
    
    \node at (2.9, 2.5){ \small{$m_{1,1}$}}; 
      \node at (2.9, 7.5){\small {$m_{1, d_2}$}}; 
       \node at (5.9,4.5){\small {$m_{2,2}$}};    
        \node at (8.9,2.5){\small {$m_{d_1, 1}$}};   \node at (8.9,4.5){\small {$m_{d_1, 2}$}};  
         \node at (8.9,7.5){\small {$m_{d_1, d_2}$}};  

 \end{tikzpicture}
 \caption{The contingency matrix $M(\bz)$ associated to $\bf z\in\Sym^n(\RR^2)$.}\label{fig:S-a-Im}
 \end{figure}

  \vskip .2cm
  
  Let now $M\in \CM^p(L)$. Define $U_M = U_{M}^ {\RR^p}\subset Z(\RR^p, L)$ to be the
  set of all $0$-cycles $\ba$ with $M(\bz)=M$. 
  
  \begin{prop}\label{prop:S2}
  (a) $U_M$ is a topological cell,  i.e., it is homeomorphic to a Euclidean space.
  If $M\in \CM_n(d_1,\cdots, d_p; L)$, then $\dim U_M = d_1+\cdots + d_p$.

  \vskip .2cm
  
  (b) Let $\bn(M)$ be the partition of $|M|$ into the nonzero entries of $M$.
  Then $U_M$ lies in the stratum $Z(X,L)_{\bn(M)}$ of the stratification
   $\Sc^{(0)}$, see Proposition \ref{prop:start-S0}. 
  
  \end{prop}

 \noindent{\sl Proof:}  
 
 To see (a), note that $\bz\in U_M$ is completely determined by the
 data of subsets $I_\nu\subset \RR$, $\nu=1,\cdots, p$, with $|I_\nu|=d_\nu$.
 Therefore
 \[
 U_M\, \=\,  \RR^{d_1}_< \times\cdots\times \RR^{d_p}_<, 
 \text{ where } \RR^d_< = \{(t_1, \cdots, t_d)\in\RR^d |\, t_1 <\cdots < t_d\}
 \]
  is indeed a cell. Part (b) is clear \qed

  \vskip .2cm

 We will  denote the cell decomposition of $Z(\RR^p, L)$ into the $U_M$
 by $\Sc^{(2)}=\Sc^{(2)}_{\RR^p}$ and call it the
 {\em contingency cell decomposition}. The cells $U_M$ will be called
 {\em contingency cells}. 
 
  \paragraph{$\Sc^{(2)}$ and $\Sc^{(0)}$.}
  Recall the stratification $\Sc^{(0)}$ of $Z(\RR^p, L)_n$ labelled by unordered
  partitions, see \S \ref{sec:0-cyc}\ref{par:strata-part}
  By Proposition \ref{prop:S2}(d), $\Sc^{(2)}$ refines $\Sc^{(0)}$,
 which we write as $\Sc^{(2)}\prec \Sc^{(0)}$. 
 Let us  describe this refinement more explicitly.

 \begin{prop}
 Assume $p\geq 2$. Then,  $U_M$ and $U_{M'}$ lie in the same stratum of $\Sc^{(0)}$
 if an only if $M\equiv_\Ano N$, see \S \ref{sect:cont}\ref{par:ano}
 \end{prop}
 
 \noindent{\sl Proof:} Let $\equiv$ be the equivalence relation on $\CM_n^P(L)$ given by
 $M\equiv N$ iff $U_M$ and $U_N$ lie in the same stratum of $\Sc^{(0)}$. We need to prove
 that $M\equiv_\Ano N$ if and only if $M\equiv N$.
 
 In one direction, if $M\to N=\del_i^\nu M$ is an elementary anodyne contraction, then the
 nonzero entries of $M$ and $N$ are in bijection, so the unordered partitions $\bn(M)$ and
 $\bn(N)$ are the same, which means that $U_M$ and $U_N$ lie in the same stratum of
 $\Sc^{(0)}$, i.e., $M\equiv N$. Since $\equiv_\Ano$ is the equivalence relation generated
 by the existence of an elementary anodyne contraction, we conclude that $M\equiv_\Ano N$
 implies $M\equiv N$. 
 
 In the other direction, let $C$ be an equivalence class of $\equiv$. Let us prove that
 $C$ is also an equivalence class of $\equiv_\Ano$. Call a matrix $M\in C$ {\em maximal},
 if there is no anodyne contraction $M'\to M$. In this each slice $\sigma_i^\nu M$ contains exactly
 one non-zero entry. Indeed, otherwise we can split this slice into two parts and still get
 a contingency matrix $M'$ with an anodyne contraction $M'\to M = \del_j^\nu M$. 
 This means that $M$ is a ``permutation'' $p$-dimensional matrix of size $r\times\cdots\times r$
 for some $r$, with the only nonzero entries being $m_{i, w_1(i), \cdots, w_{p-1}(i)}$,
 $i=1,\cdots, r$, where $W-\nu\in S_r$ are some permutations of $\{1,\cdots, r\}$.
 By the above, each $M\in C$ is $\equiv_\Ano$-equivalent to a maximal matrix. It remains
 to show that any two maximal matrices are $\equiv_\Ano$-equivalent. To see this,
 note that any maximal matrix $M$ can be contracted in an anodyne fashion to
 a matrix of size $r\times 1\times\cdots\times 1$ which can be viewed as an ordered
 partition $(a_1,\cdots, a_r)$. The $a_i$ are therefore all the nonzero entries of $M$,
 in some order. It remains to show that any two $r\times 1\times\cdots\times 1$ contingency
 matrices different by a reordering, e.g., $(a_1,\cdots, a_r)$ and $(a_{w(1)}, \cdots, a_{w(r)})$,
 $w\in S_r$, are $\equiv_\Ano$-equivalent. But this is a standard manipulation involving
 matrices of size $r\times 2\times 1\times\cdots\times 1$,.e.g.,
 \[
 (a_1, a_2, a_3) \leftarrow 
 \begin{bmatrix} a_1&0&a_3\\0&a_2&0
 \end{bmatrix} \to 
 \begin{bmatrix} a_1&a_2\\ a_2& 0
 \end{bmatrix} \leftarrow
 \begin{bmatrix}
 0&a_1&a_3\\a_2&0&0
 \end{bmatrix} \to
 (a_2, a_1, a_3),
  \]
 the arrows representing anodyne contractions. \qed

 \begin{ex}
 Let $L=\TT=\{\0bar, \1bar\}$. Elements of $\CM^p(\TT)$ will be called
 ($p$-dimensional) {\em Boolean contingency matrices}. Such matrices
 consist of $\0bar$'s and $\1bar$'s and the simplicial contractions
 \eqref{eq:def-del} are defined using Boolean addition. The content of
 any Boolean contingency matrix other then $\emptyset$ is $\1bar$.
 Taking into account Example \ref{exas:L}(e), we obtain a cell decomposition
 of the Ran space $\Ran(\RR^p)=Z(\RR^p, \TT)_\1bar$ labelled by nonempty
 $p$-dimensional Boolean contigency matrices. 
  \end{ex}
  
  \paragraph{$\Sc^{(2)}$ as a twisted realization.}\label{par:S2=realiz}
  We upgrade the argument in the proof of Proposition \ref{prop:S2}(a) into viewing
  $Z(\RR^p, L)_n$ with its cell decomposition $\Sc^{(2)}$ as a sort of
  geometric realization of the $p$-fold augmented semi-simplicial set
  $\ul\CM_n^p(L)$. Consider the closed cone
  \[
  \RR^m_{\leq} \,=\,\{(t_1,\cdots, t_m)\in \RR^m |\, t_1\leq\cdots\leq t_m\}
  \]
  with its stratification into $2^{m-1}$ faces on which some of the inequalities become
  equalities. Note that
  \[
  \RR^{\bullet+2}_\leq : [d] \mapsto \RR^{d+2}_\leq, \quad d=-1,0,1,\cdots
  \]
  is an augmented co-semi-simplicial object in the category $\STop$ of stratified spaces,
  i.e., a functor $\Delta_+\to\STop$. 
  
  For each $\bd = (d_1,\cdots, d_p) \in  \ZZ_{\geq -1}^p = \Ob(\Delta_+^p)$ we denote
  $\RR_\leq^{\bd+2} = \RR^{d_1+2}_\leq \times\cdots\times\RR^{d_p+2}_\leq$. 
  Let $X: (\Delta_+^p)^\op\to\Set$ be any $p$-fold augmented semi-simplicial set, whose
  value on an object $[\bd]\in\Ob (\Delta_+^p)$ will be denoted $X_\bd$. We define the
  {\em augmented realization} of $X$ to be the coend \cite{maclane}
  \[
  \vertiii{X}\,=\,\int^{\bd\in \Delta_+^p} \RR^{\bd+2}_\leq\times X_\bd \, =\,
  \varinjlim_{ (\Delta_+^p[\bd]\to X) \in\Gro (X)^\op}\nolimits^\Top\,\,  \RR^{\bd+2}_\leq, 
  \]
  the colimit being taken over the (opposite to the) Grothendieck construction. Note that
  $\vertiii{X}$ comes with a natural stratifiction $\Sc_{|||}$ glued out of the face stratifications
  of the $\RR^{\bd+2}_\leq$. 
  
  \begin{prop}\label{prop:S2=realiz}
  We have a stratified homeomorphism
  \[
  \bigl( Z(\RR^p,L)_n, \Sc^{(2)}\bigr) \, \= \, \bigl( \vertiii{\ul\CM_n^p(L)}, \Sc_{|||}\bigr). 
  \]
  \end{prop}
  
  \noindent{\sl Proof:} As observed in the proof of Proposition  \ref{prop:S2}(a), the cell
  $U_M \subset Z(\RR^p,L)_n$, $M\in \CM_n(d_1,\cdots, d_p; L)$, is identified with
  $\RR^{d_1}_<\times\cdots\times\RR^{d_p}_<$, the open stratum of 
  $\RR^{d_1}_\leq\times\cdots\times\RR^{d_p}_\leq$. It remains to notice that under these identifications, the embedding
  $U_{\del_i^\nu M} \subset \ol U_M$ is identified with the embedding
  \[
  \RR^{d_1}_\times\cdots\times \RR^{d_i-1}_< \times\cdots\times \RR^{d_p}_< \,\subset \,
  \RR^{d_1}_\leq\times \cdots\times\RR^{d_p}_\leq \,=\, \ol{
  \RR^{d_1}_< \times\cdots\times \RR^{d_p}_< },
  \]
  and similarly for iterated contractions. \qed
  
  \paragraph{The exit path category of the contingency cell decomposition.} 
  Let $(Z,\Sc)$ be a (Whitney) stratified space. We have a partial order or $Z$ given by $z\leq z'$ if
  the stratum containing $z$ lies in the closure of the stratum containing $z'$.
  An {\em exit path} for $(Z,\Sc)$ is a continuous path $\gamma: [0,1]\to Z$
  which is monotonous, i.e., $\gamma(t) \leq \gamma(t')$ for $t\leq t'$.  We refer
  to \cite{treumann, curry, curry-patel} for the definition of the {\em exit path category}
  $\Exit(Z,\Sc)$ (or simply $\Exit(Z)$, if $\Sc$ is clear from the context). Thus, objects
  of $\Exit(Z,\Sc)$ are points of $Z$ and morphisms are appropriate isotopy classes
  of exit paths. 
  
  It was a fundamental observation of MacPherson that morphisms of $\Exit(Z,\Sc)$ act
  on the stalks of $\Sc$-constructible sheaves (see below).
  
  The part of $\Exit(Z,\Sc)$ formed by exit paths staying in a fixed stratum is just the
  fundamental groupoid of this stratum. So one can replace $\Exit(Z,\Sc)$ by an equivalent
  subcategory with the objects corresponding to connected components of the strata. 
  
  The concept of exit paths categories is applicable, in particular, to cell decompositions.
  In this case the fundamental groupoid of each stratum (cell) is trivial, so $\Exit(Z,\Sc)$
  describes different ways (``branches'') in which a cell of higher dimension can approach
  a cell of lower dimension. We are especially interested in the contingency cell
  decomposition $\Sc^{(2)}$ of the space $Z(\RR^p, L)_n$.
  
  \begin{prop}\label{prop:exit-S-2=gro}
  We have an equivalence
  \[
  \Exit(Z(\RR^p, L)_n, \Sc^{(2)}) \,\=\, \Gro(\ul\CM_n^p(L))^\op. 
  \]
  \end{prop}
  
  \noindent{\sl Proof:} We establish the following general fact which in virtue of Proposition 
  \ref{prop:S2=realiz}, will imply our statement.
  
  \begin{prop}\label{prop:exit-real-gro}
  Let $X: (\Delta_+^p)^\op\to\Set$ be any $p$-fold augmented semi-simplicial set. Then we have
  an equivalence
  \[
  \Exit(\vertiii{X}, \Sc_{|||}) \,\=\, \Gro (X)^\op. 
  \]
  \end{prop}
  
  \noindent {\sl Proof of Proposition \ref{prop:exit-real-gro}:} Consider first the case when $X$
  is a standard multi-simplex, i.e., a representable functor
  \[
  \Delta_+^p [\bd] \,=\, \Hom_{\Delta_+^p} (-, [\bd]), \quad \bd\in \ZZ_{\geq -1}^p. 
  \] 
  In this case $\vertiii{\Delta^p_+[\bd]} = \RR^{\bd+2}_{\leq}$ with its stratification (cell
  decomposition) by faces. This cell decomposition is quasi-regular, see \cite{KS-arr} ???
  and so the corresponding category of exit paths is just the poset of cells ordered by 
  inclusion of closures. This poset is canonically identified with $\Gro(\Delta^p_+[\bd])^\op$.
  Indeed, the objects of the latter category are, by definition, all the multi-simplices
  of $\Delta_+^p[\bd]$, i.e., faces of $\RR^{\bd+2}_{\leq}$, and the morphism between
  two such multisimplices (a contraction) is unique if it exists, which happens precisely
  when we have an inclusion of face closures. 
  
  Let now  $X: (\Delta_+^p)^\op\to\Set$ be arbitrary. Note that
  \[
  \Gro(X)^\op \,=\, \int^{\bd\in \Delta_+^p} X_\bd \times \Gro(\Delta_+^p[\bd])^\op 
  \,=\,
    \varinjlim_{ (\Delta_+^p[\bd]\to X) \in\Gro (X)^\op}\nolimits^\Cat\,\, \Gro(\Delta_+^p[\bd])^\op,
  \]
  the coend of colimit being taken in the category $\Cat$ of small categories and $X_\bd$
 begin considered as a discrete category. This is a general fact exhibiting the Grothendieck
 construction of any $\Set$-valued functor as the colimit (over itself) of the Grothendieck
 constructions of representable functors. 
 
 Now, considering $\Exit$ as a functor $\STop\to\Cat$, the above analysis of
 $\Exit\vertiii{\Delta^p_+[\bd]}$ gives, by the universal property of $\varinjlim$, a functor
 $\xi: \Gro(X)^\op\to \Exit \vertiii{X}$. We claim that $\xi$ is an equivalence. This is a
 van Kampen-type statement: that the functor $\Exit$ commutes with colimits of a certain type.
 It does not seem to be covered by the standard formulations of the ``van Kampen theorem
 for exit paths'' \cite{treumann, curry, lurie-HA} which treat  colimits of
 spaces arising from open covers. So we give a direct proof.
 
 Indeed, $\xi$ is a bijection on objects (we take one base point in each cell of $\vertiii{X}$,
 and these cells are in bijection with multi-simplices of $X$). Let us prove that $\xi$ is bijective
 on morphisms.  Let $\sigma\in X_\bd$ be a multisimplex of $X$ and $U_\sigma \subset
 \vertiii{X}$ be the corresponding cell. Then $\ol U_\sigma \subset \vertiii{X^{\leq\sigma}}$,
 where $X^{\leq\sigma}\subset X$ is the p-fold augmented semi-simplicial subset whose
 multisimplices have the form $\phi(\sigma)$, $\phi\in\Hom_{(\Delta_+^p)^\op} ([\bd], [\bd'])$.
 Thus, for $\tau \in X_{\bd'}$ we have $U_\tau\subset\ol U_\sigma$ if and only if 
 $\tau\in X_{\bd'}^{\leq\sigma}$. Now, $\Hom_{\Exit \vertiii{X}}(U_\tau, U_\sigma)$
 can be calculated inside $\ol U_\sigma$. So to prove that
 \[
 \xi_{\tau,\sigma}: \Hom_{\Gro(X)^\op}(\tau,\sigma) \lra \Hom_{\Exit \vertiii{X}}(U_\tau, U_\sigma)
 \]
 is a bijection, we can assume $X=X^{\leq\sigma}$. Thus $\vertiii{X}$ is obtained from
 $\RR^\bd_{\leq} = \vertiii{\Delta^p_+[\bd]}$ by identifying some faces (according to the
 identifications of the faces of $\Delta_+^p[\bd]$ in $X^{\leq\sigma}$). Let
 $\rho: \RR^\bd_\leq \to\vertiii{X}$ be the resulting
 map. Since in our setting
 there are no degeneracies (all morphisms in $\Delta_+$ are monomorphisms), 
 $\rho$ does not contract any faces to lower dimension
 but can only identify several faces of the same dimension. In particular, $\rho$ is proper with
 finite fibers. 
 
 Let us now construct an inverse $\eta_{\tau,\sigma}$ to $\xi_{\tau,\sigma}$. 
 Let $\gamma: [0,1]\to\vertiii{X}$ be an exit path with $\gamma(0)\in U_\tau$,
 $\gamma(1)\in U_\sigma$. By definition, this means that there are  finitely many ``exit times''
 $0=t_0 < t_1 <\cdots < t_k=1$ such that $\gamma|_{(t_i, t_{i+1}]}$ lies in the same cell, say, $U_{\tau_i}$,
 $\tau_i\in X_{\bd^{(i)}}$, and 
 \[
 U_\tau \subset \ol U_{\tau_1}, \, U_{\tau_1}\subset \ol U_{\tau_2},\,
  \cdots, U_{\tau_{k-1}}\subset \ol U_{\tau_k}, \,\,\tau_k=\sigma.
 \]
 Now, since $\rho$ induces a homeomorphism $\RR^\bd_<\to U_\sigma$ and $\gamma(t)\in U_\sigma$
 for $t\in (t_{k-1}, 1]$, we consider the path $\rho^{-1}(\gamma): (t_{k-1}, 1]\to \RR^\bd_<$. Because
 the the above properties of $\rho$ we conclude that $\rho^{-1}(\gamma)$ has a well defined limit
 $\lim_{t\to t_{k-1}^+} \rho^{-1}(\gamma(t))\in \RR^\bd_\leq$, which we denote $(\rho^{-1}\gamma)(t_{k-1})$.
 It lies in some face of $\RR^\bd_\leq$, which corresponds to some $\phi_{k-1}: \Delta^p_+[\bd^{(k-1)}]\to \Delta^p_+[\bd]$,
 so we denote this face $\phi_{k-1}(\RR^\bd_\leq)$. Further, $\rho$ induces a homeomorphism of the interior of this
 face and the face $U_{\tau_{k-1}}$, so we extend the lift of $\gamma$ to $\gamma_{k-2}: (t_{k-2}, 1]\to \RR^\bd_{\leq}$
 such that $\rho(\gamma_{k-2}(t))=\gamma(t)$. Then, as before, we find $\lim_{t\to t_{k-2}^+} \gamma_{k-2}(t)\in\RR^\bd_\leq$
 and so on. In this way we lift $\gamma$ to $\gamma_0: [0,1]\to \RR^\bd_\leq$ which is an exit path for the standard
 (face) stratification of $\RR^\bd_\leq$. By the above, such path gives a morphism $\phi: [\bd']\to[\bd]$ in $(\Delta_+^p)^\op$
 with $\phi(\tau)=\sigma)$. So $\phi$ can be seen as a morphism from $\tau$ to $\sigma$ in $\Gro(X)^\op$. 
 We define $\eta_{\tau,\sigma}(\gamma)=\phi$.
 
 It is clear that this construction is invariant under admissible isotopies of exit paths, so we get a well defined map
 \[
 \eta_{\tau, \sigma}: \Hom_{\Exit \vertiii{X}} (U_\tau, U_\sigma) \lra \Hom_{\Gro(X)^\op}(\tau,\sigma).
 \]
It is further clear that $\eta_{\tau, sigma}$ is inverse to $\xi_{\tau,\sigma}$. 

\begin{exas}\label{exas:cont-cells-sym-ran}
(a) Let $L=\ZZ_+$, so $Z(\RR^p,\ZZ_+)_n=\Sym^n(\RR^p)$. In this case the contingency cell decomposition $\Sc^{(2)}$
is quasi-regular (see \cite{KS-cont} for the case $p=2$, the case $p>2$ is similar). So $Exit(\Sym^n(\RR^p), \Sc^{(2)})$
is identified with the poset of cells. This is in agreement with Example \ref{ex:S-2-sym} which identifies
 $\Gro(\ul\CM_n^p(\ZZ_+))^\op$ as a poset.
 
 \vskip .2cm
 
 (b) Let $L=\TT=\{\0bar, \1bar\}$ and $p=1$. For $d\geq 1$ let $(\1bar^d)$ be the ordered partition $(\1bar,\cdots, \1bar)$
 ($d$ times). The cell $U_{(\1bar^d)}\subset Z(\RR, \TT)_\1bar = \Ran(\RR)$ consists of all $d$-element subsets
 $\{t_1, \cdots, t_d\}\subset \RR$, where we can assume $t_1 <\cdots < t_d\}$. We see that $U_{(\1bar^d)}$ can approach
 $U_{(\1bar^{d-1})}$ along $d-1$ ``branches'' which correspond to $t_i$ merging with $t_{i+1}$, $i=1, \cdots, d-1$.
 These ``branches'' correspond to face morphisms $[d-3]\to[d-2]$ in $\Delta_+ = \Gro(\ul\CM_\1bar^1(\TT))^\op$. 
\end{exas}
  

    \section{$0$-cycles in $\RR^2=\CC$: further features}\label{sec:X=C}
  
  \paragraph{The case $p=2$, or $X=\CC$.}\label{par:X=C} 
  
   From now on we further specialize to $X=\RR^2$ which we identify with the complex plane $\CC$
  via $(x_1, x_2) \mapsto x_1+ix_2$. Accordingly, all contingency matrices will be
  assumed $2$-dimensional, unless indicated otherwise and we use
  the notation $\CM_n(L)$ for $\CM^2_n(L)$ etc. We will  denote
  the simplicial contractions $\del_i^1, \del_i^2$ by $\del'_i, \del''_i$
 and write $\Gro'(\ul\CM_n(L))$, $\Gro''(\ul\CM_n(L))$ for $\Gro^1(\ul\CM_n^2(L))$ and $\Gro^2(\ul\CM_n^2(L))$ 
  respectively. Similarly, the equivalence relations $\equiv^1_\Ano, \equiv^2_\Ano$ will be denoted by
  $\equiv'_\Ano, \equiv''_\Ano$. 
    This agrees with
  the notation of \cite{KS-cont, KS-prob} as well as of \cite{KS-hW}.  
 We  adapt some further
   considerations of \cite[\S 1]{KS-hW} to our situation. 
   
     \vskip .2cm
     
     Let $\phi: [d]\to [d']$, $d, d'\geq -1$, be a morphism in $\Delta_+^\op$ (i.e., a composition of the face maps $\del_j$).
     Let $M\in\CM_n(d+2,q; L)$, $N\in \CM_n(d'+2, q; L)$ be contingency matrices such that
      $(\phi\times\Id)\in \Mor(\Delta_+^\op
     \times\Delta_+^\op)$ takes $M$ to $N$. We will denote the morphism $M\to N$ in $\Gro'(\ul\CM_n(L))$ represented
     by $\phi\times \Id$ by $M\buildrel\phi'\over\to N$. Similarly, if $M\in \CM_n(p, d_2; L)$ abd $N\in\CM_n(p, d'+2; L)$
     we such that $(\Id\times\phi)\in \Mor(\Delta_+^\op\times\Delta_+^\op$ takes $M$ to $N$, then we denote the
     resulting morphism $M\to N$ in $\Gro''(\ul\CM_N(L))$ by $M\buildrel\phi''\over\to N$. 
     This is in agreement with the notation $\del'_j, \del''_j$ when $\phi=\del_j$ is a generating morphism of $\Delta_+^\op$.

     \paragraph{Mixed forks and their suprema.}\label{par:forks}
      By a {\em mixed fork} in $\Gro(\ul\CM_n(L))$ we will mean a diagram
     $P=\{ M'\buildrel \phi'\over\leftarrow N'\buildrel \psi''\over\to N\}$ in which the first morphism lies in $\Gro'(\ul\CM_n(L))$
     and the second one in $\Gro''(\ul\CM_n(L))$. We call the {\em mixed supremum} of $P$ and denote $\Sup(P)$
     the set of contingency matrices $M$ fitting into a commutative diagram in $\Gro(\ul\CM_n(L))$:
     \be\label{eq:sup-diag}
     \xymatrix{
     M\ar[r]^{\phi'} \ar[d]_{\psi''}
      & N \ar[d]^{\phi''}
     \\
     M'\ar[r]_{\phi'} & N. 
     }
     \ee
     This means that $N$ is obtained from $M$ by the same horizontal simplicial contraction (namely $\phi\times\Id$ corresponding
     to $\phi\in\Mor(\Delta_+^\op)$) as $N$ is obtained from $M'$, and similarly in the vertical direction. Thus, if $N'$ is of size
     $p\times q$, then $M'$ is of size $p_1\times q$, $p_1\geq p$, while $N$ is of size $p\times q_1$, $q_1\geq q$ and $M$ must of
     of size $p_1\times q_1$. 
 The following is an analog of \cite[Prop. 1.11]{KS-hW}. 
 
 \begin{prop}
 Suppose that in a mixed fork $P=\{ M'\buildrel \phi'\over\leftarrow N'\buildrel \psi''\over\to N\}$ one of the
 morphisms is anodyne. Then $\Sup(P)$ consists of precisely one element, i.e., there is a unique diagram
 \eqref{eq:sup-diag}. Further, if $M'\buildrel\phi'\over\to N'$ is anodyne, then $M\buildrel \phi'\over\to N$
 is anodyne. If $N\buildrel \psi''\over\to N'$ is anodyne, then $M\buildrel\psi''\over\to M'$ is anodyne.
 \end{prop}
 
 \noindent{\sl Proof:} We first consider the case when $\phi$ and $\psi$ are generating morphisms of
 $\Delta_+^\op$. That is, $\phi=\del_i$ for some $i$, so $\phi'=\del'_i$ sums the $(i+1)$st and $(i+2)$nd rows of 
 a matrix, while $\psi=\del_j$ for some $j$, so $\psi''=\del''_j$ sums the $(j+1)$st and $(j+2)$nd columns.
 Suppose that $M'\buildrel \del'_i \over\to N'$ is anodyne (the other case is treated similarly). Note that all the entries
 of any possible $M\in\Sup(P)$ are determined uniquely except for the $2\times 2$ submatrix
 \be\label{eq:2x2-indet-mat}
 \begin{gathered}
 \begin{bmatrix}\
 m_{i,j}& m_{i,j+1} \\ m_{i+1, j} & m_{i+1, j+1}
 \end{bmatrix}, \quad 
 \begin{matrix}
 m_{i,j} + m_{i, j+1} = n_{i,j} ,
 \\
 m_{i+1, j+1} + m_{i+1, j+1} = n_{i+1, j},
 \end{matrix}
 \\
 m_{i,j} + m_{i+1, j} = m'_{i,j}, \quad m_{i, j+1} + m_{i+1, j+1} = m'_{i+1, j},
 \\
 m'_{i,j} + m'_{i+1, j} = n_{i,j} + n_{i, j+1} = n'_{i,j}. 
 \end{gathered}
 \ee
 By the assumption that $M'\buildrel \del'_i\over\to N'$ is anodyne, one of the $m'_{i,j}, m'_{i+1, j}$ is $0$.
 Suppose $m'_{i,j}=0$, the other case is treated similarly. Then $m_{i,j} + m_{i,+1 j}= m'_{i,j}=0$
 which implies $m_{i,j}=m_{i+1, j}=0$. The conditions \eqref{eq:2x2-indet-mat} imply that 
 $m_{i, j+1}=n_{i,j}$ and $m_{i+1, j+1} = n_{i, j+1}$. This determines the $2\times 2$ matrix uniquely. 
 Further, if we define it by the above rules and so define the matrix $M$, then $M$ will be a contingency
 matrix. Indeed, the only possible thing to check is that the $i$th column of $M$, where we put two zero entries, is not entirely $0$. But the rest of this column is identical with the rest of the $i$th column of $M$, and this
 rest cannot be $0$ since $m'_{i,j}=0$ and $M'$ is a contingency matrix. 
 This shows that there exists exactly one $M$ fitting into a diagram  \eqref{eq:sup-diag}.  Further, the fact that
 $M\buildrel \del'_i\over\to N$ is anodyne, follows from $m_{i,j}= m_{i, j+1}=0$.
 This proves the proposition in the case when $\phi=\del_i, \psi=\del_j$.
 
 In the general case we write $\phi=\del_{i_1} \cdots \del_{i_p}$, $\psi=\del_{j_1}\cdots\del_{j_q}$.
 If $M'\buildrel\phi'\over\to N'$ is anodyne (the case when $N\buildrel \psi''\over\to N'$ is anodyne,
 is treated similarly), then all the contraction induced by the $\del_{i_\nu}$ are anodyne. So using the
 particular case treated already we recover, in a unique way, a commutative diagram in the form of a
 $p\times q$ grid of rectangles with horizontal arrows being anodyne $\del_{i_\nu}$ and
 vertical arrows being the $\del_{j_\mu}$. Our required diagram is then recovered as the outside border
 of this grid with arrows in each direction composed. \qed

 \paragraph{The Fox-Neuwirth-Fuchs cell decomposition $\Sc^{(1)}$ of $Z(\CC,L)$.}\label{par:FNF}

 For an ordered partition $\bm =(m_1,\cdots, m_d)\in \CM^1(L)$
 we write $l(m)=d$ and call this number the {\em length} of $\alpha$. 

Let such an $\bm$ be given.  
  Further, for each $i=1,\cdots, d$ let $\bn^{(i)}\in \CM^1_{\alpha_i}(L)$ be
  an ordered partition of $n_i$. We write $\vbn$
  for the sequence (vector) of ordered partitions
   $ (\bn^{(1)}, \cdots, \bn^{(d)})$.
  
  Denote by $\CM_{[\bm|\vbn]}(L)$ the set of contingency matrices $M=\| m_{ij}\|\in \CM_n(L)$
  with the following properties:
  \begin{itemize}
  \item[(1)] $M$ has exactly $d$ rows (and some  number $q$ of columns)
  and for each $i=1,\cdots, d$ the $i$th row sum is equal to $\sum_j m_{ij}=m_i$. 
  
  \item[(2)] The ordered partition of $m_i$ given by omitting the zero entries from among the
  $m_{i1}, \cdots, m_{iq}$, is equal to $\bn^{(i)}$. 
    \end{itemize}

    Let $U_{[\bm|\vbn]}\subset Z(\CC,L)$ be the union of the contingency cells $U_M$ for 
    $M\in \CM_{[\bm|\vbn]}(L)$. 
    
    \begin{prop}
    Each $U_{[\bm| \vbn]}$ is a topological cell of dimension  $l(\bm)+\sum l(\bn^{(i)})$. 
    \end{prop}
    
    \noindent{\sl Proof:} Indeed, $U_{[\bm|\vbn]}$ fibers over the
    contingency cell $U_{\bm}^{\RR}\subset Z(\RR,L)$
    with fibers  homeomorphic to the product $\prod U_{\bn^{(i)}}^{ \RR}$
    of similar cells, i.e., to a cell. Therefore $U_{[\bm|\vbn]}$ is
    a cell itself.  \qed
    
    \vskip .2cm
    
    Extending the terminology of 
     \cite{KS-cont}, we  will call the  $U_{[\bm|\vbn]}$ the {\em Fox-Neuwirth-Fuchs cells}
    of $Z(\CC,L)$. The stratification (cell decomposition) of $Z(\CC,L)$ given by these cells
    will be denoted $\Sc^{(1)}$.  
    
      \begin{prop}\label{prop:equiv-ano}
   The sets $\CM_{[\bm|\vbn]}(L)$ are precisely the equivalence classes of $\equiv_\Ano'$. 
    \end{prop}
    
    \noindent {\sl Proof:} We first prove that $M\equiv'_\Ano N$ implies that
    $M$ and $N$ lie in the same $\CM_{[\bm|\vbn]}(L)$. For this it suffices
    to assume that there is anodyne contraction $N\buildrel\phi'\over\to M$, $\phi\in\Mor(\Delta_+^\op)$. 
 This implies  that
  $M$ is obtaines from $N$ by adding some columns.
    This does not affect the row sums, so such sums for $M$ and $N$ are
    equal, giving the same ordered partition $\bm$.
    Further, any anodyne contraction $\del'_i$ does not
    change the pattern of the nonzero entries of the rows of the matrix.
    Thus the ordered partitions $\bn^{(i)}$ associated to the rows of $M$ and
    $N$ are the same. This shows that $M,N$ lie in the same $\CM_{[\bm|\vbn]}(L)$.
    
     Let us prove the converse statement: that whenever $M,N\in \CM_{[\bm|\vbn]}(L)$,
     we necessarily have $M\equiv'_\Ano N$. For this, we look at 
       the cell decomposition
     of $U_{[\bm|\vbn]}$ into the contingency cells $U_M$
     (labelled by the set $\CM_{[\bm|\vbn]}(L)$).  By Proposition \ref{prop:exit-S-2=gro}
       any inclusion of
     the cell closures $U_M\subset \ol U_N$ for $M,N\in \CM_{[\bm|\vbn]}(L)$
     means that there is a morphism from $N$ to $M$ in $\Gro(\ul\CM_n(L))$.
     Because $M,N$ both lie in $\CM_{[\bn | \vbn]}(L)$, such a morphism must
     fave the form of an anodyne contraction $N\buildrel\phi'\over\to M$, $\phi\in\Mor(\Delta_+^\op)$. 
       Since  $U_{[\bm|\vbn]}$ (being itself a cell)
     is connected, the equivalence relation on cells of its cell decomposition
     generated by the relation of
     inclusion of the cell closures, has only one equivalence class. \qed

    \begin{rem}
   One can generalize the Fox-Neuwirth-Fuchs cell decomposition to $Z(\RR^p, L)$ for any
   $p\geq 2$, following the approach of Vassiliev \cite[\S I.3.3]{vassiliev}. Since in this paper
   we are interested only in the case $p=2$, we do not pursue this generalization. 
    \end{rem}

  \paragraph{Imaginary strata of $Z(\CC,L)$.}

  A continuous map $f: X\to Y$ of topological spaces gives the direct image map
  $f_*: Z(X,L) \to Z(Y,L)$ given by
  \[
  f_*\biggl( \sum_{x\in X} n_x \cdot x\biggr) \,=\, \sum_{y\in Y} m_y\cdot y, \quad m_y = \sum_{f(x)=y} n_x. 
  \]
In particular, the imaginary part map $\Im: \CC\to\RR$  gives $\Im_*: Z(\CC,L)\to Z(\RR,L)$. 
  Now, $Z(\RR,L)$ has its own cell decomposition $\Sc^{(2)} = \Sc^{(2)}_\RR$, see \S\ref{par:S2} 
  It consists of cells
    $U^\RR_\alpha$ labelled by ordered partitions 
  ($1$-dimensional contingency matrices) $\bn = (n_1,\cdots, n_d)$. We denote
  \be
  Z_\bn^\Im \,=\, \Im_*^{-1}(U_\bn^\RR)\,\subset \, Z(\CC,L)
  \ee
  the preimage of the cell $U_\bn$ under $\Im_*$.
  It is a locally closed subset but not generally a manifold. 
  We will call  the $Z_\bn^\Im$ 
  the {\em imaginary strata} of $Z(\CC,L)$. If $n=|\bn|=\sum n_i$, then $Z_\bn^\Im$
  lies in the component $Z(\CC,L)_n$. 
  We denote the stratification (not satisfying the Whitney conditions)
  of $Z(\CC,L)$ into imaginary strata by $\Ic$. 
  
  \vskip .2cm
  
 For $\bn$ as above let $\CM_\bn(L)\subset \CM(L)$ be the set
 of contingency   matrices $M=\|m_{ij}\| \in\CM(L)$ with exactly $d$ rows
  summing to   
  $\sum_j m_{ij} = n_i$.

  \begin {prop}
  (a) Suppose $\bn=(n_1,\cdots, n_d)$. The set $Z_\bn^\Im$ is the union
  of  contingency cells $U_M$ for $M\in\CM_\bn(L)$  
  \vskip .2cm
  
  (b) Let $\equiv'$ be the equivalence relation on $\CM(L)$
  generated b y the relation of existence of a morphism in $\Gro'(\ul\CM_n(L))$. 
    The sets $\CM_\bn(L)$ for all $\bm\in\CM^1(L)$
  are precisely the equivalence classes of  $\equiv'$, i.e., the connected components of $\Gro'(\ul\CM_n(L))$. 
  
  \vskip .2cm
  
  (c) We have the following refinement relations
  \[
  \Sc^{(2)} \prec\Sc^{(1)} \prec \Ic, \Sc^{(0)}.  
  \]
 \end{prop}
 
 \noindent{\sl Proof:} (a) is clear from the definition. Part (b) is proved similarly to
 Proposition \ref{prop:equiv-ano}. Finally, (c) follows from (a) and from the
 definition of the cells $U_{[\bm|\vbn]}$ comprising $\Sc^{(1)}$.  \qed

  \begin{exas}  (a) The smallest imaginary stratum in $Z(\CC,L)_n$ is
  $Z_{(n)}^\Im = \RR\times Z(\RR,L)_n$. It corresponds to the ordered partition $(n)$ consisting of
  a single element $n$ and consists of $0$-cycles $\sum_{x\in\CC} n_x\cdot x$, $\sum n_x=n$,
  with the property that all the $x$ such that $n_x\neq 0$ have the same imaginary part. 
  
  \vskip .2cm
  
  (b) Let $L=\TT=\{\0bar, \1bar\}$, so $Z(\CC,\TT)_\1bar = \Ran(\CC)$ is the Ran space of $\CC$. 
  Similarly to Example \ref{ex:partitions-T}, ordered partitions of $\1bar$ have the form $\1bar^d=(\1bar, \cdots, \1bar)$
  ($d$ times), $d\geq 1$. The  corresponding imaginary stratum $Z_{\1bar^d}^\Im$ will be denoted $\Ran(\CC)_d^\Im$.
  It consists of finite nonempty subsets $I\subset \CC$ such that  the  
  set of imaginary parts $\Im(I)\subset \RR$ has cardinality
  exactly $d$.  Note that each $\Ran(\CC)_d^\Im$ is infinite-dimensional and has infinite codimension
  in $\Ran(\CC)$. 
   \end{exas}
   

\section{Bialgebras graded by a monoid and their PROB}\label{sec:bi-PROB}

\paragraph{$L$-graded bialgebras.}\label{par:L-gr-bi}
  Let $L$ be a monoid satisfying Assumptions \ref{ass:L}. 
 Let  also $(\Vc, \x, \b1, R)$ be a $\k$-linear  braided monoidal category.
By an {\em $L$-graded bialgebra}  we mean a collection $A=(A_n)_{n\in L}$ of objects of $\Vc$
together with multiplication and comultiplication morphisms
\[
\mu_{m,n}: A_m\otimes A_n \lra A_{m+n}, \quad \Delta_{m,n}: A_{m+n}\lra A_m\otimes A_n, \quad m,n\in L,
\]
satisfying the following axioms;
\begin{itemize}

\item[(LB0)] (Primitivity)  $A_0=\b1$  and the morphisms
\[
\mu_{0,n}: \b1 \otimes A_n\to A_n, \quad \mu_{n,0}; A_n\otimes \b1\to A_n, \quad
\Delta_{n,0}: A_n\to A_n\otimes \b1, \quad \Delta_{0,n}: A_n\to \b1\otimes A_n
\]
are the identities. 

\item[(LB1)] The $\mu_{m,n}$ are associative.

\item[(LB2)] The $\Delta_{m,n}$ are coassociative.

\item[(LB3)] (Compatibillity) Suppose $ l_1+l_2 \,=\, m_1+m_2\,=\, n$. Then 
  \[
 \Delta_{l_1, l_2} \,\mu_{m_1, m_2} \,=\,  \sum_{O}
 (\mu_{o_{11}, o_{21}}\otimes \mu_{o_{12}, o_{22}})\circ
  (\Id\otimes R_{A_{o_{12}}, A_{o_{21}}}
 \otimes\Id) \circ (\Delta_{o_{11}, o_{12}}\otimes\Delta_{o_{21}, o_{22}}), 
 \]
 the summation being over $2\times 2$
 matrices $O=\begin{pmatrix}
o_{11}&o_{12} \\ o_{21}&o_{22}
\end{pmatrix}$, $o_{ij}\in L$, such that
\[
o_{11}+o_{21}=l_1, \,\, o_{12}+o_{22}=l_2, \quad o_{11}+o_{12}=m_1,\,\, 
o_{21}+o_{22}
=m_2. 
\]
\end{itemize}  

\begin{rems} (a) Note that Assumptions \ref{ass:L} imply that the summation in (LB3) is over a finite set. 

\vskip .2cm

(b) Suppose $\Vc$ has countable direct sums (which we do  not assume in general).
  Then, instead of understanding $A$ as the collection $(A_n)$,
we can form the direct sum $A^\oplus = \bigoplus_{n\in L} A_n$ which is a bialgebra in $\Vc$ graded by $L$.
 The data of $A^\oplus$ as a bialgebra with  grading determines the $\mu_{m,n}$ and $\Delta_{m,n}$
 as the graded components of the multiplication and comultiplication in $A^\oplus$. In this way we can view
 $L$-graded bialgebras as bialgebras graded by $L$ and  satisfying the additional condition (LB0). 
\end{rems}

\begin{ex}\label{ex:A-to-A+}
Suppose that $\Vc$ is abelian, with $\x$ being exact in each argument. Let $B$ be an ungraded bialgebra
in $\Vc$ with unit $e: \b1\to B$ and counit $\eta: B\to \b1$. Let $B_+ = \Ker(\eta)$. Then $B_+$ is an
ideal in $B$ as an algebra, since $\eta$ is an algebra homomorphism. Let us
define $\ol\Delta: B\to B\otimes B$ by $\ol\Delta = \Delta -\Id\otimes e - e\otimes \Id$, i.e., 
$\ol\Delta(x) = \Delta(x) - x\otimes 1 -1\otimes x$ in the
element notation. Then $\ol\Delta$ takes $B_+$ to 
\[
B_+\otimes B_+ \,=\, \Ker \{\eta\otimes \Id\} \cap \Ker\{\Id\otimes\eta\}. 
\]
Indeed, using the element notation again, $(\eta\otimes \Id)(\Delta(x))=x$ by the axiom of counit. At the same
time, if $x$ is an element of $A_+$, then
\[
(\eta\otimes \Id)(x\otimes 1) = \eta(x)=0, \quad (\eta\otimes\Id)(1\otimes x) = x,
\]
so $\Delta(x) - x\otimes 1 -1\otimes x$ is annihilated by $\eta\otimes\Id$. Similarly for $\Id\otimes\eta$. 

Now, putting $A_\0bar = \b1$, $A_\1bar=B_+$, we associate to $B$ a $\TT$-graded bialgebra $A$, where
$\TT$ is the Boolean algebra of truth values, see Example \ref{exas:L}(e). Its only nontrivial data are:
\[
\mu_{\1bar, \1bar} = \mu: B_+\otimes B_+\lra B_+, \quad \Delta_{\1bar, \1bar} = \ol\Delta: B_+\lra B_+\otimes B_+. 
\]
This defines an equivalence
of categories
\[
\bigl\{ \text{Ungraded bialgebras in } \Vc \text{ with unit and counit}\bigr\} \buildrel \=\over\lra
\bigl\{ \TT\text{-graded bialgebras in } \Vc\bigr\}. 
\]
The  functor in the opposite direction is given by
\[
(\b1=A_\0bar, A_\1bar) \,\mapsto \, B=\b1\oplus A_\1bar.
\]
The two functors are quasi-inverse to each other  because for any unital and counital bialgebra
 $B$ the endomorhism $P=e\eta: B\to B$
is a projector with image $B_+$ and $Q=\Id-P$ is a projector with image $\Im(e)\= \b1$. 
\end{ex}

\paragraph{The PROB of $L$-graded bialgebras.}\label{par:bi-PROB}
 We denote by $\Ben^L$
the $k$-linear braided monoidal category generated by the {\em universal $L$-graded bialgebra}
$\ba=(\ba_n)_{n\in L}$. 
In other words, for any other $\k$-linear braided category $\Vc$, we have a bijection between
braided monoidal functors $F: \Ben^L\to\Vc$ and $L$-graded bialgebras $A$ in $\Vc$ given by $A=F(\ba)$. 

Thus, similarly to  \cite {habiro-bottom} \cite{KS-prob},   $\Ben^L$ is generated, as a braided monoidal category,
by the generating objects which are formal symbols $\ba_n$, $n\in L$, with $\ba_0=\b1$. In other words,
objects of $\Ben^L$ are formal tensor products $\ba_\bn=\ba_{n_1} \otimes\cdots \otimes \ba_{n_d}$ corresponding to all ordered partitions $\bn = (n_1,\cdots, n_d)\in\CM^1(L)$,
see \S \ref{sect:cont}\ref{par:cont}
Morphisms are generated by the elementary formal morphisms
\[
\mu_{m,n}: \ba_{m}\otimes \ba_n\to \ba_{m+n}, \quad \Delta_{m,n}: \ba_{m+n}\to \ba_m\otimes\ba_n
\]
as well as the braidings, which are subject only  to the relations following from (LB0)-(LB3) and the axioms
of a braided category. 

We refer to $\Ben^L$ as the {\em (colored) PROB governing $L$-graded bialgebras},
see  \cite{KS-prob} for the discussion of terminology. 
By construction, $\Ben^L$ splits, as a $\k$-linear category, into an orthogonal
direct sum of blocks $\Ben^L_n$, $n\in L$, where $\Ben^L_n$
is the full subcategory on objects 
 $\ba_\bn=\ba_{n_1} \otimes\cdots \otimes \ba_{n_d}$ corresponding to partitions 
 $\bn=(n_1,\cdots, n_d)$ from $\CM^1_n(L)$, i.e., such that $\sum n_i=n$. 
 The tensor product gives functors 
 $\otimes: \Ben^L_m \times\Ben^L_n\to \Ben^L_{m+n}$. 

\begin{exas}
(a) Let $L=\TT$ be the Boolean algebra of truth values, see Example  \ref{exas:L}(e). Since $\TT$ has
only one non-zero element $\1bar$, the braided category $\Ben^\TT$ is generated by one object $\ba_\1bar$.
Thus $\Ben^\TT$ is an uncolored PROB: its objects are $\ba_\1bar^{\otimes d}$, $d\geq 0$.  

\vskip .2cm

(b) Let also $\Ben^\cu$ be the uncolored PROB governing (ungraded) counital and unital bialgebras. 
Thus $\Ben^\cu$ is generated by one formal object $\bb$ and morphisms
\[
\mu: \bb\otimes\bb\to\bb, \quad \Delta: \bb\to\bb\otimes\bb, \quad \eta: \bb\to \b1, \quad e: \b1\to\bb
\]
subject to the axioms of a counital and unital bialgebra. 

It is instructive to compare $\Ben^\TT$ and $\Ben^\cu$. Let $\Ben^{\TT \oplus}$ be the additive envelope of $\Ben^\TT$, i.e., the braided monoidal category
formed by formal direct sums of objects of $\Ben^\TT$. 
\end{exas}

\begin{prop}
We have a braided monoidal functor 
\[
\Ben^\cu \to \Ben^{\TT\oplus}, \quad \bb \mapsto 1 \oplus \ba_\1bar,
\]
which identifies $\Ben^\cu$ with the full subcategory in $\Ben^{\TT\oplus}$ on the objects $(1\oplus \ba_\1bar)^{\otimes d}$. 
\end{prop}

\noindent{\sl Proof:} This is an upgrade of the reasoning of Example \ref{ex:A-to-A+}. That is, the idempotents
$P=e\eta$, $Q=\Id-P$ in   $\Hom_{\Ben^\cu}(\bb, \bb)$ identify $\bb$ with $\b1\oplus \ba_\1bar$.  
This gives an explicit identification
 \be
\Hom_{\Ben^\TT} (\ba_\1bar^{\otimes m}, \ba_\1bar^{\otimes n}) =
Q^{\otimes n} \Hom_{\Ben^\cu}(\bb^{\otimes m}, \bb^{\otimes n}) 
\subset  \Hom_{\Ben^\cu}(\bb^{\otimes m}, \bb^{\otimes n}). 
\ee
 \qed


   \section{Perverse sheaves on the spaces of $0$-cycles in $\CC$}\label{sec:perv-0-cyc}
   
   \paragraph{The constructible derived category.}\label{par:constr-der-cat}
    We fix $n\in L$
   and denote the component $Z(\CC,L)_n$ simply by $Z$. Note that $Z$ can
   be infinite dimensional and we consider it with the exhaustion by finite dimensional skeleta,  see \eqref{eq:exhaustion1}
   \[
   Z_1\subset Z_2\subset \cdots \subset Z \,=\, \varinjlim\nolimits^\Top \, Z_r,
   \quad Z_r := Z_r(\CC,L)\cap Z(\CC,L)_n.
 \]
 We denote by $i_r: Z_N\hra Z$   the (closed) embeddings.
 
 \vskip .2cm
 
   Let $\Vc$ be a $\k$-linear Grothendieck abelian category (not assumed
   monoidal). We can speak about $\Vc$-valued sheaves on $Z$  on any
   topological space $T$. The category of such sheaves will be denoted
   $\Sh(T,\Vc)$. 
   
   Note that $i_{r*}:\Sh(Z_r, \Vc)\to  \Sh(Z,\Vc)$ is an embedding of a full subcategory
   (whose objects are sheaves supported on $Z_r$). Let $\Fc\in\Sh(Z,\Vc)$.
   We say that $\Fc$ is {\em finitely supported}, if $\Fc$ lies in the image of one of the $i_{r*}$.
   We denote $\Sh^{<\oo}(Z,\Vc)$ the full category in $\Sh(Z,\Vc)$ formed by finitely
   supported sheaves. Thus $\Sh^{<\oo}(Z,\Vc)=\Sh(Z,\Vc)$, if $Z$ is finite dimensional.
   
   \vskip .2cm
   
   For any abelian category $\Ac$ we denote by $C^b(\Ac)$ the category of bounded
   complexes over $\Ac$ and by $D^b(\Ac)$ the bounded derived category of $\Ac$.
   
   \vskip .2cm
   
   The category $\Sh^{<\oo}(Z, \Vc)$ is abelian and we denote $D^{<\oo}(Z,\Vc) =
   D^b\, \Sh^{<\oo}(Z, \Vc)$ its bounded derived category.

   \vskip .2cm
   
   Recall the stratifications $\Sc^{(i)}$ of $Z(\CC,L)$, $i\in\{0,1,2\}$. 
   We denote $\Sh^{<\oo}_{\Sc^{(i)}}(Z,\Vc) \subset \Sh^{<\oo}(Z,\Vc)$
   the full subcategory
   formed by $\Sc^{(i)}$-{\em constructible sheaves}, i.e., sheaves
   whose restriction to any stratum of $\Sc^{(i)}$ is locally constant. 
    Let further $D^{<\oo}_{\Sc^{(i)}}(Z,\Vc)
    \subset D^{<\oo}(Z,\Vc)$ and
      be the full subcategory formed by    
 $\Sc^{(i)}$- {\em constructible complexes}, i.e., complexes whose cohomology
 sheaves are $\Sc^{(i)}$-constructible. Since $\Sc^{(2)}\prec \Sc^{(1)} \prec\Sc^{(0)}$,
 we have embeddings of full subcategories
 \[
 \Sh^{<\oo}_{\Sc^{(0)}}(Z,\Vc) \subset  \Sh^{<\oo}_{\Sc^{(1)}}(Z,\Vc) \subset  \Sh^{<\oo}_{\Sc^{(2)}}(Z,\Vc),\quad
  D^{<\oo}_{\Sc^{(0)}}(Z,\Vc)  \subset D^{<\oo}_{\Sc^{(1)}}(Z,\Vc)  \subset D^{<\oo}_{\Sc^{(2)}}(Z,\Vc).   
 \]

 \paragraph{Cellular description of constructible categories.}
 Recall the category $\Gro(\ul\CM_n(L))$, see \S \ref{sect:cont}\ref{par:gro}
 Let $\Fun^{<\oo}(\Gro(\ul,\CM_n(L))^\op, \Vc)$ be the category of functors 
 $\Gro(\ul,\CM_n(L))^\op \to \Vc$ which vanish on almost all objects. 
 
 \begin{prop}\label{prop:constr-sh-gro}
 (a) We have equivalences of categories
 \[
 \Sh_{\Sc^{(2)}}(Z,\Vc) \,\=\, \Fun(\Gro(\ul\CM_n(L))^\op, \Vc), \quad 
  \Sh^{<\oo}_{\Sc^{(2)}}(Z,\Vc) \,\=\, \Fun^{<\oo}(\Gro(\ul,\CM_n(L))^\op, \Vc).
 \]
 (b) The category $\Sh^{<\oo}_{\Sc^{(0)}}(Z,\Vc)$ is equivalent to the full subcategory in 
 $\Fun^{<\oo}(\Gro(\ul,\CM_n(L))^\op, \Vc)$ consisting of functors that take any anodyne
 morphism to an isomorphism.
 
 \vskip .2cm
 
 (c) The category $\Sc^{<\oo}_{\Sc^{(1)}}(Z,\Vc)$ is equivalent to the full subcategory in 
 $\Fun^{<\oo}(\Gro(\ul,\CM_n(L))^\op, \Vc)$ consisting of functors  that take any
 anodyne morphism in the subcategory $\Gro'(\ul\CM_n(L))^\op$ to an isomorphism. 
 \end{prop}
 
 \noindent {\sl Proof:} (a) The first equivalence follows from Proposition \ref{prop:exit-S-2=gro} and
 the general result \cite{treumann, curry, curry-patel} describing constructible sheaves on any stratified
 space in terms of representations of the category of exit paths. The second equivalence is the restriction
 of the first one.  Part (b) follows from (a) and from Proposition \ref{prop:S2}, since the strata of $\Sc^{(0)}$
 are connected. Part (c) follows from (a) and from Proposition \ref{prop:equiv-ano} in a similar way. \qed
 
 \begin{prop}\label{prop:D-constr=compl-Gro}
 (a) The category $D^{<\oo}_{\Sc^{(2)}}(Z,\Vc)$ is equivalent to $D^b \Fun^{<\oo}(\Gro(\ul\CM_n(L))^\op, \Vc)$,
 the bounded derived category of the abelian category $ \Fun^{<\oo}(\Gro(\ul\CM_n(L))^\op, \Vc)$. 
 
 \vskip .2cm
 
 (b) The category $D^{<\oo}_{\Sc^{(0)}}(Z, \Vc)$ is equivalent to the full subcategory in  $D^b \Fun^{<\oo}(\Gro(\ul\CM_n(L))^\op, \Vc)$
 consisting of complexes of functors (or, equivalently, functors $\Gro(\ul\CM_n(L))^\op\to C^b(\Vc)$)
 taking all anodyne morphisms to quasi-isomorphisms. 
 
 \vskip .2cm
 
 (c) The category $D^{<\oo}_{\Sc^{(1)}}(Z, \Vc)$ is equivalent to the full subcategory in  $D^b \Fun^{<\oo}(\Gro(\ul\CM_n(L))^\op, \Vc)$
consisting of functors  $\Gro(\ul\CM_n(L))^\op\to C^b(\Vc)$ taking anodyne morphisms from $\Gro'(\ul\CM_n(L))^\op$
to quasi-isomorphisms. 
 \end{prop}
 
 \noindent{\sl Proof:} It suffices to establish (a), as (b) and (c) are formal consequences. Now, the proof of (a)  goes along the
 standard lines of comparing $\Ext$'s. we present it in some detail since we are dealing with sheaves with values
 in an arbitrary Grothendieck abelian category $\Vc$. 
 
 Any complex of $\Sc^{(2)}$-constructible sheaves has, obviously $\Sc^{(2)}$-constructible cohomology sheaves. This gives
 an exact functor of triangulated categories
 \[
 \xi: D^b \Fun^{<\oo}(\Gro(\ul\CM_n(L))^\op, \Vc) \,\=\, D^b \Sh^{<\oo}(Z,\Vc) \lra D^{<\oo}_{\Sc^{(2)}}(Z, \Vc),
 \] 
 the first identification  being a consequence of Proposition \ref{prop:constr-sh-gro}(a). Let us prove that $\xi$ is
 an equivalence. We separate this into two lemmas.
 
 \begin{lem}\label{lem:I}
 $\xi$ is fully faithful.
 \end{lem}
 
 \begin{lem}\label{lem:II}
 $\xi$ is essentially surjective.
 \end{lem}
 
 \paragraph{Proof of Lemma \ref{lem:I}.}\label{par:proof-lem-I}
 The statement means that for any bounded complexes $\Fc^\bullet, \Gc^\bullet$ from $\Sh^{<\oo}_{\Sc^{(2)}}(Z, \Vc)$ the canonical map
 \[
 \xi_{\Fc^\bullet, \Gc^\bullet}: \Ext^i_{\Sh^{<\oo}_{\Sc^{(2)}} (Z, \Vc)} (\Fc^\bullet, \Gc^\bullet) 
 \lra
 \Hom_{D^{<\oo}_{\Sc^{(2)}}(Z, \Vc)}(\Fc^\bullet, \Gc^\bullet[i])
 \]
 is an isomorphism for all $i$. 
 Since $\Fc^\bullet$ and $\Gc^\bullet$ are supported on finitely many cells, d\'evissage reduces this question to
 the case when $\Fc^\bullet$ and $\Gc^\bullet$ are sheaves in degree $0$ and, moreover, are single-cell sheaves. 
 That is, for $M\in \CM_n(L)$ let $j_M: U_M\to Z$ be the embedding of the corresponding contingency cell.
 Then we can assume that
 \[
 \Fc^\bullet=\Fc = j_{M*} \ul V_{U_M}, \quad \Gc^\bullet = \Gc= j_{N*}\ul W_{U_N}, \quad M,N\in\CM_n(L),
 \]
 where $V,W\in\Ob(\Vc)$ and $\ul V_{U_M}$, $\ul W_{U_N}$ are the corresponding constant sheaves on $U_M$
 and $U_N$ respectively. 
 
 \begin{lem}\label{lem:V^h^M}
 Under the equivalence of Proposition \ref{prop:constr-sh-gro}(a), the sheaf $j_{M*}\ul V_{U_M}$ corresponds to the functor
 \[
 V^{h^M} \,=\, \Map (\Hom_{\Gro(\ul\CM_n(L))}(M, -), V); \Gro(\ul\CM_n(L))^\op\lra \Vc. 
 \]
 \end{lem}
 
 \noindent To explain the notation, for a set $I$ and an object $V\in \Vc$ we have the object
 \[
 V^I \,=\, \Map(I, V) \, := \, \prod_{i\in I} V  \quad\in\quad \Vc,
 \]
 the product of $I$ copies of $V$. So the covariant set-valued functor
 \[
 h^M\,=\,\Hom(M, -): \Gro(\ul\CM_n(L)) \lra\Set
 \]
 gives a contravariant $\Vc$-valued functor $V^{h^M} : \Gro(\ul\CM_n(L))^\op \lra\Vc$. 
 
 \vskip .2cm
 
 \noindent{\sl Proof of Lemma \ref{lem:V^h^M}:} This is a generality. Denote for short $\Ec = \Gro(\ul\CM_n(L))$ (the ``entrance
 path'' category of $\Sc^{(2)}$). Let $G: \Ec^\op\to\Vc$ be the  functor corresponding to $j_{M*} \ul V_{U_M}$
 by Proposition   \ref{prop:constr-sh-gro}(a). Note that $j_{M*} \ul V_{U_M}$ is characterized by the natural
 identification
 \[
 \Hom(\Fc, j_{M*} \ul V_{U_M}) \,\=\, \Hom(j_M^* \Fc, \ul V_{U_M})
 \]
 for any $\Fc\in\Sh_{\Sc^{(2)}}(Z, \Vc)$. Such $\Fc$ are in equivalence with functors $F: \Ec^\op\to\Vc$ and
 \[
 \Hom(j_M^*\Fc, \ul V_{U_M}) \,=\, \Hom_\Vc (F(M), V). 
 \]
 So $G$ is characterized by the natural (in $F$) isomorphism
 \[
 \Hom_{\Fun(\Ec^\op, \Vc)}(F,G) \, \= \, \Hom_\Vc(F(M), V).
 \]
 The functor $G=V^{h^M}$ obviously admits such an isomorphism
 (a standard property of representable functors). \qed
 
 \vskip .2cm
 
 Continuing with the proof of Lemma \ref{lem:V^h^M}, we note that
 \be\label{eq:Ext=Map}
 \Ext^\bullet_{\Fun(\Gro(\ul\CM_n(L))^\op, \Vc)} (V^{h^M}, W^{h^N}) \,=\,
 \Map \bigl( \Hom_{\Gro(\ul\CM_n(L))}(M,N), \Ext^\bullet_\Vc(V,W)\bigr). 
 \ee
 Indeed, by passing to injective resolutions (which is possible since $\Vc$, being
 Grothendieck, has enough injectives), this statement reduces to the case when $W$
 is injective, in which case $W^{h^N}$ is also an injective object of the functor category.
 So one can replace $\Ext^\bullet$ by $\Hom$ on both sides of \eqref{eq:Ext=Map}, after
 which the statement becomes a standard property of representable functors.
 
 So to prove that $\xi$ is fully faithful, it suffices to show, again under the assumption
 that $W$ is injective, the identification
 \[
 \Ext^q_{\Sh(Z,\Vc)} (j_{M*} \ul V_{U_M}, j_{N*} \ul W_{U_N}) \,=\,
 \begin{cases}
 \Hom_\Vc(V,W)^{\Hom_{\Gro(\ul\CM_n(L))}(M,N)}, & \text{ if } q=0; 
 \\
 0, & \text{ if } q>0. 
 \end{cases}
 \]
 Here, the case $q=0$ is a rephrasing of the fact that the (abelian) category
 $\Sh_{\Sc^{(2)}}(Z, \Vc)$ is equivalent to $\Fun(\Gro(\ul\CM_n(L))^\op, \Vc)$.,
 i.e., of Proposition \ref{prop:constr-sh-gro}(a). The vanishing of $\Ext^q$ for $q>0$
 is a consequence of the following fact which alone uses the specifics of the situation.
 
 \begin{lem}\label{lem:RjN*=0}
 The sheaf $j_{N*} \ul W_{U_N}$ is quasi-isomorphic to the complex $Rj_{N*} \ul W_{U_N}$, 
 i.e., the higher direct images $R^q j_{N*} \ul W_{U_N}$, $q>0$, are zero. 
 \end{lem}
 
 To see the deduction,  notice that by adjunction in the derived category,
 \[
 \Ext^q_{\Sh(Z,\Vc)}(j_{M*} \ul V_{U_M}, Rj_{N*} \ul W_{U_N}) \,=\,
 \Ext^q_{\Sh(U_N, \Vc)} (j_N^* j_{M*} \ul V_{U_M}, \ul W_{U_N}) \,=\,
 \Ext^q_\Vc(V',W) =0,
 \]
 where $V'$ is the stalk at $U_N$ of the constant  sheaf $j_N^* j_{M*} \ul V_{U_M}$. This 
 last $Ext^q$ vanishes since $W$ is injective. 
 
 \vskip .2cm
 
 \noindent {\sl Proof of Lemma \ref{lem:RjN*=0}:} Choose a finite set $I$ of generators of
 the monoid $L$, so we have a surjective morphism of monoids $p: \ZZ_+^I\to L$.
 For any $\wt n\in \ZZ_+^I$ we then have the induced 
  map of cycle spaces $\pi_{\wt N}:  Z(\CC, \ZZ_+^I)_{\wt n} \to Z(\CC, L)_{p(\wt n)}$. 
 Note that $\pi_{\wt n}$ is a proper morphism with finite fibers.  The space $Z(\CC,  \ZZ_+^I)_{\wt n}$
 has its own contingency cell decomposition into cells $U_{\wt N}$, $\wt N\in \CM_{\wt n}(\ZZ_+^I)$.
 As in Example \ref{exas:cont-cells-sym-ran}(a), this decomposition is quasi-regular, i.e.,
 $\ul U_{\wt N}$ is homeomorphic to a part of a closed Euclidean ball. In particular,
 $R^q j_{\wt N *} \ul W_{U_{\wt N}} =0$ for $q>0$. Now, given $N\in\CM_n(L)$, we can find
 $\wt n\in\ZZ_+^I$ with $p(\wt n)=n$ and $\wt N\in\CM_{\wt n}(\ZZ_+^I)$ such that
 $\pi_{\wt n}$ maps $U_{\wt N}$ to $U_N$ in a homeomorphic way and $\ol U_{\wt N}$
 to $\ol U_N$ in a surjective, finite-to one (by the above) way. 
 Thus 
 \[
 Rj_{N*} \ul W_{U_N}\,=\, R\pi_{\wt n *} Rj_{\wt N *} \ul W_{U_{\wt N}} \,=\, 
 \pi_{\wt n *} j_{\wt N *} \ul W_{ U_{\wt N}}
 \]
 is situated in degree $0$. \qed
 
 This finishes the proof of Lemma \ref{lem:I}. 
 
 \paragraph{Proof of Lemma \ref{lem:II}.}
 By construction, the essential image of $\xi$ contains all single sheaves from 
 $\Sh_{\Sc^{(2)}}(Z, \Vc)$ as well as theirn shifts. Now, every $\Fc^\bullet
 \in D^{<\oo}_{\Sc^{(2)}}(Z, \Vc)$ admits cohomological truncations $\tau_{\leq r} \Fc^\bullet$,
 $r\in \ZZ$, fitting into exact triangles
 \[
 \ul H^r(\Fc^\bullet)[-r-1]\buildrel \eps_r\over\to \tau_{\leq r-1}\Fc^\bullet \to
 \tau_{\leq r} \Fc^\bullet \to \ul H^r(\Fc^\bullet) [-r], 
 \]
 with $\tau_{\leq r}\Fc$ being $0$ for $r\ll 0$ and $\Fc^\bullet$ for $r\gg 0$. 
 Using Lemma \ref{lem:I} we represent, inductively, the object $\tau_{\leq r-1} \Fc^\bullet$
 and the morphism $\eps_r$ as the images under $\xi$ of an object and a morphism 
 (denote this morphis m $e_r$) in the source of $\xi$, after which $\tau_{\leq r} \Fc^\bullet$
 is exhibited as the image under $\xi$ of the cone of $e_r$. This finishes the proof of
 Lemma \ref{lem:II} and of Proposition \ref{prop:D-constr=compl-Gro}.

 \paragraph{More general abelian categories. Criterion of $\Sc^{(0)}$-constructibility.}
 Propositions \ref{prop:constr-sh-gro} and \ref{prop:D-constr=compl-Gro} allow us to eliminate,
 in the remainder of the paper,  the assumption that $\Vc$
 is Grothendieck (needed, a  priori, to handle descent in infinite coverings). 
 Indeed, we can (and will) treat all the constructible sheaves or complexes on $Z$ 
 in a purely combinatorial way,
 as certain functors from $\Gro(\ul\CM_n(L))^\op$ to $\Vc$
 or to  the category of complexes over $\Vc$. 
 For this, $\Vc$
 can be  any abelian category.
 
 Firther, the functors $Rf_*, f^!$ etc. associated to various embeddings $f$ of cellular
 subspaces in $(Z, \Sc^{(2)})$ can be defined in purely cellular terms, i.e., directly at
 the level of functors from $\Gro(\ul\CM_n(L))^\op$ or its subcategories, so we will
 use them, often without further explanation.

  Let $\tau:\CC\to\CC$ be the permutation of coordinates: $\tau(x+iy)=y+ix$. The stratifications
  $\Sc^{(0)}$ and $\Sc^{(2)}$ are preserved by $\tau$ but $\Sc^{(1)}$  is taken into a different
  stratification $\tau\Sc^{(1)}$. Similarly to \cite{KS-cont}, Cor. 5.9, we have:
  
  \begin{prop}\label{prop:S2=S1+tS1}
  Let $\Fc$ be an $\Sc^{(2)}$-constructible complex on $Z$. The following are equivalent:
  \begin{itemize}
  \item[(i)] $\Fc$ is $\Sc^{(0)}$-constructible.
  
  \item[(ii)] $\Fc$ is both $\Sc^{(1)}$-constructible and $\tau\Sc^{(1)}$-constructible. 
  \end{itemize}
 
  \end{prop}
  
  \noindent{\sl Proof:} This is a direct consequence of Proposition
  \ref{prop:equiv-ano} and the fact that the class of anodyne morphism in $\Gro(\ul\CM_n(L))^\op$
  is generated by the subclasses of horizontal and vertical anodyne morphisms. \qed

  \paragraph{Verdier duality.}

 Let $\Vc^\op$ be the category opposite to $\Vc$. It is again abelian.
 Since $(Z_r, \Sc^{(i)})$, $i=0,1,2$,
 is a finite-dimesional stratified space, it possesses Verdier duality
 which is an anti-equivalence between the bounded derived categories of constructible
 complexes of sheaves with values in $\Vc$ and $\Vc^\op$. Because $D^{<\oo}_{\Sc^{(i)}}(Z,\Vc)$
 consists of bounded complexes of sheaves supported on  finite-dimensional skeleta, these
 dualities unite into equivalences
 \[
 \DD: D^{<\oo}_{\Sc^{(i)}}(Z, \Vc) \lra D^{<\oo}_{\Sc^{(i)}}(Z, \Vc^\op)^\op.
 \]
 Let us describe $\DD$ at the level of functors on the exit path category.
 
 \vskip .2cm
 
 Let $F: \Gro(\ul\CM_n(L))^\op\to C^b(\Vc)$ be a functor equal to $0$ on almost all objects.
 Let us define a functor $F^!:  \Gro(\ul\CM_n(L))\to C^b(\Vc)$. First, introduce the following notation.
 For $M\in \CM_n(L)$ let $d_M=\dim_\RR U_M$, i.e., $d_M=p+q$, if $M$ is of size $p\times q$. 
 Let aslo $\OR_M=H^{d_m}_c(U_M, \k)$ be the $1$-dimensional orienttion space of $U_M$.
 Further, let us write   $\phi: M' \mathop{\to}_{d} M$
  to signify that $\phi$ is a morphism $M'\to M$ in $\Gro(\ul\CM_n(L))$ and $d_{M'}=d_M+d$. 
 Now put
 \[
 F^!(M) \,=\,\Tot\biggl\{\OR_M\otimes F(M) \buildrel \delta\over\to 
 \bigoplus_{\phi_1: M_1\mathop{\to}_{1} M} \OR_{M_1}\otimes F(M_1) 
  \buildrel \delta\over\to 
   \bigoplus_{\phi_2: M_1\mathop{\to}_{2} M} \OR_{M_2}\otimes F(M_2) 
  \buildrel \delta\over\to 
\cdots \biggr\},
  \]
 where:
 \begin{itemize}
 \item The double complex (i.e., a complex of complexes) under the sign $\Tot$ has horizontal
 grading such that $\OR_M \otimes F(M)$ is in degree $d_M$. 
 
 \item The matrix element
 \[
 \delta_{\phi_p, \phi_{p+1}}: \OR_{M_p}\otimes F(M_p) \lra \OR_{M_{p+1}}F(M_{p+1})
 \] 
 of $\delta$ is equal to
 \[
 \sum_{\phi: M_{p+1}\to M_p, \,\, \phi_p\phi=\phi_{p+1}} \eps_\phi\otimes F(\phi), 
 \]
 where $\eps_\phi: \OR_{M_p}\to \OR_{M_{p+1}}$ is the identification (coorientation)
 induced by the embedding of the $(d_M+p)$-dimensional cell $U_{M_p}$
 to the closure of the $(d_M+p+1)$-dimensional cell $U_{M_{p+1}}$. 
 \end{itemize}
 
 \begin{prop}
 The differential $\delta$ thus defined, satisfies $\delta^2=0$. 
 \end{prop}

 \noindent{\sl Proof:} If $N\in\CM_n(L)$ has size $r\times s$, then we have an identification
 $U_N\= \RR^r_< \times \RR^s_< \subset \RR^{r+s}$, see the proof of Proposition
 \ref{prop:S2}(a). For each such $N$ let us trivialize $\OR_{U_N}$, i.e.,  identify it with $\k$
 using the above identification of $U_N$ and the standard orientation of $\RR^{r+s}$. 
 After this, we write the summand of $F^!(M)$ corresponding to $\phi_p: M_p \mathop{\to}_p M$
 as $F(M_p)$. 
 
 Let now $\phi_p: M_p \mathop{\to}_p M$, $\phi_{p+1}: M_{p+1} \mathop{\to}_{p+1} M$
 be given. Suppose $M_p$ is of size $r\times s$. The existence of at least one
 $\phi: M_{p+1}\to M_p$ means that one of the two cases hold:
 \begin{itemize}
 \item[(1)] $M_{p+1}$ is of size $(r+1)\times s$ and $\phi=\del'_i$ for some $i=0,\cdots, r$.
 
 \item[(2)] $M_{p+1}$ is of size $r\times (s+1)$ and $\phi=\del''_j$ for some $j=0,\cdots, s$. 
 \end{itemize}
 For such $\phi$ the identification $\eps(\phi)$ of already trivialized orientation spaces is given by
 the standard simplicial alternating sign pattern:
 \[
 \eps(\del'_i) = (-1)^i, \quad \eps(\del''_j) = (-1)^{r+j}. 
 \]
 So we can write $\delta_{\phi_p, \phi_{p+1}}: F(M_p)\to F(M_{p+1})$ as
 \[
 \delta_{\phi_p, \phi_{p+1}} \,=
 \sum_{ i:\,  \del'_i M_{p+1}= M_p \atop \phi_p \del'_i = \phi_{p+1}} (-1)^i\delta'_i \,+ 
 \sum_{j: \, \del''_j M_{p+1}=M_p\atop \phi_p \del''_j = \phi_{p+1}} (-1)^{r+j} \delta''_j,
 \]
 where $\delta'_i=F(\del'_i)$ and $\delta''_j=F(\del''_j)$ are the values of $F$ on elementary
 morphisms in $\Gro(\ul\CM_n(L))$. By the above, only one of the two sums can be nonzero.
 
 Now the property $\delta^2=0$ folllows, in the usual way, from the standard simplicial identities
 $\del'_{i_1} \del'_{i_2} = \del'_{i_2-1} \del_{i_1}$, $i_1 < i_2$, from the similar identities for the
 $\del''_j$ and from the commutativity of $\del'_i$ and $\del''_j$. \qed
 
 \vskip .2cm
 
 Let $\psi: M \mathop{\to}_q N$ be a morphism in $\Gro(\ul\CM_n(L))$. Any $\phi: M_p \mathop{\to}_p M$ gives
 $\psi_{p+q} = \psi\phi_p: M_p \mathop{\to}_{p+q} N$. Because of the grading convention of the complex $F^!(M)$, the
 correspondence $\phi_p \mapsto \psi_{p+q}$ defines a morphism of cmplexes $\psi_*: F^!(M) \to F^!(N)$. This
 defines a functor $F^!: \Gro(\ul\CM_n(L))\to C^b(\Vc)$. 
 
 \begin{prop}\label{prop:verdier-gro}
 Suppose $\Fc\in D^{<\oo}_{\Sc^{(2)}}(Z,\Vc)$ is represented by $F: \Gro(\ul\CM_n(L))^\op\to C^b(\Vc)$.
 Then the Verdier dual complex $\DD(\Fc)\in D^{<\oo}_{\Sc^{(2)}}(Z,\Vc^\op)$ is represented by $F^!$ considered as a functor
 $\Gro(\ul\CM_n(L))^\op\to C^b(\Vc^\op)$. 
 \end{prop}

 \noindent{\sl Proof:} For a point $x\in Z$ let $j_x: \{x\}\to Z$ be the embedding. By general properties of Verdier duality, the stalk of
 $\DD(\Fc)$ at $x$ is $j_x^!\Fc^\bullet$, the ``$!$-stalk'' of $\Fc$. If $M\in\CM_n(L)$, then the $j_x^!\Fc$
  for all $x\in U_M$
  can be identified with
  \[
  R\Gamma_c(U_M, j_M^! \Fc) \,=\, R\Gamma_c(U_M, \k) \otimes R\Gamma(U_M, j_M^!\Fc) \,=\,
  \OR_{U_M}[-d_M]\otimes R\Gamma(U_M, j_M^!\Fc). 
  \]
 Now, $j_M^!\Fc = \ul{R\Gamma}_{U_M}(\Fc)$ is the complex calculating the (hyper)cohomology with support in $U_M$. 
 If $\Fc$ is a single sheaf, corresponding to $F: \Gro(\ul\CM_n(L))^\op\to \Vc$, then the space of {\em sections} with support in $U_M$
 is found as
 \[
 \Gamma_{U_M}(Z,\Fc) \,=\, \Ker\bigl\{ F(M) \lra \bigoplus_{\phi: M_1\to M} F(M_1)\bigr\}, 
 \]
 where the direct sum is over all non-unit morphisms $M_1\to M$ in $\Gro(\ul\CM_n(L))$ (not necessarily $M_1 \mathop{\to}_1 M$). 
 We take the derived functor of this by using the standard simplicial resolution of any functors in terms of injective functors
 $V^{h^M}$ discussed in \S \ref{par:proof-lem-I} This gives an explicit complex
 \[
 R\Gamma_{U_M}(Z,\Fc) \,=\,\biggl\{ F(M) \to \bigoplus_{M^{(1)}\to M} F(M^{(1)}) 
 \to
 \bigoplus_{M^{(2)} \to M^{(1)} \to M} F(M^{(2)}) \to\cdots\biggr\}, 
 \]
 where the direct sum in the $q$th term is over all chains of non-unit morphisms $M^{(q)} \to\cdots\to M$. Denote this complex by $R^\bullet$.
 It has a decreasing filtration $V$ where $V^p R^\bullet$ consists of summands corresponding to chains $M^{(q)} \to\cdots\to M$ with
 $d_{M^{(q)}}\geq d_M+p$. The quotient $V^p R^\bullet/ V^{p+1} R^\bullet$ is the direct sum, over all $\phi_p: M_p \mathop{\to}_p M$, of tensor 
 products $O^\bullet_{\phi_p}\otimes F(M_p)$, where $O^\bullet_{\phi_p}$ is a complex of vector spaces defined as follows.
 Let $\Dec_\bullet(\phi_p)$ be the augmented semi-simplicial set with $\Dec_n(\phi_p)$ being the set of
 decomposition sof $\phi_p$ into the composition of $n+1$ non-unit morphisms $M_p = M^{(n+1)} \to M^{(n)} \to\cdots\to M^{(0)} = M$
 (so $\Dec_{-1}(\phi_p)$ consists of one element). Then $O^\bullet_{\phi_p} = \k[\Dec_\bullet(\phi_p)]^* [-2]$ is the cochain complex of 
 $\Dec_\bullet(\phi_p)$ with coefficients in $\k$, shifted so that its lowest term (dual to $\k[\Dec_{-1}(\phi_p)]=\k$)
 appears in degree $+1$. 
 
 We note that $O^\bullet_{\phi_p}$ is quasi-isomorphic to $\OR_{M_p/M}[-d_{M_p} + d_M]$, the relative orientation
 space $\OR_{M_p/M}= \OR_{M_p} \otimes\OR_M^{-1}$ in degree $D_M-d_{M_p}$. Indeed, suppose $M_p$ is
 of size $a\times b$. Then morphisms in $\Gro(\ul\CM_n(L))$ with source $M_p$ are in bijection with
 faces of the polyhedral cone $\RR^a_\leq \times \RR^b_\leq$, see \S\ref{sec:cont-cells} \ref{par:S2=realiz}
Thus $\phi_p: M_p\to M$ corresponds to some face $\Gamma$, and $\Dec_\bullet(\phi_p)$ is the
augmented nerve of the poset of faces containing $\Gamma$, i.e., the augmented simplicial set whose $n$-simplices,
$n\geq -1$, are chains of faces
$\Gamma\subsetneq \Gamma_0\subsetneq \cdots\subsetneq\Gamma_n \subsetneq \RR^a_\leq \times \RR^b_\leq$
(the empty chain for $n=-1$). Now, removing the augmentation, we obtain a semi-simplicial set whose realization
is a sphere of dimension $\codim(\Gamma)-2$. The top cohomology of this realization is identified with
the relative orientation space $\OR_{(\RR^a_\leq\times \RR^b_\leq)/\Gamma} = \OR_{M_p/M}$
while the $0$th cohomology is killed by the augmentation. Thus
\[
\gr^\bullet_V R^\bullet \,\=\, \bigoplus_{p\geq 0} \bigoplus_{\phi_p: M_p \mathop{\to}_p M} O^\bullet_{\phi_p}\otimes F(M)
\,\buildrel{ \on{qis}.} \over \simeq \, 
\bigoplus_{p\geq 0} \bigoplus_{\phi_p: M_p \mathop{\to}_p M} \OR_{M_p/M}[-d_{M_p}+d_M]
\otimes F(M_p) 
\]
 which, as a graded object (without differential) is identified with the complex $F^!(M)$ but without the
 factor $\OR_M[-d_M]$. Further, it is straightforward to check that the induced differential on the
 associated graded (i.e., the ``differential $\delta_1$ in the $E_1$-term of the spectral sequence of the filtered complex'')
 is given by the same formula $\delta_{\phi_p, \phi_{p+1}}=\sum \eps_\phi\otimes F(\phi)$ as in the definition
 of $F^!(M)$. So restoring the factor  $\OR_M[-d_M]$, we get
 \[
 R\Gamma_c(U_M, j_M^!\Fc) \,\=\, \OR_M[-d_M]\otimes H^\bullet(gr^\bullet_V R^\bullet, \delta_1) \,=\, F^!(M).
 \]
 We omit further details. \qed
 
 \begin{rem}
  Proposition  \ref{prop:verdier-gro} is an extension of the well known cellular interpretation of $\DD$ for
  cellular spaces for which the exit path category is a poset, see, e.g., \cite[Ch.12]{curry}. Our argument
  is likewise general though we restricted the attention to our concrete situation. 
 \end{rem}
 
   
  \paragraph{The perverse t-structure.}\label{par:perv-t-str}
   For any stratum $S$ of the stratification $\Sc^{(0)}$ of $Z$
  we denote by $j_S: S\hra Z$ the embedding. We denote by $\dim(S)$ the complex dimension of $S$. 
  Consider the following standard perversity conditions on $\Fc\in D_{\Sc^{(0)}}(Z,\Vc)$:
  \begin{itemize}
  \item[(P$^{-}$)] For any stratum $S$ of $\Sc^{(0)}$ we have $\ul H^q(j_S^*\Fc)=0$ for $q>-\dim(S)$.
  
    \item[(P$^{+}$)] For any stratum $S$ of $\Sc^{(0)}$ we have $\ul H^q(j_S^!\Fc)=0$ for $q< -\dim(S)$.
  \end{itemize}

  \begin{prop}\label{prop:perv-D-finite}
  (a) Let $D^{<\oo}_{\Sc^{(0)}}(Z,V))^{-}, D^{<\oo}_{\Sc^{(0)}}(Z,V))^{+}\subset D^!_{\Sc^{(0)}}(Z,V))$
  be the full subcategories consisting of objects satisfying the conditions (P$^{-}$) and (P$^{+}$)
  respectively. The pair $(D^{<\oo}_{\Sc^{(0)}}(Z,V))^{-}, D^{<\oo}_{\Sc^{(0)}}(Z,V))^{+})$
  is a t-structure on the triangulated category $D^{<\oo}_{\Sc^{(0)}}(Z,V))$. 
  
  \vskip .2cm
  
  (b) Verdier duality $\DD$ sends $\Dc^{<\oo}_{\Sc^{(0)}}(Z,\Vc)^-$ to $\Dc^{<\oo}_{\Sc^{(0)}}(Z,\Vc)^+$
  and vice versa. 
  
  \end{prop}
  
  \noindent{\sl Proof:} (a) This follows, in a formal and standard way, from the
  fact that the conditions (P$^{-}$) and (P$^{+}$)
  define a t-structure on each $D_{\Sc^{(0)}}(Z_r, \Vc)$ and that
  these t-structures are compatible under the direct images with
  respect to the embeddings $Z_r\hra Z_{r'}$, $r\leq r'$. 
  Part (b) is deduced from the  statements about the $Z_r$ in a similar way. 
  \qed
  
  \vskip .2cm
  
  We denote by $\Perv^{<\oo}(Z, \Vc)$ the heart of the above t-structure
  on $D^{<\oo}_{\Sc^{(0)}}(Z,\Vc)$, i.e., the full subcategory of objects
  satisfying both  (P$^{-}$) and (P$^{+}$).
   It is an abelian category whose
  objects will be called {\em  perverse sheaves on} $Z$ {\em with finite support}. 
  By Proposition \ref{prop:perv-D-finite}(b), Verdier duality identifies $\Perv^{<\oo}(Z, \Vc)$ with
  $\Perv^{<\oo}(Z, \Vc^\op)^\op$.

  \begin{rem}\label{rem:Perv=?}
  The conditions (P$^{\pm}$)  make sense for arbitrary (not necessarily bounded or
  finitely supported) $\Sc^{(0)}$-constructible complexes on $Z$. We do not know whether they
  define a t-structure on the triangulated category formed by such complexes. 
   More specifically, it is not clear whether
  any such complex $\Fc$ includes into a triangle
  \[
  \Fc^{\leq 0} \lra \Fc \lra \Fc^{\geq 1}
  \]
  with $\Fc^{\leq 0}$ satisfying (P$^{-}$) and $\Fc^{\geq 1}[-1]$
  satisfying (P$^{+}$). The standard argument for this (in the case of finite-dimensional
  stratified spaces) starts from the highest (open) stratum and goes downward, each step making correction
  at lower and lower strata, see, e.g.,  \cite[Th. 8.1.27]{hotta}.
   In particular, we can retain the name
  ``perverse sheaves'' for complexes satisfying
   (P$^{-}$) and (P$^{+}$)
but it is not clear whether the full subcategory formed by such objects is abelian. 

The question about the ``correct'' category of perverse sheaves not necessarily supported on a finite
skeleton can be studied for more general ind-stratified spaces and we plan to address it
is a separate paper. 

  \end{rem}

    \paragraph{The main result:   perverse sheaves with finite support on
    $Z(\CC,L)_n$ and the block $\Ben^L_n$.}
    Recall also the block $\Ben_n^L$ of the PROB $\Ben^L$ of $L$-graded bialgebras,
    see \S \ref{sec:bi-PROB}\ref{par:bi-PROB}
    A functor $F: \Ben_n^L\to\Vc$ will be said to have {\em finite support}, 
    if $F$ vanishes on almost all objects of $\Ben_n^L$. We denote $\Fun^{<\oo}(\Ben_n^L, \Vc)
    \subset \Fun(\Ben_n^L, \Vc)$ the full subcategory formed by functors with finite support. 
    The following is our first main result whose proof will occupy the rest of the paper.
        
    \begin{thm}\label{thm:Bn-perv}
    We have an equivalence of categories $\Perv^{<\oo} (Z, \Vc) \= \Fun^{<\oo} (\Ben^L_n, \Vc)$. 
    \end{thm}

    \begin{rem}
    Continuing Remark \ref{rem:Perv=?}, it is natural to expect the existence of a ``correct''
    category $\Perv(Z(\CC,L)_n, \Vc)$ (with no assumption on finite support)
    which is identified with the category $\Fun(\Ben_n^L, \Vc)$ of all functors $\Ben_n^L\to\Vc$.
    As we shall see, the arguments proving Theorem \ref{thm:Bn-perv} will provide
    a certain complex $\Fc^\bullet_F$ of sheaves for an arbitrary functor $F: \Ben_n^L\to\Vc$
   but we do not know a good  sense in which it is perverse. In particular, if $\Vc$ is braided monoidal
    and $B$ is an $L$-graded bialgebra in $\Vc$, then we have a braided functor $F: \Ben^L\to\Vc$
    which gives, for each $n$, a functor $F_n: \Ben_n^L\to\Vc$.
    These functors typically do not have finite support. 
    The complexes  $\Fc^\bullet_{F_n}$, taken together, carry a natural factorization structure
    (see \cite{KS-shuffle}    for the case $L=\ZZ_+$ when all the $Z(\CC,L)_n$ are finite dimensional
    and so the concept of a perverse sheaf is unambiguous). Thus the role of the correct
    category $\Perv(Z(\CC,L)_n, \Vc)$ is that it should contain the components
    of the ``factorizable perverse sheaves" associated to $L$-graded bialgebras. 
    \end{rem}
      
   
   \section{Janus sheaves, graded bialgebras and perverse sheaves}\label{sec:janus}

    \paragraph{Janus sheaves.   } Here we adapt \cite[Def. 2.1]{KS-hW}
    to our situation.  We use the notation and conventions  
     of \S \ref{sec:X=C}\ref{par:X=C} 
 Fix $n\in L$, so $\CM_n(L)$ is the set
    of contingency matrices of content $n$ and $\ul\CM_n(L)$ is the corresponding double augmented
    semi-simplicial set. We recall the notations $\Gro'(\ul\CM_n(L))$ and $\Gro''(\ul\CM_n(L))$
    for the horizontal and vertical Grothemdieck constructions of $\ul\CM_n(L)$.  These are categories
    with the same set of objects, namely $\CM_n(L)$. 
    
    \vskip .2cm
    
    In addition to the standard notation $f: M\to N$ for a morphism in $\Gro'(\ul\CM_n(L))$ we will
    use the notation expained in  \S \ref{sec:X=C}\ref{par:X=C} That is, if $N=(\phi\times\Id)(M)$ 
    where $\phi: [m]\to [n]$ is a morphism in $\Delta_+^\op$ and $M$ has $m+2$ columns, then we denote the
    morphism $f: M\to N$ associated to $\phi\times\Id$ by $f=\phi'$. Similarly, if $N=(\Id\times\phi)M$ where $M$
    has $m+2$ rows, then we denote the morphism $f: M\to N$ in $\Gro''(\ul\CM_n(L))$ associated
    to $\Id\times\phi$ by $f=\phi''$. 
    
     \vskip .2cm

    Let $\Vc$ be a $\k$-linear category
    (not assumed abelian or monoidal).

\begin{Defi}\label{def:MBS}
A {\em Janus presheaf} of type $(L,n)$  with values in $\Vc$ is  a set $E$ of data   consisting 
 of  an object $E(M)\in\Vc$ for each  
$M\in \CM_n(L)$ and morphisms
\begin{itemize}
\item $\delta'_f= \delta'_{f,E}: E(M)\to E(N)$ for any morphism $f: N\to M$ in $\Gro'(\ul\CM_n(L))$.

\item $\delta''_g=\delta''_{g,E}: E(N)\to E(M)$ for any morphism $g: N\to M$ in $\Gro''(\ul\CM_n(L))$,
\end{itemize}
\noindent satisfying the conditions:
\begin{itemize}
\item[(JS1)] The $\delta'_{f}$  form a  functor  
 $\Gro'(\ul\CM_n(L))^\op\lra \Vc$. Similarly, the
$\delta''_{g}$  form a   functor   $\Gro''(\ul\CM_n(L))\lra \Vc$. 

\item[(JS2)] The $\delta'$-  and $\delta''$-maps commute with each other. That is, suppose
$P=\{M'\buildrel \phi'\over\to N'\buildrel \psi''\over\leftarrow N\}$ be a mixed fork in $\Gro(\ul\CM_n(L))$,
see \S  \ref{sec:X=C}\ref{par:forks}  Then
\[
\delta'_{\phi': M'\to N} \delta''_{\psi'': N\to N'} \,=\sum_{M\in\Sup(P) }\delta''_{\psi'': M\to M'} \delta'_{\phi': M\to N}.
\]
where the sum os over $M$ fitting into a diagram \eqref{eq:sup-diag}. 
\end{itemize}

A  Janus presheaf is called a {\em  Janus sheaf},
if it satisfies the following:

\begin{itemize}

\item[(JS3)] If $f: N\to M$ is an anodyne morphism in $\Gro'(\ul\CM_n(L))$, then $\delta'_f$
is an isomorphism. If $g: N\to M$ is is an anodyne morphism in $\Gro''(\ul\CM_n(L))$, then $\delta''_g$
is an isomorphism.
\end{itemize}
 \end{Defi}
 
We denote by  $\JPS(\Vc)=\JPS_{L,N}(\Vc)$ the category of Janus presheaves of type $(L,n)$
with values in $\Vc$ and by $\JS(\Vc)=\JS_{L,N}(\Vc)$ the full subcategory formed by Janus sheaves. 
We say that a Janus sheaf $E$ has {\em finite support}, if $E(M)=0$ for almost all $M\in\CM_n(L)$.
Denote $\JS^{<\oo}(\Vc)=\JS^{<\oo}_{L,n}(\Vc)$ the full subcategory in $\JS_{L,n}(\Vc)$ 
formed by Janus sheaves with finite support.

\begin{ex}
The concept of  a Janus sheaf is analogous to that of  a mixed
Bruhat sheaf \cite[Def. 2.1]{KS-hW} which is defined for a pair $(W,S)$
consisting of a  finite reflection group $W$ (acting on a real vector space $\hen_\RR$)
and a set $S$ of simple reflections generating $W$. In fact, the two concept
have a nonempty intersection, appearing when $L=\ZZ_+$, so $n\in \ZZ_+$. A Janus
of type $(\ZZ_+, n)$ is the same as a mixed Bruhat sheaf
of type $(S_n, S)$ where $S_n$ is the symmetric group acting on $\RR^n$
by permutation and $S$ is the set of simple transpositions $(i,i+1)$. 
\end{ex}

\paragraph{Janus sheaves and graded bialgebras.} Let $\Vc$ be a braided monoidal category.
let $I$ be a finite set and $(V_i)_{i\in I}$ be a family of objects of $\Vc$ labelled by $I$.
Unless $\Vc$ is symmetric, the notation $\bigotimes_{i\in I} V_i$ is ambiguous: although the
tensor products of the $V_i$ in different orders are isomorphic via the braidings, the isomorphisms
are not canonical since there are many braids lifting a given permutation. It has been pointed out by
Deligne that a choice of an embedding $\phi: I\hra\CC$ (a ``$2$-dimensional ordering'') fixes the issues,
defining an object which can be denoted $\bigotimes_{i\in I}^\phi V_i$ uniquely up to a unique
isomorphism. The braiding isomorphisms in the approach arise from the monodromy of the
local system on the space of embeddings $\phi: I\hra \CC$ formed by the $\bigotimes_{i\in I}^\phi V_i$.
A detailed treatment of this construction can be found in \cite[\S 3.2A]{KS-shuffle} and \cite[\S2B]{KS-prob}. 
When the embedding $I\hra\CC$ is clear from the context or explained verbally, we will use the simple notation
$\bigotimes_{i\in I} V_i$. 

\vskip .2cm

Let now $A=(A_n)_{n\in L}$ be a graded bialgebra in $\Vc$. For any contingency matrix $M=\|m_{ij}\|$
of size $p\times q$ we define
\[
A_M \,=\bigotimes_{(i,j)\in \{1,\cdots, p\}\times \{1,\cdots, q\}} A_{m_{ij}}, 
\]
where the set $\{1,\cdots, p\}\times \{1,\cdots, q\}$ is embedded into $\CC$ by $(i,j) \mapsto i+j \sqrt{-1} $. 

\vskip .2cm

Let $i\in\{0,\cdots, p-2\}$, so $\del'_i M$ is the contingency matrix of size $(p-1)\times q$ obtained from $M$ by
adding the $(i+1)$st ant $(i+2)$nd columns:
\[
(\del'_i M)_{rs} \,= \begin{cases}
m_{rs}, & \text{ if } r\leq i'
\\
m_{i_1, s} + m_{i+2, s},& \text{ if } r=i+1,
\\
m_{r-1, s}<& \text{if } r\geq i+2. 
\end{cases}
\]
We define the. {\em  $i$th  comultiplication map}
\[
\delta'_{i,M}: A_{\del'_i M} = \bigotimes_{(r,s)} A_{ (\del'_i M)_{rs}} \lra \bigotimes_{(r,s)} A_{m_{rs}} = A_M
\]
as the tensor product of the identities of the $A_{m_{rs}}$ for $r\leq i$ or $r\geq i+2$ and of the comultiplication
maps
\[
\Delta_{m_{i+1, s}, m_{i+2, s}}: A_{m_{i+1, s} + m_{i+2, s}} \lra A_{m_{i+1, s}}\otimes A_{m_{i+2, s}}.
\]
See \cite[\S2F]{KS-prob} for a more detailed treatment in the case $L=\ZZ_+$. The details in the general case
are identical. 

Similarly, let $j\in \{0,\cdots, q-2\}$ so $\del_j'' M$ is obtained from $M$ by adding the $(j+1)$st and $(j+2)$nd
columns:
\[
(\del''_j M)_{rs} \,= \, \begin{cases}
m_{rs},& \text{ if } s\leq j,
\\
m_{r, j+1} + m_{r, j+2},& \text{ if } s=j+1, 
\\
m_{r, s-1},& \text{ if } s\geq j+2. 
\end{cases}
\]
We define the {\em  $j$th multiplication map}
\[
\delta''_{j,M}: A_M = \bigotimes_{(r,s)} A_{m_{rs}} \lra \bigotimes_{(r,s)} A_{ (\del''_j M)_{rs}} = A_{\del''_j M}
\]
as the tensor product of the identities of the $A_{m_{rs}}$ for $s\leq j$ or $s\geq j+2$ and  of the multiplication maps
\[
\mu_{m_{s, j+1}, m_{s, j+2}}: A_{m_{s, j+1}} \otimes A_{m_{s, j+2}} \lra A_{m_{s, j+1} + m_{s, j+2}}.
\]
Again, the details are identical to those for $L=\ZZ_+$ explained in \cite[\S2F]{KS-prob}. 

\begin{prop}\label{prop:bialg-Janus}
Let $n\in L$. The correspondence taking:
\begin{itemize}
\item each $M\in \CM_n(L)$ to $A_M$;

\item each generating morphism $M\to \del'_i M$ in $\Gro'(\ul\CM_n(L))$ to $\delta'_{i,M}$;

\item each generating morphism $M\to\del''_j M$ in $\Gro''(\ul\CM_n(L))$ to $\delta''_{j,M}$,
\end{itemize} 
extends to a Janus sheaf $E_{A,n}$ on type $(L,n)$ with values in $\Vc$. 
\end{prop}

\noindent{\sl Proof:} Similar to the proof given in  \cite{KS-prob} to Prop. 2.7 and Cor. 3.2 for the case $L=\ZZ_+$.
To summarize the main points: the functoriality in $\Gro'(\ul\CM_n(L))$ follows from coassociativity, while
the functoriality in $\Gro''(\ul\CM_n(L))$ follows from associativity  of $A$. This shows that we get
 a datum satisfying (JS1). The axiom (JS2) follows by induction from the axiom (LB3) (compatibility of
 $\Delta$ and $\mu$) for the graded bialgebra $A$, see \S 
 \ref {sec:bi-PROB}
\ref{par:L-gr-bi} In fact, the sum over $2\times 2$ matrices in (LB3) is nothing but an instance of (JS2)
corresponding to the mixed  fork
\[
(l_1, l_2) \buildrel \del'_0\over\lra (n) \buildrel \del''_0\over\lla \begin{pmatrix} m_1 \\ m_2
\end{pmatrix}
\]
in $\Gro(\ul\CM_n(L))$. Finally, (JS3) is obvious since $A_0=\b1$ and so the tensor product for $A_M$ and $A_N$
for $M,N$ connected by an anodyne morphism, contain the same non-unit factors. \qed

 \paragraph{Janus sheaves and perverse sheaves.} The role of Janus sheaves for us stems from the following
 statement which is the second main result of this paper and will be used to prove Theorem 
  \ref{thm:Bn-perv}. 

\begin{thm}\label{thm:perv=cont}
We have an equivalence of categories $\Perv^{<\oo}(Z(\CC,L)_n, \Vc) \= 
\JS^{<\oo}_{L,n}(\Vc)$.
\end{thm}

The proof will occupy the  next section.


\section{Janus sheaves and Cousin complexes}\label{sec:janus=cousin}

In this section we prove Theorem \ref{thm:perv=cont}. by identifying Janus sheaves of type $(L,n)$ with
some explicit complexes of sheaves on $Z(\CC,L)_n$ which turn out to represent perverse sheaves in
a one-to-one way.

\paragraph{Sheaves on imaginary strata.}
Recall that $Z=Z(\CC,L)_n$ is subdivided into the imaginary strata
$Z_\bn^\Im$ labelled by ordered partitions $\bn=(\bn_1,\cdots, \bn_d)\in\CM_n^1(L)$ of $n$. 
For such $\bn$
let $\eps_\bn: Z_\bn^\Im\hra Z$ be the embedding. 
Recall that $\CM_\bn(L)\subset \CM_n(L)$ denotes the set of contingency matrices
of with exactly $d$ rows and the sum of the $i$th row equal to $n_i$.
The contingency cells $U_M$, $M\in\CM_\bn(L)$ form a
  cell decomposition of $Z_\bn^\Im$ which we still denote $\Sc^{(2)}$ or by  $\Sc^{(2)}_\bn$. 
 Combining them together, we form a disjoint union
 \[
 Z^{\{\Im\}} \,=\bigsqcup_{\bn\in\CM^1_n(L)} Z_\bn^{\Im} \buildrel\eps\over\lra Z
 \]
with the map $\eps$ being a set-theoretic bijection but not a homeomorphism. We denote
$\Sc^{(2)} = \Sc^{(2)}_{\{\Im\}}$ the cell decomposition of $Z^{\{\Im\}}$ formed by the $\Sc^{(2)}_\bn$. 

\vskip .2cm

Let $\Vc$ be a Grothendieck abelian category. let us introduce the following notation, analogous to that of 
\S \ref {sec:perv-0-cyc}\ref{par:constr-der-cat}

\vskip .2cm

$\Sh_{\Sc^{(2)}}(Z^{\{\Im\}}, \Vc)\subset \Sh(Z^{\{\Im\}}, \Vc)$: the full subcategory of $\Sc^{(2)}$-cellular
sheaves.

\vskip .2cm

$\Sh^{<\oo}(Z^{\{\Im\}}, \Vc)\subset \Sh(Z^{\{\Im\}}, \Vc)$: the full subcategory of sheaves supported on
a finite skeleton.
\[
\Sh^{<\oo}_{\Sc^{(2)}}(Z^{\{\Im\}}, \Vc) \, := \, \Sh^{<\oo}(Z^{\{\Im\}}, \Vc) \cap \Sh_{\Sc^{(2)}}(Z^{\{\Im\}},  \Vc).
\]
$D^{<\oo}_{\Sc^{(2)}}(Z^{\{\Im\}}, \Vc)$: the full subcategory in $D^b\Sh^{<\oo}(Z^{\{\Im\}}, \Vc)$ formed by
complexes whose cohomology sheaves lie in $\Sh_{\Sc^{(2)}}(Z^{\{\Im\}},  \Vc)$.

\vskip .2cm

Note that an object of $\Sh(Z^{\{\Im\}}, \Vc)$ is simply a collection $(\Fc_\bn)_{\bn\in \CM^1_n(L)}$
of sheaves $\Fc_\bn$ on the $Z_\bn^{\Im}$. 

\begin{prop}\label{prop:exit-Z-Im}
(a) The exit path category of $(Z^{\{\Im\}}, \Sc^{(2)})$ is identified with $\Gro'(\ul\CM_n(L))^\op$,
the opposite of the horizontal Grothendieck construction of $\ul\CM_n(L)$. In particular, 
$\Sh_{\Sc^{(2)}}(Z^{\{\Im\}}, \Vc)$ is identified with $\Fun(\Gro'(\ul\CM_n(L))^\op, \Vc)$. 

\vskip .2cm

(b) The category $D^{<\oo}_{\Sc^{(2)}}(Z^{\{\Im\}}, \Vc)$ is identified with 
$D^b \Fun^{<\oo}(\Gro'(\ul\CM_n(L))^\op, \Vc)$, where $\Fun^{<\oo}$ means the category of functors
equal to $0$ on almost all objects. 
\end{prop}

\noindent{\sl Proof:} (a) This is a consequence of Proposition \ref{prop:exit-S-2=gro} and the
following fact: if $M,N\in\CM_\bn(L)$ are such that $U_N\subset \ol U_M$ , then any
morphism $M\to N$ in $\Gro(\ul\CM_n(L))$ actually lies in $\Gro'(\ul\CM_n(L))$. 

(b) This is a consequence of Proposition  \ref{prop:D-constr=compl-Gro}. \qed

\vskip .2cm

As before, we use this proposition to eliminate assumption that $\Vc$ is Grothendieck, so in the
sequel $\Vc$ is an arbitrary abelian category. 

\begin{prop}\label{prop:Rqeps*F=0}
Let $\Fc \in \Sh_{\Sc^{(2)}}(Z^{\{\Im\}}, \Vc)$ be any $\Sc^{(2)}$-cellular sheaf on $Z^{\{\Im\}}$. Then $R^q\eps_*\Fc=0$ for $q>0$. 
\end{prop}

Equivalently, for any $\Sc^{(2)}$-cellular sheaf $\Fc_\bn$ on $Z_\bn^\Im$ we have 
$R^q\eps_{n*}\Fc_\bn=0$ for
$q>0$. Here, as before, $\eps_\bn: Z_\bn^\Im\to Z$ denotes the embedding. 

\vskip .2cm

\noindent{\sl Proof:} Let
$e: \Gro'(\ul\CM_n(L))^\op\to \Gro(\ul\CM_n(L))^\op$ be the embedding of the
horizontal Grothendieck construction to the full one (a bijection on objects).  In the identifications of
Propositions \ref{prop:exit-S-2=gro}, \ref{prop:D-constr=compl-Gro} and \ref{prop:exit-Z-Im}, the
pullback functor
\[
\eps^*: \Sh_{\Sc^{(2)}}(Z,\Vc) \lra \Sh_{\Sc^{(2)}}(Z^{\{\Im\}}, \Vc)
\]
corresponds to the functor given by the pullback (composition with $e$):
\[
e^*: \Fun(\Gro(\ul\CM_n(L))^\op, \Vc) \lra \Fun(\Gro'(\ul\CM_n(L))^\op, \Vc).
\]
Therefore $R^q\eps_*$ is given by the $q$th derived right Kan extension functor
\[
R^q e_*:  \Fun(\Gro'(\ul\CM_n(L))^\op, \Vc) \lra  \Fun(\Gro(\ul\CM_n(L))^\op, \Vc).
\]
Explicitly, the value of $R^qe_* F$ on $M\in\CM_n(L)$ is given by the
$q$th derived functor of the projective limit functor
\be\label{eq:Rqe*F}
(R^q e_* F)(M) \,=\, R^q \varprojlim\nolimits_{\{ M\to e(N)\}\in M/e} F(M),
\ee
where $M/e$ is the comma category whose objects are pairs formed by $N\in \CM_n(L)$
and a morphism $M\to e(N)$ in $\Gro(\ul\CM_n(L))^\op$ or, in plainer terms,
morphisms $\phi$ from $e(N)=N$ to $M$ in $\Gro(\ul\CM_n(L))$. In the latter interpretation,
morphisms from $(N_1, \phi_1: N_1\to M)$ to $(N_2, \phi_2: N_2\to M)$  are morphisms $\psi: N_2\to N_1$
in $\Gro'(\ul\CM_n(L))$ such that $\phi_1\psi=\phi_2$. 

Let us use the notation $N\buildrel '\over\to M$ resp. $N\buildrel ''\over\to M$ to indicate a morphism
$N\to M$ in $\Gro'(\ul\CM_n(L))$ resp. $\Gro''(\ul\CM_n(L))$. Any object of $M/e$ figuring in 
\eqref{eq:Rqe*F} can be uniquely represented as a composition $N\buildrel '\over\to P\buildrel ''\over\to M$.
This means that the non-derived ($q=0$) instance of \eqref{eq:Rqe*F}  can be represented as a double
projective limit
\be\label{eq:p-->M<--N}
\varprojlim\nolimits_{\{P\buildrel ''\over\to M\}} \varprojlim\nolimits_{\{N\buildrel '\over\to P\}} F(N).
\ee
To prove that the derived functor in \eqref{eq:Rqe*F} is zero for $q>0$, it suffices to show that each
limit in \eqref{eq:p-->M<--N} has higher derived functors equal to $0$. Now, the inner limit over 
$N\buildrel '\over\to P$ is over a conical diagram, since $F$ is a functor on $\Gro'(\ul\CM_n(L))$.
Therefore $\varprojlim_{\{N\buildrel '\over\to P\}} F(N) = F(P)$, while higher derived functors of this $\varprojlim$
vanish. So it remains to show that $R^q\varprojlim_{\{ P\buildrel ''\over\to M\}} F(P)=0$
for $q>0$. But the comma category formed by arrows $P\buildrel ''\over\to M$ is discrete.
Indeed, morphisms between such arrows are given by morphisms $P_1\to P_2$ in $\Gro'(\ul\CM_n(L))$
making a triangle commutative. So the $R^q\varprojlim$ over this discrete category vanish for
$q>0$.  \qed

\vskip .2cm

Since the limit over a discrete category is the direct product, we arrive at the following complement to
Proposition \ref{prop:Rqeps*F=0}. 

\begin{prop}\label{prop:e*F=prod}
Let $\Fc\in\Sh_{\Sc^{(2)}}(Z^{\{\Im\}}, \Vc)$ and $F: \Gro'(\ul\CM_n(L))^\op\to\Vc$ be the corresponding functor
so that $F(M)$ is the stalk of $\Fc$ at $U_M$.  Then $\eps_*\Fc$ corresponds to the functor $e_*F: 
\Gro(\ul\CM_n(L))^\op\to\Vc$ given on objects by
\[
(e_*F)(M) \,= \prod_{N\buildrel ''\over\to M} F(N). 
\]
The values of $e_*F$ on morphisms from $\Gro'(\ul\CM_n(L))$ are
induced by those of $F$, while the value on a morphism $g: M_1\to M_2$ from $\Gro''(\ul\CM_n(L))$
is the morphism
\[
(e_*F)(g): \prod_{\phi_2: N_2\buildrel ''\over\to M_2} F(N_2) \lra \prod_{\phi_1: N_1\buildrel ''\over\to M_1}
F(N_1) 
\]
with matrix elements
\[
(e_*F)(g)_{\phi_1, \phi_2} = \begin{cases}
\Id_{F(N_2)}, & \text{ if } N_1=N_2 \text { and } \phi_2=g\phi_1,
\\
0, & \text{ otherwise.}
\end{cases}
\]
\qed
\end{prop}

\paragraph{Cousin-type complexes.} 
\begin{Defi}
A {\em Cousin-type complex} on $Z$ (with values in $\Vc$)  is a datum consisting of
$\Sc^{(2)}$-cellular sheaves $\Fc_\bn$ on $Z_\bn^\Im$, $\bn\in \CM^1_n(L)$
and a complex of sheaves of the form
\[
\Fc^\bullet =\biggl\{ \cdots \buildrel d\over \to
\bigoplus_{l(\bn)=r} \eps_{\bn*} \Fc_\bn
\buildrel d\over \to  \cdots\buildrel d\over \to 
\bigoplus_{\bn=(n_1, n_2)} \eps_{(n_1, n_2)*} \Fc_{(n_1, n_2)}  \buildrel d\over\to\eps_{(n)*} \Fc_{(n)}
\biggr\}
\]
graded so that the sum over $l(\bn)=r$ is in degree $-r$ (so the complex
terminates in degree $-1$). Here $\eps_{\bn *}$ is understood as
the non-derived direct image functor, i.e., as $R^0\eps_{\bn *}$. 
\end{Defi}

\begin{prop}\label{prop:Cousin=Janus}
The category formed by Cousin-type complexes (with morphisms being
morphisms of complexes of sheaves) is equivalent to the 
category of  Janus presheaves of type $(L,n)$. 
\end{prop}

\noindent{\sl Proof:} Let $E=(E(M), \delta'_{f}, \delta''_{g})$ be a  Janus
 presheaf  of type $(L,n)$. For $\bn\in \CM_n(L)$ we form
the cellular sheaf $\Fc_\bn$ on $Z_\bn^\Im$ via Proposition
\ref {prop:exit-Z-Im}(a). That is, the stalk of $\Fc_\bn$
at the cell $U_M, M\in\CM_\bn(L)$, is $E(M)$ and the generalization maps
are given by the $\delta'_f$.

Next, let $l(\bn)=r$ and $N\in\CM_\bn(L)$. For any $i=0,\cdots, d-2$
we have the ordered partition $\del_i\bn$ obtained by adding the $(i+1)$st
and $(i+2)$nd entries of $\bn$ and the contingency matrix 
$\del''_iN\in \CM_{\del_i\bn}(L)$ obtained by adding he $(i+1)$st
and $(i+2)$nd rows of $N$. Thus $N\to \del''_iN $, an elementary
 morphism in $\Gro''(\ul\CM_n(L))$. Let 
 \[
\delta''_{i,N} = \delta''_{N\to\del''_n N}: E(N) \lra E(\del''_i N).
\]
Let $\ul\CM_\bn(L)\subset\ul\CM_n(L)$ be the (single, not double) augmented
semi-simplicial set formed by $N\in\CM_\bn(L)$ and the face morphisms $\delta'_i$ and their
compositions. Let $\Gro'(\ul\CM_\bn(L))\subset\Gro'(\ul\CM_n(L))$ be the corresponding
Grothendieck construction. Taken for all $\bn\in\CM^1_n(L)$, the $\Gro'(\ul\CM_\bn(L))$
form the connected components of $\Gro'(\ul\CM_n(L))$. 

Let $e_\bn: \Gro'(\ul\CM_\bn(L))\to \Gro(\ul\CM_n(L))$ be the embedding functor.
The functor $E_\bn: \Gro'(\ul\CM_\bn(L))^\op\to\Vc$ representing the cellular sheaf $\Fc_\bn$ is
the restriction of $(E,\delta')$ to $\Gro'(\ul\CM_\bn(L)$. So $\eps_{\bn *} \Fc_\bn$ is represented by the right
Kan extension $e_{\bn *} E_\bn: \Gro(\ul\CM_n(L))^\op\to\Vc$ whose values on objects and morphisms are
given by Proposition \ref{prop:e*F=prod}. Let us define a morphism of functors
\[
D_{i, \bn}: e_{\bn *} E_\bn \lra (e_{\del_i\bn })_* E_{\del_i\bn}
\]
by defining the functor
\[
D^*_{i, \bn}: e^*_{\del_i\bn} e_{\bn *} E_\bn \lra E_{\del_i\bn}
\]
corresponding to it by adjunction. The value of $D^*_{i,\bn}$ on an object $M\in\CM_{\del_i\bn}(L)$ is the morphism
\[
(D^*_{i,\bn})_M: \prod_{ h: N\buildrel ''\over\to M \atop N\in\CM_\bn(L)} E(N) \lra E(M), 
\]
whose component $(D^*_{i,\bn})_{M,h}: E(N)\to E(M)$ is $0$ unless $M=\del''_i N$ and $g=\del_i''$, in which case it is
$\delta''_{i,N}$.

Let us show that $D^*_{i, \bn}$ is indeed a morphism of functors, i.e., that it commutes with
the action of each morphism $f: M_1\to M_2$ from $\Gro'(\ul\CM_{\del_i \bn}(L))$. In other words, we prove
that the diagram
\[
\xymatrix{
\prod_{h_1: N_1\buildrel ''\over\to M_1 \atop N_1\in\CM_\bn(L)} E(N_1) \ar[rrr]^{(D^*_{i,\bn})_{M_1} } &&& E(M_1)
\\
&&&
\\
\prod_{h_2: N_2\buildrel ''\over\to M_2 \atop N_2\in\CM_\bn(L)} E(N_2)
\ar[uu]^{(e^*_{\del_i\bn} e_{\bn *} E_\bn)(f)}
\ar[rrr]_{(D^*_{i,\bn})_{M_2} } &&& E(M_2)
\ar[uu]_{\delta'_f}
}
\]
is commutative. Let us check it on a factor $E(N_2)$ in the lower left corner, corresponding to any given $h_2$. 
The composition $\delta'_f (D^*_{i,\bn})_{M_2}$ takes $E(N_2)$ to $E(\del''_i N_2)$ by $\delta''_{i,N}$, if
$M_2=\del''_i N_2$ (and to $0$ otherwise). After this, application of $\delta'_f$ gives $\delta'_f\delta''_{i,N}$.
This gives a mixed fork $P=\{N_2\buildrel \del''_i\over\rightarrow M_2\buildrel f\over\leftarrow M_1\}$ in the
sense of \S \ref{sec:X=C}\ref{par:forks} 
Now, the other path in the diagram corresponds to the sum over $\Sup(P)$ as in the condition (JS2). 
So this condition implies that $D^*_{i,\bn}$ and thus $D_{i,\bn}$ is indeed a morphism of functors.

\vskip .2cm

Let $d_{i,\bn}: \eps_{\bn *}\Fc_\bn \to (\eps_{\del_i\bn})_*\Fc_{\del_i\bn}$ be the morphism of sheaves corresponding
to $D_{i,\bn}$ and let
\[
d\,=\sum_{l(\bn)=r} \sum_{i=0}^{r-2} (-1)^i d_{i,\bn}: \bigoplus_{l(\bn)=r} \eps_{\bn *} \Fc_\bn
\lra\bigoplus_{l(\bn)=r-1} \eps_{\bn *} \Fc_\bn. 
\]
We then have $d^2=0$ because $E$ is a functor $\Gro''(\ul\CM_n(L))\to \Vc$. So we have a Cousin-type
complex associated to $E$. Denote it $\Fc^\bullet(E)$. 

\vskip .2cm

Conversely, let $\Fc^\bullet$ be a Cousin-type complex. We recover the Janus sheaf $E$ from $\Fc^\bullet$ by 
reversing the above considerations. That is, we define $E(M)$ as the stalk of $\Fc_\bn$ at $U_M$, where $\bn$ is
such that $M\in\CM_\bn(L)$. The data of the $\delta'_f$, i.e., the 
contravariant functoriality of $E$ on $\Gro'(\ul\CM_n(L))$, are given by the sheaf structure of the $\Fc_\bn$.
The differential $d$ gives the covariant functoriality on $\Gro''(\ul\CM_n(L))$ for elementary morphisms
$\del''_i$. The condition $d^2=0$ translates into simplicial identities for the morphisms $\delta''_{i,N}$ associated
by $E$ to the $\del''_o$, i.e., they define $E$ as a functor on all morphisms from $\Gro''(\ul\CM_n(L))$.
The condition that $d$ is a morphism of sheaves, translates into (JS2). Proposition \ref{prop:Cousin=Janus} is
proved. 

\paragraph{The Cousin-type complex as the horizontal Verdier dual.} The Cousin-type complex $\Fc^\bullet =
\Fc^\bullet(E)$ associated to a Janus presheaf $E$, is a complex of $\Sc^{(2)}$-cellular sheaves. As such, it
is represented by a functor which we denote $\DD''E: \Gro(\ul\CM_n(L))^\op\to C^b(\Vc)$. 
In particular, $\DD''E$ is contravariant on both subcategories $\Gro'(\ul\CM_n(L))$ and $\Gro''(\ul\CM_n(L))$,
while $E$ itself is contravariant on $\Gro'$ and covariant on $\Gro''$ (so it is a ``sheaf in the horizontal
direction'' and a ``cosheaf in the vertical direction''). The following result describes $\DD''E$ as a kind of
relative Verdier dual to $E$ in the vertical direction so that the cosheaf--type behavior in that direction
is transformed into the sheaf-type behavior. 

\vskip .2cm

Given $M\in\CM_n(L)$ of size $p\times q$, the cell $U_M$ is identified with $\RR^p_<\times \RR^q_<$,
see \S \ref{sec:cont-cells}\ref{par:S2} So the orientation space $\OR_M$ splits naturally into tensor product
$\OR_M=\OR'_M\otimes\OR''_M$, where $\OR'_M$ is the orientation space of $\RR^p_<$ and $\OR''_M$ is that
of $\RR^q_<$. We also write $d'_M=p$, $d''_M=q$, so $d_M=\dim_\RR U_M=d'_M + d''_M$. 
 Further, we will write $f: N\mathop{\buildrel  ''\over\to}_r M$ to signify that $f$ is a morphism in
 $\Gro''(\ul\CM(L))$ and $d_N=d_M+r$ (or, equivalently, $d''_N=d''_M+r$). 
 We also write $\OR''_{N/M}=\OR''_M \otimes (\OR''_M)^{-1}$. 
 
 \begin{prop}\label{prop:cousin=rel-D}
  We have
 \[
 (\DD''E)(M) \,= \, \biggl\{ \cdots \buildrel d\over \to \bigoplus_{g_2: M_2 \mathop{\buildrel ''\over\to}_2 M} E(M_2)
 \otimes\OR''_{M_2} \buildrel d\over \to 
  \bigoplus_{g_1: M_1 \mathop{\buildrel ''\over\to}_1 M} E(M_1)
 \otimes\OR''_{M_1}  \buildrel d\over \to 
 E(M)\otimes\OR''_M\biggr\},
 \]
 where:
 \begin{itemize}
 \item The rightmost term $E(M)\otimes \OR''_M$ is in degree $(-d''_M)$.
 
 \item The matrix element
 \[
 d_{g_r, g_{r-1}}: E(M_r) \otimes\OR''_{M_r} \lra E(M_{r-1})\otimes\OR''_{M_{r-1}}
 \]
 is equal to 
 \[
 \sum_{\phi: M_r \buildrel ''\over\to M_{r-1}\atop g_{r-1}\phi=g_r} \eps_\phi\otimes \delta''_\phi,
 \]
 where, in its turn: 
 \begin{itemize}
 \item $\delta''_\phi$ is the value of $E$ on $\phi$; 
 
 \item $\eps_\phi: \OR''_{M_r} \to \OR''_{M_{r-1}}$ is the identification (coorientation) induced by the embedding
 of $U_{M_{r-1}}$ as a codimension $1$ cell in the closure of $U_{M_r}$. 
 \end{itemize}
 
 \end{itemize}
 \end{prop}
 
 \noindent{\sl Proof:} Direct comparison with the construction of $\Fc^\bullet(E)$. We need to note that any $\phi$
 in the sum for $d_{g_r, g_{r-1}}$ has the form $\del''_i$ for some $i$. Furtherm the space $\OR''_{M_r}$ 
 can be trivialized using the identification $U_{M_r} \= \RR^{d'_{M_r}}_< \times \RR^{d''_{M_r}}_<$ and the
 canonical orientation of $\RR^{d''_{M_r}}$. After this, the identification $\eps_\phi$ is transformed into the
 standard simplicial sign $(-1)^i$, thus matching with the differential constructed in Proposition 
 \ref{prop:Cousin=Janus}. \qed
 
 \paragraph{Cousin-type complexes and (absolute) Verdier duality.} Let $\Vc$ be an abelian category
 and $\Vc^\op$ its opposite category. For convenience we will denote by $V^*$ the object of $V^\op$
 corresponding to $V\in\Ob(\Vc)$ and by $f^*: W^*\to V^*$ the morphism in $\Vc^\op$ corresponding
 to a morphism $f: V\to W$ in $\Vc$. 
 
 \vskip .2cm
 
 Let $E=(E(M), \delta'_{f,E}, \delta''_{g,E})$ be a Janus presheaf of type $(L,n)$ with values in $\Vc$. We define a
 Janus presheaf $E^\tau = (E^\tau(M), \delta'_{f, E^\tau}, \delta''_{g, E^\tau})$ with
 values in $\Vc^\op$ as follows:
 \begin{itemize}
 \item $E^\tau(M) = E(M^t)^*$, where $M^t = \|m_{ji}\|$ is the transpose of $M=\{ m_{ij}\|$. 
 Note that a morphism $f: M\to N$ in $\Gro'(\ul\CM_n(L))$ resp. $\Gro''(\ul\CM_n(L))$
 gives a morphism $f^t: M^\tau\to N^\tau$ in $\Gro''(\ul\CM_n(L))$ resp. $\Gro'(\ul\CM_n(L))$. 
 
 \item For $f: M\to N$ a morphism in $\Gro'(\ul\CM_n(L))$ we put $\delta'_{f, E^\tau}=(\delta''_{f^t, E})^*$. 
 
 \item Similarly, for $g: M\to N$ a morphism in $\Gro''(\ul\CM_n(L))$ we put
 $\delta''_{g, E^\tau}= (\delta'_{g^t, E})^*$. 
 \end{itemize}
 It is clear that $E^\tau$ is again a Janus presheaf. Further, if $E$ is a Janus sheaf, then so is $E^\tau$. 
 
 \vskip .2cm
 
 As before, we denote by $\tau$ the involution on $Z=Z(\CC,L)_n$ induced by permutation of coordinates
 $(x+iy)\mapsto (y+ix)$ on $\CC$. Clearly $\tau$ preserves the cell dcomposition $\Sc^{(2)}$,
 with $\tau(U_M)=U_{M^t}$. Recall that we have Verdier duality which is an anti-equivalence
 \[
 \DD: D^{<\oo}_{\Sc^{(2)}}(Z,\Vc) \lra  D^{<\oo}_{\Sc^{(2)}}(Z,\Vc^\op). 
 \]
 
 \begin{prop}\label{prop:D=tau}
 Let $E$ be a Janus presheaf of type $(L,n)$ with values in $\Vc$ and with finite support. Let $\Fc^\bullet(E)$
 be the corresponding Cousin-type complex. Then we have an isomorphism in the derived category:
 \[
 \DD(\Fc^\bullet(E)) \,\=\, \tau^* \Fc^\bullet(E^\tau). 
 \]
 \end{prop}
 
 \noindent{\sl Proof:} Let us denote by the same letter $\DD$ the effect of Verdier duality on functors from
 the category of exit paths, i.e., the horizontal arrow in the following diagram:
 \[
 \xymatrix{
 D^b\Fun^{<\oo}(\Gro(\ul\CM_n(L))^\op, \Vc) \ar[r]^\DD 
 \ar[d]^{\DD''} & D^b\Fun^{<\oo}(\Gro(\ul\CM_n(L)), \Vc) 
 \\
 D^b \JPS^{<\oo}(\Vc) \ar[ur]_{\DD'}
  & 
 }
 \]
The description of $\DD$ given in Proposition \ref{prop:verdier-gro} implies that it can be decomposed as
the composition $\DD=\DD'\circ\DD''$ of the horizontal and vertical Verdier dualities. Here $\DD''$ changes
contravriant functoriality on $\Gro''(\ul \CM_n(L))$ to covariant one, thus producing a complex of Janus
sheaves while $\DD'$ changes contravariant functoriality on $\Gro'(\ul\CM_n(L))$ to covariant one,
thus producing a complex of cosheaves with values in $\Vc$, i.e., of sheaves with values in $\Vc^\op$. 
More  precisely, $\DD'$ is the analog of the functor described on objects in Proposition \ref{prop:cousin=rel-D}
 but with $\Gro'(\ul\CM_n(L))$ playing the role of $\Gro''(\ul\CM_n(L))$, while $\DD''$ involves $\Gro''$. 
 
 This means that the complex of $\Vc^\op$-valued sheaves associated to $\DD'E$, $E\in\JPS^{<\oo}(\Vc)$ is 
 $\tau^* \Fc^\bullet(E^\tau)$, as the rolws of the real and imaginary strata are interechanged. 
 Our statement reduces to the claim that vertical Verdier duality is involutive, i.e., the composition
 \[
 D^b \JPS^{<\oo}(\Vc)\buildrel \DD''\over\lra D^b\Fun^{<\oo}(\Gro(\ul\CM_n(L))^\op, \Vc) 
 \buildrel \DD''\over\lra  D^b \JPS^{<\oo}(\Vc)
 \]
 is isomorphic to the identity. This is proved by the same formal argument as the involutivity of the ordinary
 (absolute) Verdier duality at the cellular level. More precisely, let $E\in D^b \JPS^{<\oo}(\Vc)$. Then, $\DD'' E$
 being contravariant in $\Gro''(\ul\CM_n(L))$,  the value $(\DD'' \DD''E)(M)$ is the total complex of
 \[
 (\DD''E)(M) \otimes\OR''_M \to \bigoplus_{N_1\mathop{\buildrel ''\over\to}_1 M} (\DD''E)(N_1) \otimes\OR''_{N_1}
 \to 
 \bigoplus_{N_2\mathop{\buildrel ''\over\to}_2 M} (\DD''E)(N_2) \otimes\OR''_{N_2}\to\cdots
\]
 the horizontal grading starting at $d''_M$. For uniformity, the leftmost term corresponds to the unique arrow
 (identity) $M= N_0 \mathop{\buildrel ''\over\to}_0 M$. Substituting here the expression for
 each $(\DD''E)(N_r)$ from Proposition \ref{prop:cousin=rel-D}, we represent $(\DD''\DD''E)(M)$ as the total complex
 of the double complex
 \[
 \xymatrix{
 E(M) & & 
 \\ 
 \bigoplus\limits_{M_1\mathop{\to}_1 N_0 = M} E(M_1)\otimes \OR''_{M_1/N_0} \ar[r]  \ar[u]& \bigoplus\limits_{M_1=N_1
 \mathop{\to}_1 M} E(M_1) &
 \\
 \bigoplus\limits_{M_2\mathop{\to}_2 N_0 = M} E(M_2)\otimes \OR''_{M_2/N_0} \ar[r] \ar[u] &
 \bigoplus\limits_{M_2 \mathop{\to}_1 M_1\mathop{\to}_1 N_0 = M} E(M_2)\otimes \OR''_{M_2/N_1} 
 \ar[r]\ar[u] & \bigoplus\limits_{M_2=N_2
 \mathop{\to}_2 M} E(M_2) 
 \\
 \vdots \ar[u] & \vdots \ar[u]&\vdots\ar[u]
 }
 \]
 where $E(M)$ is in degree $(0,0)$ and terms that are $0$, are not depicted. All arrows indexing direct sums,
 are in $\Gro''(\ul\CM_n(L))$, so we use the notation $\mathop{\to}_p$ instead of 
 $\mathop{\buildrel ''\over \to}_p$.
 The vertical differential involves the covariant functoriality of $E$ on $\Gro''(\ul\CM_n(L))$ while the horizontal
 differential is purely combinatorial. More precisely, the $(-p)$th term is a direct sum
 $\bigoplus_{f: M_p \mathop{\buildrel '' \over\to}_p M} C^\bullet_f \otimes_\k E(M_p)$, where $f$
 runs over all morphisms $M_p\to M$ in $\Gro''(\ul\CM_n(L))$ of relative dimension $p$, and $C^\bullet_f$
 is the complex of $\k$-vector spaces with
 \[
 C^r_f \,= \bigoplus_{M_p \mathop{\buildrel h\over\to}_{p-r} N_r \mathop{\buildrel g\over\to}_r M \atop gh=f}
 \OR''_{M_p/N_r},
 \]
 the sum being over all factorizations of $f$ into the compostion $gh$ where $g$ has relative dimension $r$ and $h$
 has relative dimension $p-r$. 
 
 We now proce that each row of our double complex except the top one (consisting of $E(M)$ alone) is exact.
 For this it suffices to rpove that each $C^\bullet_f$ is exact, if $p>0$.  To see this, note that $f$, being
 a morphism in the Grothendieck construction, corresponds to a morphism $\phi: [d''_M+p-2] \to [d''_M-2]$
 in $\Delta_+^\op$ (here $d''_M$ is the number of columns of $M$). In particular, $M=\del''_\phi (M_p)$
 with respect to the double augmented semi-simplicial structure on $\ul\CM_n(L)$. This means that factorizations
 of $f$ are in bijection with factorizations of $\phi$ in $\Delta_+^\op$. Let $m=d''_M+p-2$. Morphisms
 in $\Delta_+^\op$ with source $M$ are in bijection with faces $\Gamma\subset\Delta^m$ of the standard
 $m$-simplex.  Factorizations of $\phi$ correspond to intermediate faces $\Gamma\subset\Gamma'\subset \Delta^m$,
 where $\gamma$ is the face corresponding to $\phi$. So we find that $C^\bullet_f$ is the augmented
 cellular cochain complex of the link of $\Gamma$ in $\Delta^m$. Since this link is itself a simplex, $C^\bullet_f$
 is exact and so $(\DD''\DD'' E)(M)$ is quasi-isomorphic to $E(M)$. 
 This finishes the proof of Proposition \ref{prop:D=tau}.

 \paragraph{Fron a Janus sheaf to a perverse sheaf.} Let $E$ be a Janus sheaf with finite support.
 The associated Cousin-type complex $\Fc^\bullet(E)\in D^{<\oo}_{\Sc^{(2)}}(Z,\Vc)$ will be called the
 {\em Cousin complex} of $E$. 
 
 \begin{prop}\label{prop:cous-perv}
$\Fc^\bullet(E)\in\Perv^{<\oo}(Z,E)$. 
 \end{prop}
 
 \noindent{\sl Proof:} This statement consists of three parts:
 \begin{itemize}
 \item[(1)] $\Fc^\bullet(E)$ is an $\Sc^{(0)}$-constructible complex, i.e., the $\ul H^r\Fc^\bullet(E)$ are $\Sc^{(0)}$-constructible sheaves.
 
 \item[(2)] $\Fc^\bullet(E)$ satisfies (P$^-$),  see \S \ref{sec:perv-0-cyc}\ref{par:perv-t-str}
 
 \item[(3)]  $\Fc^\bullet(E)$ satisfies (P$^+$). 
 \end{itemize}
 
 To prove (1) it suffices, by Proposition \ref{prop:S2=S1+tS1},  that $\Fc^\bullet(E)$ is both $\Sc^{(1)}$-constructible
 and $\tau \Sc^{(1)}$-constructible. We note that each term $\Fc^{-r}(E)=\bigoplus_{l(\bn)=r} \eps_{\bn *} \Fc_\bn$
 is $S^{(1)}$-constructible (and so the cohomology sheaves of $\Fc^\bullet(E)$ are $\Sc^{(1)}$-constructible
 as well). This is because each sheaf $\Fc_\bn\in\Sh_{\Sc^{(2)}}(Z_\bn^\Im)$ is $\Sc^{(1)}$-constructible. This last 
 statement is a consequence of the axiom (JS3) for $E$ which says that anodyne morphisms in
 $\Gro''(\ul\CM_n(L))$ (corresponding to inclusions of $\Sc^{(2)}$-cells within one $\Sc^{(1)}$-cell)
 are taken to isomorphisms. 
 
 To prove that $\Fc^\bullet(E)$ is $\tau\Sc^{(1)}$-constructible, , it suffices to show that $\DD\Fc^\bullet(E)$,
 its Verdier dual, is $\tau\Sc^{(1)}$-constructible. But by Proposition \ref{prop:D=tau}, 
 $\DD\Fc^\bullet(E)\= \tau^* \Fc ^\bullet(E^\tau)$. Since $E^\tau$ is a Janus sheaf, $\Fc^\bullet(E^\tau)$
 is $\Sc^{(1)}$-constructible by the above, and so $\tau\Fc^\bullet(E^\tau)$ is $\tau\Sc^{(1)}$-constructible.
 This proves (1).
 
 \vskip .2cm
 
 To prove (2), we need to show that for each $r\geq 1$ the sheaf $\ul H^{-r}\Fc^\bullet(E)$ (already known
 to be $\Sc^{(0)}$-constructible) haas support on the union of $\Sc^{(0)}$-strata of complex dimension $\leq r$. 
 The strata $Z_\alpha$ of $\Sc^{(0)}$ are labelled by unordered partitions $\alpha$ of $n$. 
 If $\bn\in\CM_n^1(L)$ is an ordered partition of $n$, we denote $\alpha(\bn)$ the corresponding
 unordered partition.  Now, if $Z_\alpha\subset \supp\, \ul H^{-r}\Fc^\bullet(E)$
 then $Z_\alpha \subset\supp \, \Fc^{-r}(E) = \bigoplus_{l(\bn)=r} Z_\bn^\Im$. 
 By irreducibility of $Z_\alpha$ this implies that $Z_\alpha\subset Z_\bn^\Im$ for some $\bn\in\CM_n^1(L)$
 with $l(\bn)=r$. But such an ainclusion is possible only if $\alpha$ is a coarsening of $\alpha(\bn)$,
 so $\alpha$ consists of $\leq r$ parts and therefore $\dim_\CC Z_\alpha \leq r$. This proves (2).
 
 \vskip .2cm
 
 Finally, the condition  (P$^+$) is Verdier dual to  (P$^-$). So to show that $\Fc^\bullet(E)$ satisfies  (P$^+$),
 it is enough to show that $\DD\Fc^\bullet(E)$ satisfies  (P$^-$). But $\DD\Fc^\bullet(E) \= \tau^*\Fc^\bullet(E^\tau)$,
 and $\Fc^\bullet(E^\tau)$ satisfies  (P$^-$) by the above. Since $\tau: Z\to Z$ preserves the strata of
 $\Sc^{(0)}$, we see that $\tau^*\Fc^\bullet(E^\tau)$ satisfies  (P$^-$) as well. This shows (3) and
 finishes the proof of Proposition \ref{prop:cous-perv}. 
 
 \paragraph{From a perverse sheaf to a Janus sheaf.} Proposition  \ref{prop:cous-perv} defined a functor
 \be\label{eq:functor-c}
 c: \JS^{<\oo}_{L,n}(Z,\Vc) \lra\Perv^{<\oo}(Z,\Vc), \quad E\mapsto \Fc^\bullet(E). 
 \ee
 In order to prove Theorem \ref{thm:perv=cont}, we construct the inverse functor.
 
 \begin{prop}\label{prop:eps_n^!}
 Let $\Gc^\bullet\in\Perv^{<\oo}(Z,\Vc)$ and $\bn\in\CM_n^1(L)$ be an ordered partition of length $r$.
 Then the complex $\eps_\bn^!\Gc^\bullet$ on $Z_\bn^\Im$ is quasi-isomorphic to a single sheaf $\Fc_\bn$
 in degree $(-r)$. 
 \end{prop}
 
 \noindent {\sl Proof:} As in \cite[\S4.1]{KS-hW}, we deduce the statement from a result
 of \cite{KS-arr} about hyperplane arrangements which we recall for convenience of the reader. 
 
 \vskip .2cm
 
 let $\Hc$ be an arrangement of real hyperplanes (containing $0$)  in a Euclidean space $\RR^D$. It subdivides
 $\RR^D$ into a disjoint union of {\em faces} which are locally closed conic subsets. For a face $C$ we denote
 $\rho_C: \RR^D+C\sqrt{-1}  \hra \CC^D$ the emdedding of the corresponding tube domain. The complexification of the
 hyperlplanes from $\Hc$ define a complex stratification of $\CC^D$ into generic points of complex flats.
 Recall that a {\em complex flat} is the intersection of the complexifications of  some hyperplanes from $\Hc$.
 We denote $\Perv_\Hc(\CC^D)$ the category of perverse sheaves on $\CC^D$, constructible with 
 respect to this stratification. 
 
 \begin{prop}\label{prop:7.15}
 Let $C$ be a face of $\Hc$ of real dimension $r$. Then for any $\Gc^\bullet\in \Perv_\Hc(\CC^D)$
 we have $\ul H^{\neq  -r} \rho_C^!\Gc^\bullet =0$.
 \end{prop}
 
 \noindent{\sl Proof:} This is a consequence of \cite[Prop. 3.10]{KS-arr} (reduction to the case $r=0$)
 and  \cite[Prop. 4.9(a)]{KS-arr} (case $r=0$). For another
  transparent proof of the case $r=0$ see \cite[Prop. 1.3.3]{KS-shuffle}.
 These references deal with the case $\Vc=\Vect_\k$. The case of an arbitrary abelian category $\Vc$
 presents no difference. \qed 
 
 \vskip .2cm
 
 To reduce Proposition \ref{prop:eps_n^!} from \ref{prop:7.15}, we choose a set $(\xi_k)_{k\in K}$
 of monoid generators of $L$ and for each $l\geq 0$ consider the map
 $p_l: (\CC\times K)^l\to Z(\CC,L)$ described in \eqref{eq:map-p-l}. 
 The product $(\CC\times K)^l=\CC^l\times K^l$ is the disjoint union of $K^l$ copies of $\CC^l$. Denote
 the copy corresponding to $(k_1,\cdots, k_l)\in K^l$ by $\CC^l_{(k_1,\cdots, k_l)}$. The restriction of $p_l$
 to $\CC^l_{(k_1,\cdots, k_l)}$ is the map
 \[
 p_{(k_1,\cdots, k_l)}: \CC^l_{(k_1,\cdots, k_l)}\lra Z(\CC,L), \quad (z_1,\cdots, z_l) \mapsto \sum \xi_{k_i}\cdot z_i. 
\]
  Its image lies in $Z(\CC,L)_n$, where $n=\sum\xi_{k_i}$. The stratification of $(\CC\times K)^l$ by multiplicities,
  when restricted to $\CC^l_{(k_1,\cdots, k_l)}$, gives the stratification given by the arrangement $\Hc'$
  of hyperplanes $H_{ij}=\{z_i=z_j\}$ for all pairs $(i<j)$ such that $k_i=k_j$. Let $\Hc\supset\Hc'$ be the diagonal
  arrangement formed by all $H_{ij}$, $i<j$. 
  
  Let us now choose $l$ so that $\Im(p_l)\supset \Supp(\Gc)$. Since $\Gc$ is a pervese sheaf on $Z(\CC,L)_n$,
  it follows tnat there is $(k_1,\cdots, k_l)$ such that $\Supp(\Gc)\subset \Im(p_{(k_1,\cdots, k_l)})$.
  Thus 
  \[
  p^*\Gc\, \in\,  \Perv_{\Hc'} (\CC^l_{(k_1,\cdots, k_l)}) \, \subset\,   \Perv_{\Hc} (\CC^l_{(k_1,\cdots, k_l)}).
  \]
  It further follows that 
  \[
  \Supp(\eps_\bn^! \Gc) \,\subset \, Z_\bn^\Im \cap \Im (p_{(k_1,\cdots, k_l)}). 
  \]
  let $F$ be the fiber product fitting into the Cartesian diagram
  \[
  \xymatrix{
  F\ar[rr]^\pi 
  \ar[d]_\rho
  && Z_\bn^\Im \ar[d]^{\eps_\bn}
  \\
  \CC^l_{(k_1,\cdots, k_l)} \ar[rr]_{p_{(k_1,\cdots, k_l)}} && Z(\CC,L)_n.
  }
  \]
 To show that the sheaf  $\ul H^j(\eps_\bn^!)\Gc)$ vanishes when $j\neq -r$, it suffices to prove that the sheaf
 $\pi^* \ul H^j (\eps_\bn^!\Gc)= \ul H^j (\pi^* \eps_\bn^!\Gc)$ vanishes for $j\neq -r$. 
 But $p_{(k_1,\cdots, k_l)}$ is finite and so, by Cartesianity, is $\pi$. Therefore $\pi^*=\pi^!$ and
 $p_{(k_1,\cdots, k_l)}^* = p_{(k_1,\cdots, k_l)}^!$. Furthermore, $\Gc':= p_{(k_1,\cdots, k_l)}^*\Gc$ lies in
 $\Perv_\Hc(\CC^l_{(k_1,\cdots, k_l)})$, so it is enough to show that
 \[
\pi^* \eps_\bn^! \Gc\,=\, \pi^! \eps_\bn^!\Gc \,=\, \rho^! p_{(k_1,\cdots, k_l)}^!\Gc \,=\,
\rho^! p_{(k_1,\cdots, k_l)}^* \Gc \,=\, \rho^!\Gc' 
 \]
 has $\ul H^{\neq -r}=0$. But this folllows from Proposition \ref{prop:7.15}  once we notice that $F$ is a disjoint
 union of tube domains $\RR^l+ iC$ for some faces $C$ of $\Hc$ of dimension $r$. Indeed, let $(z_1,\cdots, z_l)\in
 \CC^l_{(k_1,\cdots, k_l)}$ with $z_\nu = x_\nu + y_\nu  \sqrt{-1}$, $x_\nu, y_\nu\in\RR$. The condition that
 $p_{(k_1,\cdots, k_l)}(z_1,\cdots, z_l) = \sum\xi_{k_i} \cdot z_i$ lies in $Z_\bn^\Im$ is a condition
 on $(y_1,\cdots, y_l)\in\RR^l$. More precisely, there should be exactly $r$ distinct values among the
 $y_1,\cdots, y_l$, so $\{1,\cdots, l\}$ splits into a disjoint union of $r$ nonempty subsets $I_1, \cdots, I_r$
 with $y_i=y_j$ for any $i,j\in I_\nu$. Denoting this common value $y_{I_\nu}$ and numbering the $I_\nu$
 by $y_{I_1} <\cdots < y_{I_r}$, the final condition for $p_{(k_1,\cdots, k_l)}(z_1,\cdots, z_l)\in Z_\bn^\Im$
 is 
 \be\label{eq:sum-xi-k-i}
 \sum_{i\in I_1} \xi_{k_i} = n_1, \cdots, \sum_{i\in I_r} \xi_{k_i} = n_r.
 \ee
 Now, faces of $\Hc$ of dimension $r$ are precisely labelled by ordered sequences $\bI=(I_1,\cdots, i_r)$
 of nonempty subsets of $\{1,\cdots, l\}$, forming a disjoint decomposition. The face $C_\bI$
 associated to $I$ consists of $(y_1,\cdots, y_l)$ with the $y_i$ for $i$ in each $I_\nu$ being equal to some 
 number $y_{I_\nu}$ such that $y_{I_1} <\cdots < y_{I_r}$.  Thus $F$ is the union of  the  tube domains
 $\RR^l+  C_\bI \sqrt{-1}$
 for $\bI$ satisfying \eqref{eq:sum-xi-k-i}. Proposition \ref{prop:eps_n^!} is proved. 
 
 \vskip .2cm
 
 We now construct a canonical Cousin resolution of any $\Gc\in\Perv^{<\oo}(Z, \Vc)$, $Z= Z(\CC,L)_n$ by
 the method of \cite{KS-arr}. That is, we have a filtration of $Z$ by closed subsets
 \[
 F_1\subset F_2\subset \cdots \subset Z, \quad F_r \,=\bigsqcup_{l(\bn)\leq r} Z_\bn^\Im, 
 \]
 with $\bn$ running over ordered partitions of $n$. This filtration may be infinite, but $\Gc$ being supported on a
 finite skeleton of $\Sc^{(2)}$, is supported on some finite $F_r$. Denote $\eps_r: F_r\hra Z$ the embedding.
 Then we have a Postnikov system for $\Gc$: the sequence
 \be\label{eq:post-sys}
 0\buildrel \eta_0\over\lra \eps_{1*} \eps_1^! \Gc \buildrel\eta_1\over\lra \eps_{2*}\eps_2^!\Gc \buildrel\eta_2\over\lra
 \cdots
 \ee
 stabilizes to $\Gc$ (so $\eta_r =\Id_\Gc$ for $r\gg 0$) and
 \be\label{eq:cone-eta-r}
 \Cone (\eta_r) \,\= \bigoplus_{l(\bn)=r+1} \eps_{bn *} \eps_\bn^!\Gc \,\= \bigoplus_{l(\bn)=r+1} \eps_{\bn *} \Fc_\bn [-r], 
 \ee
 where $\Fc_\bn = \ul H^{-r-1}(\eps_\bn^!\Gc)$, the last identification in \eqref{eq:cone-eta-r} following from
 Proposition  \ref{prop:eps_n^!}. In fact, if $\Gc$ is represented as a complex if injective sheaves, then 
 \eqref{eq:post-sys} is a filtration on $\Gc$, eventually stabilizing to $\Gc$ and with quotients given by
 \eqref{eq:cone-eta-r}. This gives a complex of sheaves
 \[
 \begin{gathered}
 \Fc^\bullet(\Gc) \,:= \, \bigl\{ \cdots \to \Cone(\eta_2)[-2] \buildrel d\over\lra \Cone(\eta_1)[-1]\buildrel d\over\lra \Cone(\eta_0)\bigr\} \,\=
 \\
 \= \, \bigl\{ \cdots \buildrel d\over\lra \bigoplus_{l(\bn)=3}\eps_{\bn *} \Fc_\bn \buildrel d\over\lra 
  \bigoplus_{l(\bn)=2}\eps_{\bn *} \Fc_\bn \buildrel d\over\lra \eps_{(n) *} \Fc_{(n)} \bigr\} 
 \end{gathered}
 \]
 whose differential $d$ is ``the differential $d_1$ in the spectral sequence of the filtered complex 
 \eqref{eq:post-sys}'' and which is quasi-isomorphic to $\Gc$. In other words, we constructed a canonical
 resolution $\Fc^\bullet(\Gc)$ of $\Gc$ which is a Cousin-type complex. 
 
 
 Denote by $\EE(\Gc)$  the Janus sheaf associatedto the Cousin-type complex $\Fc^\bullet(\Gc)$ by
 Proposition \ref{prop:Cousin=Janus}.
 
 \begin{prop}\label{prop:EE(F(E))=E}
 Let $E$ be a Janus sheaf and  $\Fc^\bullet(E)$ be the perverse sheaf associated to $E$ by Proposition 
 \ref{prop:cous-perv}. Then $\EE(\Fc^\bullet(E))\= E$. 
 \end{prop}
 
 \noindent{\sl Proof:} We show the identification on the level of objects of $\ul\CM_n(L)$. The comparison of
 values on morphisms is similarly straightforward. 
 
 The value $\EE(\Fc(E))(M)$, $M\in\CM_n(L)$ is the stalk at $U_M$ of the sheaf $\eps_\bn^!\Fc^\bullet(E)$. Now,
 the terms of the complex $\Fc^\bullet(E)$ have the form $\eps_{\bm *} \Fc_\bm$, where $\Fc_\bm$ is the
 sheaf on $Z_\bm^\Im$ associated to the $\delta'_f$ data of $E$. So our statement will follow from the following
 identification:
 \[
 \eps_{\bn}^! \eps_{\bm *} \Fc_\bm \,=\begin{cases}
 \Fc_\bm, & \text{ if } \bn=\bm;
 \\
 0,& \text{ otherwise}. 
 \end{cases}
 \]
 To see this last identification, consider the Cartesian diagram of closed embeddings
 \[
 \xymatrix{
 I \ar[r]^{\zeta_\bm}
 \ar[d]_{\zeta_\bn}
  & Z_\bm^\Im \ar[d]^{\eps_\bm}
 \\
 Z_\bn^\Im \ar[r]_{\eps_\bn} & Z,
 } \quad I:= Z_\bn^\Im\cap Z_\bm^\Im. 
 \]
 By base change, $\eps_\bn^! \eps_{\bm *} \Fc_m \= \zeta_{\bn *} \zeta_\bm^! \Fc_\bm$. If $\bm\neq\bn$, then $I=\emptyset$ and we get $0$. If $\bm=\bn$, then $I=Z_\bn^\Im=Z_\bm^\Im$ and we get $\Fc_\bm$. \qed
 
 \begin{prop}
 If $\Gc\in\Perv^{<\oo}(Z,\Vc)$, then the Janus presheaf $\EE(\Gc)$ is a Janus sheaf.
 \end{prop}
 
 \noindent {\sl Proof:} The stratification induced by $\Sc^{(0}$ on the imaginary strata is the restriction of $\Sc^{(1)}$,
 so each $\Fc_\bn$ is $\Sc^{(1)}$-constructible. This means that anodyne morphisms from $\Gro'(\ul\CM_n(L))$ are taken by $\EE(\Gc)$ into isomorphisms. This is one hanf of the conditions (JS3). As for the other half -- that
 anodyne morphisms in $\Gro''(\ul\CM_n(L))$ are taken by $\EE(\Gc)$ into isomorphisms, it follows from the first
 half applied to  the perverse sheaf $\tau^*\DD(\Gc)$. Indeed, Propositions \ref{prop:EE(F(E))=E} and 
 \ref{prop:D=tau} imply that $\EE(\tau^*\DD(\Gc))\= \EE(\Gc)^\tau$. So $\EE(\tau^*\DD(\Gc))$ taking anodyne
 morphisms in $\Gro'(\ul\CM_n(L))$ into isomorphisms implies that so does $\EE(\Gc)^\tau$, i.e.,
 that $\EE(\Gc)$ takes anodyne morphisms in $\Gro''(\ul\CM_n(L))$ into isomorphisms. \qed. 
 
 \vskip .2cm
 
 The above proposition means that we have a functor
 \[
 e: \Perv^{<\oo}(Z, \Vc) \lra \JS^{<\oo}_{L,n}(Z, \Vc), \quad \Gc\mapsto \EE(\Gc). 
 \]
 We claim that $e$ is quasi-inverse to the functor $c$ from \eqref {eq:functor-c}. Indeed, Proposition 
 \ref{prop:EE(F(E))=E} means that $ec\=\Id$. The isomorphism $ce\=\Id$ follows directly from the construction of 
 $\EE(\Gc)$, as $e(\Gc)$ is the Janus presheaf (in fact sheaf) associated to the Cousin resolution $\Fc^\bullet(\Gc)$,
 so $c(e(\Gc))=\Fc^\bullet(\Gc)$ is quasi-isomorphic to $\Gc$. 
 
 \vskip .2cm
 
 This finishes the proof of Theorem \ref{thm:perv=cont}.

 
   \section{The universal Janus sheaf}\label{sec:uni-janus}
   
   In this section we prove Theorem \ref{thm:Bn-perv}. The argument, reducing it to Theorem  \ref{thm:perv=cont},
   proceeds by the same route as used in \cite{KS-prob} for the case $L=\ZZ_+$. So we indicate the main
   steps but will be brief on details which are similar to \cite{KS-prob}. 
   
   \paragraph{The $\k$-linear category of contingency matrices.} Let $n\in L$. We introduce a $\k$-linear 
   category $\CMen_n(L)$ whose objects are formal symbols $[M]$, $M\in\CM_n(L)$. Morphisms in $\CMen_n(L)$
   are generated by the generating morphisms
   \begin{itemize}
   \item $\delta'_f: [M]\to [N]$ for any morphism $f: N\to M$ in $\Gro'(\ul\CM_n(L))$,
   
   \item $\delta''_g: [N]\to [M]$ for any morphism $g: N\to M$ in $\Gro''(\ul\CM_n(L))$,
   \end{itemize}
subject to the following relations:
\begin{itemize}
\item[($\CMen 1$)] If $f$ is an identity morphism, then $\delta'_f=\Id$. Further, 
$\delta'_{f_1f_2} = \delta'_{f_2} \delta'_{f_1}$ for any two composable morphisms $f_1, f_2$ in $\Gro'(\ul\CM_n(L))$.
Similarly, if $g$ is an identity morphism, then $\delta''_g=\Id$. Further, $\delta''_{g_1g_2}=\delta''_{g_1}\delta''_{g_2}$
for any two composable morphisms $g_1, g_2$ in  $\Gro''(\ul\CM_n(L))$.

\item[($\CMen 2$)] For any mixed fork
$P=\{M'\buildrel \phi'\over\to N'\buildrel \psi''\over\leftarrow N\}$  in  $\Gro(\ul\CM_n(L))$,
 we have the relation
\[
\delta'_{\phi': M'\to N} \delta''_{\psi'': N\to N'} \,=\sum_{M\in\Sup(P) }\delta''_{\psi'': M\to M'} \delta'_{\phi': M\to N}.
\]
where the sum os over $M$ fitting into a diagram \eqref{eq:sup-diag}. 

\item[($\CMen 3$)] If $f: N\to M$ is an anodyne morphism in $\Gro'(\ul\CM_n(L))$, then $\delta'_f$
is an isomorphism. If $g: N\to M$ is is an anodyne morphism in $\Gro''(\ul\CM_n(L))$, then $\delta''_g$
is an isomorphism.
\end{itemize}

 \noindent The meaning of ($\CMen 3$) is that $\CMen_n(L)$ is defined as the localization $\wt\CMen_n(L)[\Ano^{-1}]$,
  where $\wt\CMen_n(L)$ is the $k$-linear category with relations  ($\CMen 1-2$) and $\Ano$ is the set of
  morphisms in  $\wt\CMen_n(L)$     consisting of $\delta'_f$ and $\delta''_g$ for $f$ or $g$ anodyne. 

\paragraph{The universal Janus sheaf.} The relations in $\CMen_n(L)$  mimick the axioms $(JS1\text{-}3)$ of a Janus
sheaf. This gives the following fact. 

\begin{prop}\label{prop:E-un}
(a) The correspondence
\[
E^\un: M\mapsto [M], \,\, f\mapsto \delta'_f, \,\, g\mapsto \delta''_g
\]
is a Janus sheaf of type $(L,n)$ with values in $\CMen_n(L)$.

\vskip .2cm

(b) Let $\Vc$ be any $\k$-linear category and $E$ be a Janus sheaf of type $(L,n)$ with values in $\Vc$.
Then there is a unique functor $F_E: \CMen_n(L)\to \Vc$ such that $F(E^\un) = E$. \qed
\end{prop}

\paragraph{The monoidal structure.} Let $\CMen(L)=\bigoplus_{n\in L}\CMen_n(L)$ be the
orthogonal direct sum of the $\CMen_n(L)$. Explicitly
\[
\Ob(\CMen(L)) \,=\, \bigsqcup_{n\in L} \Ob(\CMen_n(L)) \,=\,\bigl\{ [M], \, M\in \CM(L)\bigr\}. 
\]
 If $M\in \CM_n(L)$ and $N\in\CM_n(L)$, then
 \[
 \Hom_{\CMen(L)}([M], [N]) \,= \begin{cases}
 0, & \text{ if } m\neq n, 
 \\
 \Hom_{\CMen_n(L)} ([M], [N]), & \text{ if } m=n. 
 \end{cases}
 \]
We introduce a monoidal structure $\otimes$ on $\CMen(L)$ to be given on objects by
\[
[M]\otimes [N] \,=\, [M\oplus N], \quad M\oplus N = \begin{pmatrix}
M&0 \\ 0&N
\end{pmatrix}. 
\]
Further, let $f_1: M_1\to N_1$ and $f_2: M_2\to N_2$ be two morphisms in
$\Gro'(\ul\CM_{n_1}(L))$ and $\Gro'(\ul\CM_{n_2}(L))$ respectively. Then we have a natural
morphism $f_1\oplus f_2: M_1\oplus M_2\to N_1\oplus N_2$ in $\Gro'(\ul\CM_{n_1+n_2}(L))$. 
More precisely, $f_1$ expresses the fact that $N_1$ is obtained from $M_1$ by summing some sets $J_1,\cdots, J_p$ of adjacent columns, while $f_2$ expresses the fact that $N_2$ is obtained from $M_2$ by summing
some sets $K_1,\cdots, K_q$ of djacent columns. Now, the set of columns of $M_1\oplus M_2$ is identified
with the disjoint union of the sets of columns of $M_1$ and $M_2$. With respect to this identification,
after summing the sets of columns $J_1,\cdots, J_p, K_1,\cdots, K_q$ in $M_1\oplus M_2$, we get $N_1\oplus N_2$
and this gives $f_1\oplus f_2$.

\vskip .2cm

Similarly, if $g_1: M_1\to N_1$ and $g_2: M_2\to N_2$ are two morphisms in
$\Gro''(\ul\CM_{n_1}(L))$ and $\Gro''(\ul\CM_{n_2}(L))$ respectively, we have a natural morphism
$g_1\oplus g_2: M_1\oplus M_2\to N_1\oplus N_2$ in $\Gro''(\ul\CM_{n_1+n_2}(L))$.

\vskip .2cm

We now define the monoidal structure $\otimes$ on generating morphisms as follows:
\begin{itemize}

\item If $f_1, f_2$ are two morphisms in $\Gro'(\ul\CM(L))$, then we put
 $\delta'_{f_1}\otimes\delta'_{f_2}=\delta'_{f_1\oplus f_2}$.
 
 \item  If $g_1, g_2$ are two morphisms in $\Gro''(\ul\CM(L))$, then we put
 $\delta''_{g_1}\otimes\delta''_{g_2}=\delta''_{g_1\oplus g_2}$.
 
 \item Let $f: M_1\to N_1$ be a morphism in $\Gro'(\ul\CM(L))$ and $g: M_2\to N_2$
  be a morphism in $\Gro''(\ul\CM(L))$. We have a commutative diagram in   $\Gro(\ul\CM(L))$:
  \[
  \xymatrix{
  M_1\oplus M_2 \ar[rr]^{f\oplus \Id_{M_2}}_{'} 
  \ar[d]_{\Id_{M_1}\oplus g} ^{''}
  && N_1\oplus M_2
  \ar[d]^{\Id_{N_1}\oplus g} _{''}
  \\
  M_1\oplus N_2 \ar[rr]^{'}_{f\oplus \Id_{N_2}} && N_1\oplus N_2
  }
  \]
which is the unique diagram \eqref{eq:sup-diag} containing the mixed form obtained by omitting $M_1\oplus M_2$.
Therefore by ($\CMen 2$)
\[
\delta''_{\Id_{M_1}\oplus g} \delta'_{f\oplus \Id_{M_2}} \,=\,\delta'_{f\oplus \Id_{M_1}} \delta''_{\Id_{N_1}\oplus g}
\]
and we define $\delta'_f\otimes \delta''_g$ to be equal to this common value.

\item The definition of $\delta''_g\otimes\delta'_f$ is similar. 
\end{itemize}

\begin{prop}
The above rules for the tensor products of generating morphisms extend to a monoidal structure on the category
$\CMen(L)$ with the unit object $[\emptyset]$, where $\emptyset$ is the empty (size $0\times 0$) contingency matrix.
\end{prop}

\noindent{\sl Proof:} Similar to \cite[Prop.3.5]{KS-prob}. \qed

\paragraph{The braiding.} We now define a braiding on the monoidal category $(\CMen(L), \otimes)$. The construction is 
an adaptation of that of \cite[\S3B-C]{KS-prob}. 

\vskip .2cm

Let $M=\| m_{ij}\|_{i=1\cdots, r}^{j=1,\cdots, s}\in\CM(L)$.
 We denote its rows and columns by
\[
M_{i,\bullet} = (m_{i1},\cdots, m_{is}), \quad M_{\bullet j} = (m_{1j},\cdots, m_{rj})^t.
\]
Denote by $\sigma''_{i, i+1}M$ the matrix obtained from $M$ by interchanging the $i$ and $(i+1)$st rows, i.e.,  $M_{i,\bullet}$ and $M_{i+1,\bullet}$. Similarly, let $\sigma'_{j, j+1}M$ be the matrix obtained from $M$ by interchanging 
the $j$th and $(j+1)$st columns, i.e., $M_{\bullet,j}$ and $M_{\bullet, j+1}$. 

\vskip .2cm

Suppose that $M_{i,\bullet}$ and $M_{i+1,\bullet}$ are {\em disjoint}, i.e., for each $j$ at least one among
$\{m_{ij}, m_{i+1, j}\}$ is equal to $0$. then we have anodyne morphisms in $\Gro''(\ul\CM(L))$:
\be\label{eq:exch-ano}
M\lra \del''_{i-1}(M) =\del''_{i-1}(\sigma''_{i, i+1} M) \lla \sigma''_{i, i+1}(M)
\ee
and we define the {\em row exhange isomorphism} 
$\eps''_{i, i+1}: [M] \to [\sigma''_{i, i+1} M]$
in $\CMen(L)$ as the composition of the $\delta''$-morphism and the inverse of another such corresponding to
\eqref{eq:exch-ano}. 

\vskip .2cm
 
Similarly, if $M_{\bullet, j}$ and $M_{\bullet, j+1}$ are disjoint, i.e., for each $i$ at least one among
$\{m_{ij}, m_{i, j+1}\}$ is equal to $0$, we have anodyne morphisms in $\Gro''(\ul\CM(L))$
\be\label{eq:each-ano-2}
M\lra \del'_{j-1}(M) = \del'_{j-1}(\sigma'_{j, j+1}M) \lla \sigma'_{j, j+1} M
\ee
and we define the {\em column exchange isomorphism} $\eps'_{j, j+1}: [M] \to [\sigma'_{j, j+1} M]$ in
$\CMen(L)$ by composing the $\delta'$-morphism and the inverse of another such corresponding to 
\eqref{eq:each-ano-2}.

\vskip .2cm

Let now $M,N\in\CM(L)$ be two contingency matrices of sizes $p\times q$ and $r\times s$ respectively. 
We define the braiding isomorphism
\[
R_{[M], [N]}:  [M] \x [N] \,=\, [M\+ N] \lra [N\+ M] \,=\, [N]\x [M]
\]
 as the composition
 \[
 [M\oplus N] =
 \left[\begin{pmatrix} 
 M&0\\ 0&N
\end{pmatrix}\right]
 \buildrel R''_{M,N}
 \over\lra
 \left[  \begin{pmatrix} 
 0&N\\ M&0
\end{pmatrix}\right] \buildrel R'_{M,N}\over\lra 
\left[
 \begin{pmatrix} N&0\\ 0&M
\end{pmatrix}\right]
=[N\oplus M], 
 \]
 where: 
 \begin{itemize}
 \item[($B''$)]  $R''_{M,N}$ is the composition of $pr$ row exchange isomorphisms
 moving $r$ rows of $(0 N)$ past the $p$ rows of $(M 0)$.

 \item[($B'$)]  $R'_{M,N}$ is the composition of $qs$ column exchange isomorphisms moving $s$ columns of
 $\begin{pmatrix} N\\0\end{pmatrix}$ past $q$ columns of 
  $\begin{pmatrix} 0\\M\end{pmatrix}$. 
 \end{itemize}
 
 Note that there are  several ways to move the two groups of rows in ($B''$) or of columns in ($B'$) past each other but
 the result is independent of the choice, similar to \cite[Prop.3.4]{KS-prob}.

\begin{prop}\label{prop:CM-braided}
The isomorphisms $R_{[M], [N]}$ define a braiding on the monoidal category $\CMen(L)$. 
\end{prop}

\noindent{\sl Proof:} Similar to \cite[Prop.3.6]{KS-prob}. \qed

\paragraph{A graded bialgebra in $\CMen(L)$.} We now define
a graded bialgebra $\aen=(\aen_n)_{n\geq 0}$ in the braided category $\CMen(L)$
with components $\aen_n=[n]$ (the object corresponding to the $1\times 1$
contingency matrix $(n)$) for $n>0$ and $\aen_0=[\emptyset]$. 
The multiplication and comultiplication are given by
\[
\begin{gathered}
\mu_{m,n}: [m]\otimes [n] \,=\left[\begin{pmatrix} m&0\\0&n
\end{pmatrix} \right] \buildrel (\delta')^{-1}\over\lra  \left[
\begin{pmatrix} m\\ n\end {pmatrix}\right] \buildrel \delta''\over\lra
[m+n],
\\
\Delta_{m,n}: [m+n] \buildrel \delta'\over\lra [(m, \, n)] 
\buildrel (\delta'')^{-1}\over\lra \left[\begin{pmatrix} m&0\\0&n
\end{pmatrix} \right] \,=\, [m]\otimes [n]. 
\end{gathered}
\]

\begin{prop}\label{prop:aen-in-CM}
The morphisms $\mu_{m,n}$, $\Delta_{m,n}$ make $\aen$ into a graded
bialgebra in $\CMen$. 
\end{prop}

\noindent{\sl Proof:} Similar to \cite[Prop.4.1]{KS-prob}. \qed

\paragraph{The category $\CMen(L)$ and the PROB $\Ben^L$.} Because $\CMen(L)$ is the target of the 
universal Janus sheaf of type $(L,n)$, to deduce Theorem 
 \ref{thm:Bn-perv} from Theorem  \ref{thm:perv=cont}, it suffices to show the following.
 
 \begin{thm}\label{thm:Janus=CM}
 We have an equivalence of braided monoidal categories $\xi: \CMen(L)\to\Ben^L$ which for any $n\in L$ induces
 an eqiuvalence of categories $\xi_n: \CMen_n(L)\to\Ben^L_n$. 
 \end{thm}
 
 \noindent {\sl Proof:} The graded bialgebra $\aen$ in $\CMen(L)$ gives, by the universal property of $\Ben^L$,
 a braided monoidal functor $F_\aen: \Ben^L\to\CMen(L)$. On the other hand, consider the universal graded
 bialgebra $\ba = (\ba_n)_{n\in L}$ in $\Ben^L$. By Proposition \ref{prop:bialg-Janus} $\ba$ gives, for each
 $n\in L$, a Janus sheaf $E_{\ba, n}$ of type $(L,n)$ with values in $\Ben^L$. By Proposition \ref{prop:E-un},
 $E_{\ab, n}$ is unduced from the universal Janus sheaf $E^\un$ by a unique (up to a unique isomorphism)
 functor $\xi_{\ba, n}: \CMen_n(L)\to \Ben^L$. It is straightforward to check that the $\xi_{\ba, n}$ for
 various $n\in L$ unite into a monoidal functor $\xi_\ba: \CMen(L)\to\Ben^L$ which is quasi-inverse to $F_\aen$.
 \qed
 
 \vskip .2cm
 
 This finishes the proof of Theorem \ref{thm:Bn-perv}.



\vskip 1cm

\small{

M.K.: 
Kavli IPMU (WPI), UTIAS, University of Tokyo, 5-1-5 Kashiwanoha, Kashiwa, Chiba, 277-8583 Japan.
 Email: 
{\tt mikhail.kapranov@protonmail.com}

\smallskip

V.S.: Institut de Math\'ematiques de Toulouse, Universit\'e Paul Sabatier, 118 route de Narbonne, 
31062 Toulouse, France; Kavli IPMU  (WPI), UTIAS, University of Tokyo, 5-1-5 Kashiwanoha, 
Kashiwa, Chiba, 277-8583 Japan. Email: 
{\tt schechtman@math.ups-tlse.fr }

}

\ed